\documentclass[11pt]{article}
\usepackage[T1]{fontenc}
\usepackage{lmodern}
\usepackage[a4paper,margin=28mm]{geometry}
\usepackage{amsmath,amssymb,amsthm,mathtools,mathrsfs,bm}
\usepackage{booktabs,array,longtable,tabularx,enumitem}
\usepackage{tikz,xcolor,graphicx,caption,needspace}
\usetikzlibrary{calc,positioning,arrows.meta,patterns,matrix,decorations.pathreplacing}
\usepackage[authoryear,round]{natbib}
\usepackage[unicode,hidelinks]{hyperref}
\newcommand{\R}{\mathbb R}

\newcommand{\A}{\mathcal A}

\newcommand{\coker}{\operatorname{coker}}
\newcommand{\rank}{\operatorname{rank}}
\newcommand{\sgn}{\operatorname{sgn}}

\newcommand{\diag}{\operatorname{diag}}
\newcommand{\Span}{\operatorname{span}}
\newcommand{\eps}{\varepsilon}
\newcommand{\pd}{\partial}
\newcommand{\dd}{\mathrm d}
\newcommand{\pos}[1]{\bigl(#1\bigr)_+}

\numberwithin{equation}{section}
\allowdisplaybreaks[2]
\setlist{itemsep=3pt,topsep=4pt,leftmargin=2em}
\newtheorem{theorem}{Theorem}[section]
\newtheorem{lemma}[theorem]{Lemma}
\newtheorem{proposition}[theorem]{Proposition}
\newtheorem{corollary}[theorem]{Corollary}
\theoremstyle{definition}

\newtheorem{example}[theorem]{Example}
\theoremstyle{remark}
\newtheorem{remark}[theorem]{Remark}
\title{Universal Dimension Formula\\for Planar $C^r$ Splines at Degree $3r+1$}
\author{Xinyu Wu\\{\small\texttt{xinyu97@ustc.edu.cn}}\and Jiansong Deng\\{\small\texttt{dengjs@ustc.edu.cn}}}
\date{}
\hypersetup{pdftitle={Universal Planar Spline Dimension at Degree 3r+1},pdfauthor={Xinyu Wu and Jiansong Deng}}
\begin{document}
\maketitle
\begin{abstract}
We prove that Schumaker's lower bound gives the dimension of $S_{3r+1}^r(\Delta)$ for every nondegenerate triangulation $\Delta$ of a polygonal disk in $\R^2$ and every integer $r\ge3$. After accounting for the local vertex contributions, we reduce the remaining dimension defect to compatible vertex functionals. Apolarity associates a unique polynomial with each face, and compatibility forces differences of adjacent face polynomials to contain powers of both endpoint forms. We eliminate the leading residual coefficient using Hermite relations when $r$ is odd and an additional coefficient comparison when $r$ is even. A uniform induction then raises the endpoint powers until all face differences vanish; the zero boundary face potentials imply that every vertex functional vanishes and the global defect is zero. Together with the known cases $r=1,2$ and Hong's higher-degree theorem, this yields the dimension formula for all $r\ge1$ and $d\ge3r+1$.
\end{abstract}
\section{Introduction}\label{sec:intro}

Let $\Delta$ be a finite triangulation of a planar polygonal disk $\Omega$. Every triangle is nondegenerate, and the intersection of two distinct triangles is either empty, a common vertex, or an entire common edge. For nonnegative integers $d,r$, define
\[
 S_d^r(\Delta):=\{s\in C^r(\Omega):s|_\tau\in\mathbb P_d(\tau)
                         \text{ for every }\tau\in\Delta_2\}.
\]
Here $\Omega=|\Delta|$, and $\mathbb P_d(\tau)$ is the space of restrictions to $\tau$ of bivariate polynomials on the plane of total degree at most $d$. No boundary conditions are imposed. The dimension depends on the incidence structure and on the supporting directions of the edges \citep{LS07}, and Schumaker's lower bound incorporates the local dependencies among smoothness conditions at each interior vertex \citep{Schumaker79}.

The dimension formula is known for every triangulation when $d\ge4r+1$ \citep{AS87}. Hong established attainment of Schumaker's lower bound on arbitrary polygonal triangulations for $d\ge3r+2$ \citep{Hong91}. At the adjacent degree $d=3r+1$, \citet{AS90} established the generic dimension formula. Local direction counts determine the number of independent relations among the equations at one vertex; the remaining problem is to exclude additional relations involving equations at different vertices on the given geometric realization. We prove attainment at degree $3r+1$ for every nondegenerate geometric realization and each integer $r\ge3$. Combined with the known cases $r=1,2$, this gives universal attainment at degree $3r+1$ for every $r\ge1$.

To state the dimension formula, write $E_I=|\Delta_1^\circ|$ and $V_I=|\Delta_0^\circ|$, where an interior edge is an edge not contained in the boundary. For an interior vertex $v$, let $t_v$ be the number of distinct unoriented supporting lines of its incident edges, and put
\begin{equation}\label{eq:sigma}
 \sigma_v(r)=\sum_{j\ge1}\pos{r+1-j(t_v-1)},\qquad
 \pos{a}=\max\{a,0\}.
\end{equation}
Opposite rays count as one supporting line. An interior vertex has $t_v\ge2$, so the sum is finite.
\begin{theorem}[Universal attainment at degree $3r+1$]\label{thm:main}
For every such triangulation and every integer $r\ge3$,
\begin{equation}\label{eq:main}
 \dim S_{3r+1}^r(\Delta)
 =\binom{3r+3}{2}+\binom{2r+2}{2}(E_I-2V_I)
       +\sum_{v\in\Delta_0^\circ}\sigma_v(r).
\end{equation}
\end{theorem}
For $r\ge1$ and $d\ge3r+1$, write the Schumaker lower bound with the supporting directions of the given realization as
\begin{equation}\label{eq:lower-bound-general}
 \begin{aligned}
 L_d^r(\Delta)={}&\binom{d+2}{2}+E_I\binom{d-r+1}{2}\\
 &-V_I\left[\binom{d+2}{2}-\binom{r+2}{2}\right]
   +\sum_{v\in\Delta_0^\circ}\sigma_v(r).
 \end{aligned}
\end{equation}
At $d=3r+1$ this is the right side of \eqref{eq:main}.
\begin{corollary}[Universal attainment in the sufficient range]\label{cor:all-orders}
For every triangulation considered here and all integers $r\ge1$, $d\ge3r+1$,
\begin{equation}\label{eq:unified-range}
 \dim S_d^r(\Delta)=L_d^r(\Delta).
\end{equation}
\end{corollary}
For $d=3r+1$, the claim follows from Theorem~\ref{thm:main} when $r\ge3$, from \citet*{APS87} when $r=1$, and from \citet{WD26} when $r=2$. The case $d\ge3r+2$ follows from \citet{Hong91}.

Our proof computes the local vertex dimensions by matrices and binary apolarity, then represents a global row relation by vertex functionals and common face polynomials. The main difficulty is to rule out a nonzero family of leading coefficients satisfying all edge compatibility conditions. For odd $r\ge3$, Hermite endpoint relations suffice. Even orders $r\ge4$ require a second coefficient comparison. It yields a separation property that contradicts a balance identity inherited from the first comparison. This contradiction eliminates the nonzero family of leading coefficients. A uniform induction then raises the endpoint powers until all face differences vanish; the zero boundary face potentials eliminate the remaining vertex functional classes and give zero global defect.

Attainment of the dimension bound and the existence of bases with prescribed local supports are separate questions. Alfeld and Schumaker proved that, for every $r\ge1$ and $d\le3r+1$, there exist triangulations admitting no basis whose members are each supported in a single vertex star \citep{AS00}. Building on this dimension formula, the next step is to develop a uniform construction of bases for $S_{3r+1}^r(\Delta)$, valid for every $r\ge1$ and every nondegenerate triangulation $\Delta$.

Furthermore, at the adjacent degree $d=3r$, the edge traces no longer determine a unique face polynomial, so the present construction cannot be applied directly; see Appendix~\ref{app:degree}. The question posed by \citet{WD26} of whether $S_6^2$ attains Schumaker's lower bound on every nondegenerate geometric realization remains unresolved.

The rest of the paper is organized as follows. Section~\ref{body:framework} computes the local vertex dimensions and expresses the global problem successively as a matrix left kernel, vertex functionals, and common face polynomials. Section~\ref{body:odd} proves the dimension formula for odd $r\ge3$, and Section~\ref{body:even} completes the argument for even $r\ge4$. Section~\ref{body:low} discusses the framework at $r=1,2$, and Section~\ref{body:conclusion} summarizes the results and the scope of the method. The appendices provide the detailed proofs and coefficient calculations.

\section{From local vertex relations to common face polynomials}\label{body:framework}

We first count the dependencies internal to one vertex, and then ask whether putting all vertex relations together introduces further dependencies. We express the latter question successively as a matrix left kernel, a family of vertex functionals, and polynomial differential conditions. This section records the results and their logical connections; Appendix~\ref{app:framework} gives the complete derivations.

\subsection{Local dimensions from distinct supporting lines}
Set $R=\R[x,y,z]$, and let $R_d$ be its homogeneous degree-$d$ part, identified with bivariate polynomials of degree at most $d$. Choose a homogeneous supporting-line equation $\ell_e$ for each interior edge, and put
\[
 J(e)=(\ell_e^{r+1}),\qquad J(v)=\sum_{e\ni v}J(e).
\]
The $C^r$ condition across $e$ is exactly $s_{\tau_+}-s_{\tau_-}\in J(e)_d$. Thus $J(v)_d$ is the span of the allowed edge differences in the closure equation at $v$. Repeated supporting lines give the same subspace; opposite rays count as one line \citep{Schumaker79,SS97,MV13}.
\begin{proposition}[Local vertex dimensions]\label{body:localdim}
Suppose there are $t=t_v$ distinct supporting lines at $v$. For every $d\ge r$,
\begin{equation}\label{body:localdim-general}
 \begin{aligned}
 \dim J(v)_d
 &=\sum_{k=r+1}^{d}\min\{k+1,t(k-r)\}\\
 &=\binom{d+2}{2}-\binom{r+2}{2}
   -\sum_{s=1}^{d-r}\pos{r+1-(t-1)s}.
 \end{aligned}
\end{equation}
In particular, at $d=3r+1$,
\[
 \dim(R/J(v))_d=\binom{r+2}{2}+\sigma_v(r),\qquad
 \dim J(v)_d=2\binom{2r+2}{2}-\sigma_v(r).
\]
\end{proposition}
After translating the vertex to the origin, split by homogeneous binary degree $k$. The coefficient matrix on each piece has rank $\min\{k+1,t(k-r)\}$. The appendix proves this both by a confluent Vandermonde matrix and by binary apolarity. The first proof computes the coefficient rank directly; the second also converts multiplication by line powers into differentiation, which will be used in the global problem.

\subsection{The global question is a matrix left kernel}
Orient every interior edge. For $g_e\in J(e)_d$, place $+g_e$ at the head and $-g_e$ at the tail, keeping only interior vertex components. This gives the endpoint map
\[
 \partial_J:\bigoplus_{e\in\Delta_1^\circ}J(e)_d
 \longrightarrow\bigoplus_{v\in\Delta_0^\circ}J(v)_d,
 \qquad (g_e)_e\longmapsto
 \left(\sum_{e\ni v}\varepsilon_{v,e}g_e\right)_v.
\]
Choose independent bases in the local spaces, write $T_d$ for its matrix, and put $M_d=\coker\partial_J$.
\begin{proposition}[The remaining dimension term]\label{body:dimension}
For a triangulated disk and $d\ge r+1$,
\[
 \dim S_d^r(\Delta)=\binom{d+2}{2}+E_I\binom{d-r+1}{2}
 -\sum_{v\in\Delta_0^\circ}\dim J(v)_d+\dim M_d.
\]
\end{proposition}
This follows by integrating compatible edge differences over the disk \citep{Billera88}. The local dimensions have already counted the dependencies at each single vertex. The main theorem is therefore equivalent to full row rank of $T_{3r+1}$: no new row relation appears when all vertex equations are assembled.

If $a^{\mathsf T}T_d=0$, partition $a$ into vertex blocks to obtain a functional $\varphi_v$ on each $J(v)_d$. This definition allows the local dimensions to vary with $v$. Varying only the input on $e=uv$ gives
\[
 \varphi_u|_{J(e)_d}=\varphi_v|_{J(e)_d}.
\]
Assign the zero functional to every boundary vertex. The task is to deduce $\varphi_v=0$ for all $v$, hence $a=0$, from these edge conditions.

\subsection{Common edge traces and zero-trace spaces}
Let $A=\R[X,Y,Z]$, with the perfect pairing
\[
 \langle f,H\rangle=f(\partial_X,\partial_Y,\partial_Z)H\big|_0,
 \qquad f\in R_d,\ H\in A_d.
\]
For $\ell=ax+by+cz$, write $D_\ell=a\partial_X+b\partial_Y+c\partial_Z$. The adjoint identity
\[
 \langle\ell^{r+1}q,H\rangle=\langle q,D_\ell^{r+1}H\rangle
\]
represents every vertex functional by some $H_v\in A_d$. Define
\[
 \mathcal Z_e=\ker(D_{\ell_e}^{r+1}:A_d\longrightarrow A_{d-r-1}),
 \qquad \A_v=\bigcap_{e\ni v}\mathcal Z_e.
\]
Here $\mathcal Z_e$ is the zero-trace space of one edge. We call $\A_v$ the local annihilator: it consists of the polynomials whose dual trace vanishes on every edge incident to $v$.
\begin{proposition}[Compatible representatives and edge traces]\label{body:duality}
We have $J(v)_d^*\simeq A_d/\A_v$, and the families in $M_d^*$ satisfy
\begin{equation}\label{body:all-edge-traces}
 D_{\ell_e}^{r+1}H_u=D_{\ell_e}^{r+1}H_v=:G_e
 \qquad(e=uv).
\end{equation}
Choose a supporting-line equation $\ell_e$ for each boundary edge as well. Take $H_b=0$ at boundary vertices and $G_e=0$ on boundary edges. Every edge then has a $G_e$, which is zero whenever an endpoint is on the boundary. Replacing $H_v$ by $H_v-N_v$, $N_v\in\A_v$, preserves its functional and every $G_e$. All changes below are made at interior vertices; the boundary representatives remain $H_b=0$.
\end{proposition}
For adjacent interior vertices, \eqref{body:all-edge-traces} is also $H_u-H_v\in\mathcal Z_e$. Thus the matrix left kernel is encoded by a family of polynomials of fixed total degree, to be analyzed using their differential orders and the actual geometry.

\subsection{A unique common polynomial on each triangle}
The actual vertex $v=(x_v,y_v)$ gives the point form $P_v=x_vX+y_vY+Z$, satisfying $D_\ell P_v=\ell(x_v,y_v,1)$. On a nondegenerate triangle $\tau=(u,v,w)$, set $U=P_u,V=P_v,W=P_w$. They are independent, and the derivative on $uv$ is $\ell_{uv}(w,1)\partial_W$.
\begin{proposition}[Common face polynomial]\label{body:face}
At $d=3r+1$, each face $\tau$ has a unique $K_\tau\in A_d$ such that
\[
 D_{\ell_e}^{r+1}K_\tau=G_e\qquad(e\subset\tau).
\]
Every face containing a boundary vertex has $K_\tau=0$.
\end{proposition}
We also call the common face polynomial $K_\tau$ a face potential. Existence follows from mixed-derivative compatibility of the three traces. Uniqueness uses $3r+1>3r$: a monomial killed by all three coordinate derivatives of order $r+1$ has each exponent at most $r$. At the boundary, zero traces propagate the zero face polynomial along the open face fan; see Lemmas~\ref{lem:facepotential} and \ref{lem:boundaryfaces}.

Subtracting the residuals on adjacent faces at one vertex cancels its representative $H_v$. Each endpoint of the shared edge therefore contributes its point-form factor to the same difference $K_\tau-K_{\tau'}$. This converts the differential conditions within individual faces into restrictions on neighboring face differences.

\subsection{The residual allowed inside one face}
Keep $\tau=(u,v,w)$ and take $u$ as center. The derivatives of edges $uv,uw$ kill $H_u-K_\tau$, so the exponents of the transverse variables $V,W$ are at most $r$. Consequently,
\begin{equation}\label{body:radial-residual}
 \begin{aligned}
 H_u-K_\tau
 &=\sum_{0\le i,j\le r}c_{ij}V^iW^jU^{3r+1-i-j}\\
 &=\lambda_{u,\tau}U^{r+1}V^rW^r+U^{r+2}F_{u,\tau},
 \end{aligned}
\end{equation}
where $F_{u,\tau}\in A_{2r-1}$. At the next $U$ power, the transverse exponents sum to $2r-1$, allowing only $(r,r-1)$ and $(r-1,r)$. Thus
\begin{equation}\label{body:F-factor}
 F_{u,\tau}=V^{r-1}W^{r-1}(a_{u,\tau}V+b_{u,\tau}W)+UC_{u,\tau},
 \qquad C_{u,\tau}\in A_{2r-2}.
\end{equation}
We expand the residual in powers of the center point form $U$. A radial layer is the coefficient polynomial associated with one fixed power of $U$. The leading layer corresponds to $U^{r+1}$; the next layer, at $U^{r+2}$, is given by \eqref{body:F-factor} modulo $U$ and has two free coefficients. The following layer has three. Sections~\ref{body:odd} and \ref{body:even} compare these within-face forms with the restrictions on a shared-edge difference.

\subsection{Face differences, leading transport, and representatives}
For adjacent faces $\tau=(u,v,w)$ and $\tau'=(u,v,w')$, keep the actual point forms and write
\[
 W'=sW+pU+qV,\qquad s<0,\qquad s+p+q=1.
\]
Here $s$ is the ratio of signed doubled areas with matching vertex slots. It is negative because $w,w'$ lie on opposite sides of the edge. The coefficients sum to one because the $Z$ coefficient of each actual point form is one. Subtracting \eqref{body:radial-residual} at each endpoint gives
\begin{equation}\label{body:initialjump}
 K_{\tau'}-K_\tau=U^{r+1}V^{r+1}Q_e,\qquad Q_e\in A_{r-1}.
\end{equation}
Restricting this $Q_e$ to $U=0$ and $V=0$, and comparing the resulting expressions at $U=V=0$, gives
\begin{equation}\label{body:lambda-coupling}
 \lambda_{u,\tau'}=s^{-r}\lambda_{u,\tau},\qquad
 \lambda_{v,\tau'}=s^{-r}\lambda_{v,\tau},\qquad
 p\lambda_{v,\tau}=q\lambda_{u,\tau}.
\end{equation}
The full axis restrictions and remainder appear in \eqref{human:jump-axes}--\eqref{eq:explicitjump}. The raw $\lambda_{u,\tau}$ changes by $s^{-r}$. Multiplying it by the $r$th power of the signed doubled area, with the vertex order continued in matching slots, gives the vertex quantity independent of the face reached. Section~\ref{body:odd} checks cyclic closure and its parity signs.

At a boundary vertex, $H_b=0$ and all incident $K_\tau=0$, hence $\lambda_{b,\tau}=0$ by definition. At an interior vertex $u$ with only two supporting lines, choose neighbor forms $V,W$ on those lines. The polynomial $U^{r+1}V^rW^r$ belongs to $\A_u$. Subtracting a suitable multiple of it from $H_u$ removes the leading coefficient on every incident face. Make this representative choice first. At every remaining vertex, \eqref{body:lambda-coupling} says that vanishing on one face is equivalent to vanishing on all incident faces. Thus the nonzero leading support is a well-defined finite set of vertices.

\subsection{A normal form at an extreme support vertex}
If the normalized support $\Sigma$ is nonempty, choose a linear function with a unique maximum on $\Sigma$, and call the maximizing vertex $u$. Then $\Sigma\setminus\{u\}$ lies in a strict open half-plane through $u$. We call $u$ an extreme vertex of the support; the needed geometry follows from endpoint coupling and the cyclic neighbor order.
\begin{lemma}[Normal form at an extreme support vertex]\label{body:corner}
The neighbor cycle at $u$ consists of $n+1$ consecutive nonzero-support neighbors $p_0,\ldots,p_n$ and exactly two zero-support neighbors $p_{n+1},p_{n+2}$, where $n\ge2$. Orientation-preserving affine coordinates give their vectors from the center as
\begin{equation}\label{body:corner-rays}
 p_i-u\ \longleftrightarrow\ \rho_i(1-t_i,t_i)\quad(0\le i\le n),
 \qquad (-b,0),(0,-a),
\end{equation}
where $a,b,\rho_i>0$, $0=t_0<\cdots<t_n=1$, $\rho_0=\rho_n=1$, and the valence is $N=n+3$.
\end{lemma}
See Appendix Lemma~\ref{lem:corner-shape}. Put $e_0=p_0-u$, $e_n=p_n-u$, $E=(e_0\ e_n)$, and $D=\det E>0$. The two-way correspondence between physical and normal coordinates is
\begin{equation}\label{body:physical-change}
 \binom\alpha\beta=E^{-1}\left(\binom xy-\binom{x_u}{y_u}\right),
 \qquad \binom xy=\binom{x_u}{y_u}+E\binom\alpha\beta.
\end{equation}
In particular, the physical point is $p_i=u+\rho_i((1-t_i)e_0+t_ie_n)$. On the positive block, $\rho_i=\alpha(p_i)+\beta(p_i)$ and $t_i=\beta(p_i)/\rho_i$. The exterior neighbors are $u-be_0,u-ae_n$.

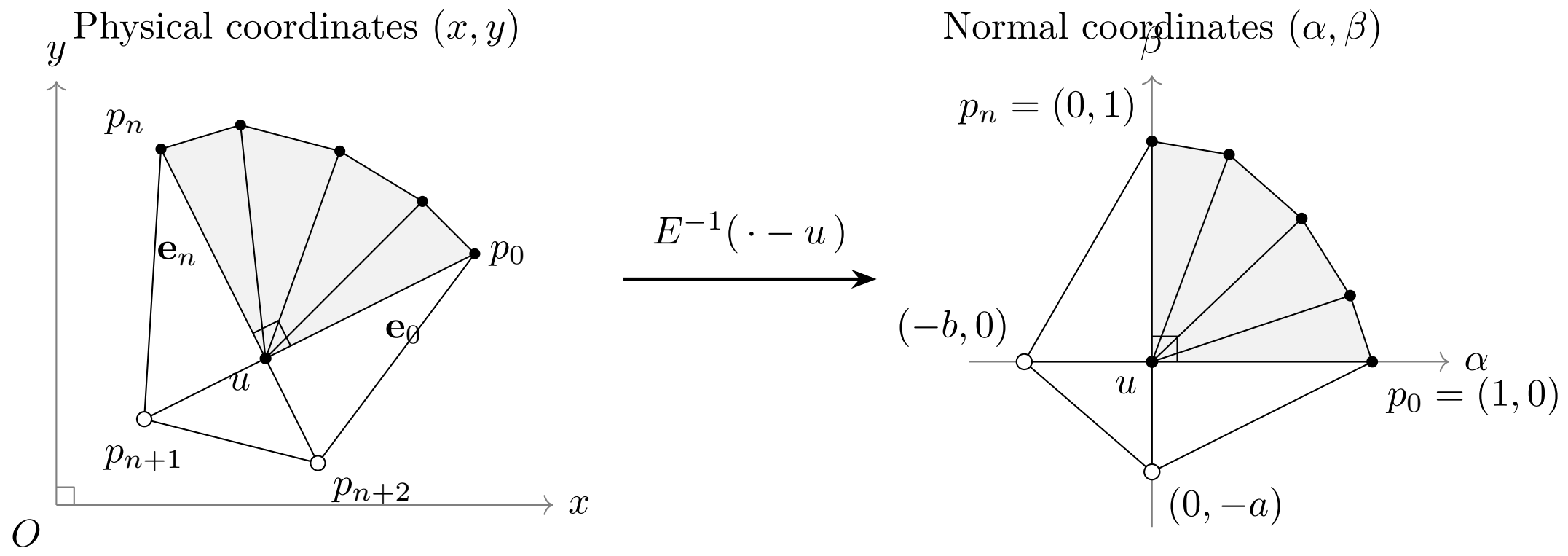
\begin{figure}[tb]
\centering
\begin{tikzpicture}[every node/.style={font=\small},line join=round]
\begin{scope}[xshift=.1cm,yshift=.2cm,scale=.95]
 \draw[->,gray] (0,0)--(4.75,0) node[right,black]{$x$};
 \draw[->,gray] (0,0)--(0,4.05) node[above,black]{$y$};
 \draw[gray] (.17,0)--(.17,.17)--(0,.17);
 \node[below left] at (0,0){$O$};
 \coordinate (P) at (2,1.4);
 \foreach \name/\x/\y in {A/4/2.4,B/3.5/2.9,C/2.71/3.38,D/1.76/3.63,E/1/3.4,F/.84/.82,G/2.5/.4}
  {\coordinate (\name) at (\x,\y);}
 \fill[gray!10] (P)--(A)--(B)--(C)--(D)--(E)--cycle;
 \foreach \name in {A,B,C,D,E,F,G}{\draw (P)--(\name);}
 \draw (A)--(B)--(C)--(D)--(E)--(F)--(G)--cycle;
 \draw (2.24,1.52)--(2.12,1.76)--(1.88,1.64);
 \foreach \name in {A,B,C,D,E}{\fill (\name) circle(1.5pt);}
 \foreach \name in {F,G}{\draw[fill=white] (\name) circle(2pt);}
 \fill (P) circle(1.6pt);
 \node[below left] at(P){$u$};
 \node[right] at(A){$p_0$};\node[above left] at(E){$p_n$};
 \node[below right] at(3,1.9){$\mathbf e_0$};
 \node[left] at(1.5,2.4){$\mathbf e_n$};
 \node[below=3pt] at(F){$p_{n+1}$};\node[below right] at(G){$p_{n+2}$};
 \node at(2.3,4.55){Physical coordinates $(x,y)$};
\end{scope}
\draw[-{Stealth[length=2.3mm]},thick] (5.25,2.25)--(7.55,2.25);
\node[above] at(6.4,2.35){$E^{-1}(\,\cdot-u\,)$};
\begin{scope}[xshift=10.05cm,yshift=1.5cm,scale=2]
 \draw[->,gray] (-.83,0)--(1.35,0) node[right,black]{$\alpha$};
 \draw[->,gray] (0,-.75)--(0,1.3) node[above,black]{$\beta$};
 \coordinate (P) at (0,0);
 \foreach \name/\x/\y in {A/1/0,B/.9/.3,C/.68/.65,D/.35/.94,E/0/1,F/-.58/0,G/0/-.5}
  {\coordinate (\name) at (\x,\y);}
 \fill[gray!10] (P)--(A)--(B)--(C)--(D)--(E)--cycle;
 \foreach \name in {A,B,C,D,E,F,G}{\draw (P)--(\name);}
 \draw (A)--(B)--(C)--(D)--(E)--(F)--(G)--cycle;
 \draw (.115,0)--(.115,.115)--(0,.115);
 \foreach \name in {A,B,C,D,E}{\fill (\name) circle(.75pt);}
 \foreach \name in {F,G}{\draw[fill=white] (\name) circle(1pt);}
 \fill (P) circle(.8pt);
 \node[below left] at(P){$u$};
 \node[below right] at(A){$p_0=(1,0)$};
 \node[above left] at(E){$p_n=(0,1)$};
 \node[above left] at(F){$(-b,0)$};\node[below right] at(G){$(0,-a)$};
 \node at(.05,1.51){Normal coordinates $(\alpha,\beta)$};
\end{scope}
\end{tikzpicture}
\caption{The same neighbor configuration at an extreme support vertex in physical coordinates and after the change $q=E^{-1}(p-u)$. Filled neighbors form the nonzero-support block; the two open neighbors are the zero block. The coordinate axes are perpendicular in both panels. The two extreme directions are also chosen perpendicular in this drawing; the proof only requires $\det E>0$.}
\label{fig:corner}
\end{figure}

To choose the corresponding basis of dual linear forms, use the original variables $X_{\rm old},Y_{\rm old},Z_{\rm old}$ and set
\[
 U=P_u,\qquad X_0=e_{0x}X_{\rm old}+e_{0y}Y_{\rm old},\qquad
 Y_0=e_{nx}X_{\rm old}+e_{ny}Y_{\rm old}.
\]
Put $\widehat X=X_0-U/b$ and $\widehat Y=Y_0-U/a$. This is invertible: $X_0=\widehat X+U/b$, $Y_0=\widehat Y+U/a$; recover the old $X,Y$ using $E^{-\mathsf T}(X_0,Y_0)^{\mathsf T}$ and then $Z_{\rm old}=U-x_uX_{\rm old}-y_uY_{\rm old}$. Therefore
\begin{equation}\label{body:corner-forms}
 L_i:=P_{p_i}/\rho_i=(1-t_i)\widehat X+t_i\widehat Y+\zeta_iU,
 \qquad \zeta_i=\rho_i^{-1}+(1-t_i)/b+t_i/a>0
\end{equation}
for $0\le i\le n$. At the exterior neighbors set $\rho_{n+1}=-b$, $\rho_{n+2}=-a$, and
\[
 (t_{n+1},\zeta_{n+1})=(0,0),\quad(t_{n+2},\zeta_{n+2})=(1,0),
 \qquad L_{n+1}=\widehat X,\quad L_{n+2}=\widehat Y.
\]
Use indices modulo $N$, with
\[
 h_i=t_{i+1}-t_i>0\ (0\le i<n),\quad\sum_{i=0}^{n-1}h_i=1,
 \qquad h_n=-1,\ h_{n+1}=1,\ h_{n+2}=-1.
\]
The positive doubled area in physical coordinates is $D\rho_i\rho_{i+1}h_i$. Both the negative exterior intervals and the negative exterior radial factors must therefore be retained. Substituting $\tau_i=\{u,p_i,p_{i+1}\}$ and $K_i=K_{\tau_i}$ into the preceding residual and face-difference formulas gives the common starting point for odd and even orders.

\subsection{From eliminating the leading coefficient to zero defect}
Eliminating the leading coefficient means choosing representatives in their existing classes so that $P_v^{r+2}\mid H_v-K_\tau$ everywhere. Sections~\ref{body:odd} and \ref{body:even} establish this for odd and even $r$, respectively. We then raise the divisibility order by one at a time, using a single statement independent of parity.

\begin{theorem}[Conditional higher-layer induction]\label{body:induction}
Let $p\ge1$ and $r\ge2p+1$. If a compatible family has representatives with
\[
 H_v-K_\tau\in P_v^{r+1+p}A_{2r-p}
\]
at every vertex and incident face, representatives of the same classes can be chosen so that
\[
 H_v-K_\tau\in P_v^{r+2+p}A_{2r-p-1}
\]
everywhere. All edge traces and face potentials are preserved. Changes are made only at interior vertices, leaving $H_b=0$ at the boundary.
\end{theorem}
The detailed statement and proof are Theorem~\ref{thm:induction}. At a fixed radial layer, the compatibility equations leave only two possible nonzero local patterns. One is common to all incident faces at a two-line vertex. We lift this common polynomial to the local annihilator before subtracting it from the representative, so the functional is unchanged. The other pattern can continue only along a straight line. Finiteness and the zero boundary states exclude this second pattern.

\begin{corollary}[Conditional completion]\label{body:completion}
For $r\ge3$, a dual defect class is zero if its representatives can be chosen so that $P_v^{r+2}\mid H_v-K_\tau$ at every vertex and incident face.
\end{corollary}
The completion is a degree argument, detailed in Theorem~\ref{thm:from-leading}. The successive increases in divisibility reach point-form exponent
\[
 a_*=r+1+\left\lceil\frac r2\right\rceil,
 \qquad 2a_*>3r+1.
\]
Every face difference is therefore zero. Sections~\ref{body:odd} and \ref{body:even} need only eliminate the leading coefficient.

\section{The dimension formula for odd smoothness orders}\label{body:odd}

Fix an odd integer $r\ge3$. By Section~\ref{body:framework}, it suffices to eliminate the nonzero leading layer of a compatible functional family and then apply the conditional higher-layer induction. We first keep $r\ge3$ general and compare two expressions in the same standard coordinates: the coefficients allowed within each face, and the restrictions imposed by the actual difference of neighboring face potentials. Parity enters only in the final sign argument.

\subsection{A common leading scale and cyclic signs}

Fix the center $u$ and its representative. On the initial ordered face $\tau=(u,p_0,p_1)$, put $\Delta_\tau=\det(p_0-u,p_1-u)$ and define
\begin{equation}\label{body:vertex-normalization}
 \Lambda_u:=\lambda_{u,\tau}\Delta_\tau^r.
\end{equation}
At a crossing, retain the slots of the center and shared neighbor and replace only the other neighbor. The transport law of Section~\ref{body:framework} gives
\begin{equation}\label{body:leading-transport}
 \Delta_{\tau'}=s\Delta_\tau,\qquad
 \lambda_{u,\tau'}=s^{-r}\lambda_{u,\tau},\qquad s<0.
\end{equation}
Thus $\lambda_{u,\tau'}\Delta_{\tau'}^r=\Lambda_u$ is invariant along this ordered face chain.
After one complete circuit, the continued signed area is $(-1)^N\Delta_\tau$. An exchange of neighbor slots leaves the leading coefficient unchanged, since the two neighbor forms both have exponent $r$. Hence
\[
 \Lambda_u=(-1)^{rN}\Lambda_u.
\]
At odd $r$, a nonzero leading layer therefore requires even valence $N$. The areas and continued orders remain signed throughout; equations \eqref{eq:lambda-transport}--\eqref{eq:area-normalized-transport} give the complete derivation.

First remove the adjustable terms at two-line vertices as in Section~\ref{body:framework}. If nonzero leading support remains, choose the unique maximizer $u$ of a linear function on this finite support and use the standard form of Lemma~\ref{body:corner}. Write $\tau_i=\{u,p_i,p_{i+1}\}$, $K_i=K_{\tau_i}$, $U=P_u$, retaining $L_i,t_i,h_i,\zeta_i$ from Section~\ref{body:framework}. The complete valence is $N=n+3$, and the exterior intervals keep the signs $h_n=-1,h_{n+1}=1,h_{n+2}=-1$.

Substituting the physical points $p_i=u+\rho_i((1-t_i)e_0+t_ie_n)$ and $p_{i+1}$ gives $P_{p_i}=\rho_iL_i$ and the actual doubled area $\Delta_i^*=D\rho_i\rho_{i+1}h_i$. Thus replacing actual neighbor forms by $L_i,L_{i+1}$ multiplies the leading coefficient by $(\rho_i\rho_{i+1})^r$.

Counterclockwise orders on the faces differ from the continued orders by $i$ reversals, so put $\sigma_i=(-1)^{ri}$. With $D=\det E>0$, in the initial radial unit $U=P_u$, set $\Lambda=\Lambda_u/D^r\ne0$. Substitution of the actual point forms into the vertex residual gives
\begin{equation}\label{body:normalized-first}
 \begin{aligned}
 R_i&:=\frac{H_u-K_i}{\Lambda}=U^{r+1}M_i+U^{r+2}F_i,\\
 M_i&=\frac{\sigma_i}{h_i^r}L_i^rL_{i+1}^r.
 \end{aligned}
\end{equation}
Equation~\eqref{body:normalized-first} normalizes the residual by the common scalar $\Lambda$; $H_u$ and $K_i$ remain fixed.

\subsection{The first calculation: endpoint coefficients within each face}

The two transverse differential conditions and the degree bound give
\begin{equation}\label{body:piece-factorization}
 F_i=L_i^{r-1}L_{i+1}^{r-1}(a_iL_i+b_iL_{i+1})+UG_i.
\end{equation}
In the basis $(\widehat X,\widehat Y,U)$, first set $U=0$ and then substitute $(-t,1-t)$. Since $L_j(-t,1-t,0)=t_j-t$,
\begin{equation}\label{body:hermite-motivation}
 \begin{aligned}
 f_i(t)&:=F_i(-t,1-t,0)\\
 &=(t_i-t)^{r-1}(t_{i+1}-t)^{r-1}
 [a_i(t_i-t)+b_i(t_{i+1}-t)].
 \end{aligned}
\end{equation}
This is the source of the Hermite structure: $f_i$ has degree at most $2r-1$ and vanishes to order at least $r-1$ at each endpoint. The two endpoint derivatives of order $r-1$ determine its remaining affine factor.

We write $[x^m]f$ for the coefficient of $x^m$ in the polynomial $f$. To compare the common node $t_i$, set $x=t-t_i$. The right piece $f_i$ and the left piece $f_{i-1}$ separately give
\begin{equation}\label{body:separate-endpoints}
 \begin{aligned}
 A_i^+&:=[x^{r-1}]f_i(t_i+x)=(-1)^{r-1}b_ih_i^r,\\
 A_i^-&:=[x^{r-1}]f_{i-1}(t_i+x)=-a_{i-1}h_{i-1}^r.
 \end{aligned}
\end{equation}
The quantities $A_i^+$ and $A_i^-$ are initially defined independently. Their equality follows from the face-difference calculation.

\subsection{The second calculation: shared data from the edge jump}

Define $J_i=(K_i-K_{i-1})/\Lambda$. Exact subtraction in \eqref{body:normalized-first} gives
\begin{equation}\label{body:first-difference}
 U^{r+2}(F_i-F_{i-1})
 =-J_i-U^{r+1}(M_i-M_{i-1}).
\end{equation}
This subtraction identity must also use the computation of the same $J_i$ at the other endpoint of the shared edge. The two endpoint factors and axis restrictions from Section~\ref{body:framework} determine its quotient modulo the product of the endpoint forms. Set
\[
 h=h_{i-1},\quad k=h_i,\qquad
 \kappa_i=\frac{\zeta_{i+1}-\zeta_i}{h_i},\qquad
 \delta_i=\kappa_i-\kappa_{i-1}.
\]
Thus $\delta_i$ is the difference between the secant slopes on the two adjacent intervals. Substitute both axis restrictions into \eqref{body:first-difference}, extract the present radial layer, and expand at $t_i$:
\begin{equation}\label{body:first-jump-expansion}
 \begin{aligned}
 f_i(t_i+x)-f_{i-1}(t_i+x)
 ={}&-r\sigma_{i-1}\delta_i x^r\\
 &+\binom r2\sigma_{i-1}\delta_i(1/k-1/h)x^{r+1}
 +O(x^{r+2}).
 \end{aligned}
\end{equation}
The undetermined edge remainder starts at $x^{r+2}$ and does not affect the three coefficients being compared. Appendix~\ref{app:coefficients}, equations \eqref{eq:axis1}--\eqref{eq:fdiff3}, gives the axis restrictions and coefficient expansion.

First compare $x^{r-1}$. Its coefficient in \eqref{body:first-jump-expansion} is zero, proving $A_i^+=A_i^-$. We may therefore define the common node value
\begin{equation}\label{body:hermite-data}
 \begin{aligned}
 A_i&:=A_i^+=A_i^-
 =\frac{f_i^{(r-1)}(t_i)}{(r-1)!}
 =\frac{f_{i-1}^{(r-1)}(t_i)}{(r-1)!},\\
 \beta_i&:=A_i/\sigma_{i-1}.
 \end{aligned}
\end{equation}
Using $(-1)^{r-1}\sigma_i=-\sigma_{i-1}$, equation~\eqref{body:separate-endpoints} becomes
\begin{equation}\label{body:ab-beta}
 a_i=-\frac{\sigma_i}{h_i^r}\beta_{i+1},\qquad
 b_i=-\frac{\sigma_i}{h_i^r}\beta_i.
\end{equation}
The two free coefficients on each face have thus been replaced by the shared data at its left and right nodes. Write
\[
 \begin{aligned}
 \Pi_i&=-\frac{\sigma_i}{h_i^r}L_i^{r-1}L_{i+1}^{r-1}
          (\beta_{i+1}L_i+\beta_iL_{i+1}),\\
 F_i&=\Pi_i+UG_i.
 \end{aligned}
\]

\subsection{Coefficient matching and the odd-order contradiction}

Substitute \eqref{body:ab-beta} into the first calculation and compare with \eqref{body:first-jump-expansion}. The coefficients of $x^r$ and $x^{r+1}$ give the two original Hermite equations:
\begin{align}
 r\delta_i&=\frac{\beta_{i-1}+r\beta_i}{h}
              +\frac{r\beta_i+\beta_{i+1}}k,\label{body:first-hermite1}\\
 \frac r2\delta_i(1/k-1/h)
 &=\frac r2\beta_i(1/k^2-1/h^2)
       +\frac{\beta_{i+1}}{k^2}-\frac{\beta_{i-1}}{h^2}.
       \label{body:first-hermite2}
\end{align}
Subtract $(1/k-1/h)/2$ times the first equation from the second:
\[
 \frac{h+k}{2hk}
 \left(\frac{\beta_{i+1}}k-\frac{\beta_{i-1}}h\right)=0.
\]
Within the block $h,k>0$. At the two ends, retain the signed exterior intervals and use $h_0,h_{n-1}<1$. This yields the two-step recurrence. At the two exterior nodes, do not cancel the vanishing factor $h+k$: substitute into the original equations instead. Summing over the full cycle then gives the following conditions.

\begin{lemma}[Hermite compatibility]\label{body:hermite}
The common endpoint data satisfy
\begin{align}
 \frac{\beta_{i+1}}{h_i}
 &=\frac{\beta_{i-1}}{h_{i-1}}
 &&(0\le i\le n),\label{body:hermite-recurrence}\\
 \beta_{n+1}&=\beta_0-r\zeta_0,\qquad
 \beta_{n+2}=\beta_n-r\zeta_n,\label{body:hermite-boundary}\\
 \sum_{i=0}^{n-1}\frac{\beta_i+\beta_{i+1}}{h_i}
 &=\beta_0+\beta_n.\label{body:hermite-flux}
\end{align}
\end{lemma}
Appendix Lemma~\ref{lem:hermite} gives the complete passage from the two original equations to the endpoint conditions and cyclic sum identity.

Parity is used at this point. For odd $r$, a nonzero leading layer requires even $N=n+3$, so $n$ is odd. The recurrence splits the $\beta_i$ into two positive-multiplier chains. The exterior conditions and $\zeta_0,\zeta_n>0$ make both chains positive. Thus every $\beta_i>0$. Since each $0<h_i<1$, the left side of the cyclic sum identity strictly exceeds $\beta_0+\beta_n$, contradicting \eqref{body:hermite-flux}. Proposition~\ref{prop:oddcorner} gives the detailed sign calculation.

\begin{lemma}[Eliminating the leading coefficient for odd orders]\label{body:odd-leading}
For every odd $r\ge3$, we may replace the representatives $H_v$ of any dual defect class by equivalent representatives modulo $\A_v$ so that
\begin{equation}\label{body:odd-gain}
 H_v-K_\tau\in P_v^{r+2}A_{2r-1}
\end{equation}
at every incident vertex and face. Auxiliary boundary representatives remain zero.
\end{lemma}
Any remaining leading support would have an extreme point with the standard form just excluded. The leading layer therefore vanishes everywhere. Normalization subtracts local annihilator terms only at two-line vertices, preserving the original functional classes, every edge trace, and every face potential; see Lemma~\ref{lem:odd-leading}.

\subsection{Higher-layer induction and the odd-order dimension formula}

Lemma~\ref{body:odd-leading} supplies the hypothesis of the completion result, Corollary~\ref{body:completion}. More explicitly, applying Theorem~\ref{body:induction} successively for $p=1,\ldots,(r-1)/2$ gives representatives satisfying
\begin{equation}\label{body:odd-depth}
 H_v-K_\tau\in P_v^{r+1+(r+1)/2}A_{2r-(r+1)/2}.
\end{equation}
The endpoint factors have total degree $3r+3>3r+1$, so Corollary~\ref{body:completion} gives $M_{3r+1}=0$.

\begin{theorem}[Every odd smoothness order]\label{body:oddcase}
The main dimension formula holds for every odd $r\ge3$.
\end{theorem}
Appendices \ref{sec:higher}--\ref{sec:stars} give the higher-coefficient transfer, classification, and representative changes; Theorem~\ref{thm:from-leading} contains the degree-completion proof. These steps require only the elimination of the leading coefficient and will also apply once it is established for even orders.

\section{Even smoothness orders: the second radial layer}\label{body:even}

Fix an even integer $r\ge4$. Use the standard configuration of Section~\ref{body:framework}, with center form $U=P_u$, neighbors $p_i$, and face potentials $K_i$. We again compare what is allowed inside each face with the face difference computed from the two endpoints of a shared edge. The first comparison supplies the node data $\beta_i$. Even orders require the same comparison once more for the remainder $G_i$. The two rounds supply a cyclic sum identity and direction constraints, which will act on the same path.

\subsection{The first Hermite compatibility conditions}

The face calculation and face-difference calculation in Section~\ref{body:odd} do not use the parity of $r$. At even orders $\sigma_i=1$, so they give
\begin{equation}\label{body:even-first-round}
 R_i:=\frac{H_u-K_i}{\Lambda}
 =U^{r+1}M_i+U^{r+2}\Pi_i+U^{r+3}G_i,
\end{equation}
where
\begin{equation}\label{body:even-Pi}
 \begin{aligned}
 M_i&=h_i^{-r}L_i^rL_{i+1}^r,\\
 \Pi_i&=-h_i^{-r}L_i^{r-1}L_{i+1}^{r-1}
       (\beta_{i+1}L_i+\beta_iL_{i+1}).
 \end{aligned}
\end{equation}
Here $F_i=\Pi_i+UG_i$ is an exact decomposition of the previous remainder. The coefficient of $x^{r-1}$ in the face difference first proves that endpoint data are shared; comparison of $x^r,x^{r+1}$ then gives the recurrence, exterior relations, and cyclic sum identity of Lemma~\ref{body:hermite}. At odd orders these exclude a nonzero leading layer. At even orders the two chains can have opposite signs and their contributions to the cyclic sum can cancel, so the remainder $G_i$ must be retained.

The common calculation excludes odd $n$, and Lemma~\ref{lem:short-corner} excludes $n=2$, leaving even $n\ge4$. The recurrence and cyclic sum identity make all block coefficients $\beta_i$ nonzero with alternating signs, and make $\beta_i\beta_{i+1}/h_i$ one common negative constant. After the positive rescaling of $U$ described before \eqref{eq:product-normalization}, we may assume
\begin{equation}\label{body:even-product}
 h_i=-\beta_i\beta_{i+1}>0\qquad(0\le i<n).
\end{equation}
Henceforth $U$ is a positive multiple of $P_u$. The actual $H_u,K_i,L_i$ and functional classes are unchanged. 

\subsection{Compatibility of the next radial coefficients}

First use the face expression. The degree of $G_i\bmod U$ is $2r-2$, and both transverse exponents are at most $r$. The only pairs are $(r,r-2)$, $(r-1,r-1)$, and $(r-2,r)$. Thus $G_i\bmod U$ contains $L_i^{r-2}L_{i+1}^{r-2}$ and has a quadratic factor left over. Each face initially has three independent coefficients; the endpoint coefficients on the two sides of a common node are still distinct data.

Now use the face difference. Subtract \eqref{body:even-first-round} for faces $i-1,i$ to obtain the exact identity
\begin{equation}\label{body:even-second-difference}
 \begin{aligned}
 U^{r+3}(G_i-G_{i-1})={}&-\frac{K_i-K_{i-1}}{\Lambda}
 -U^{r+1}(M_i-M_{i-1})\\
 &-U^{r+2}(\Pi_i-\Pi_{i-1}).
 \end{aligned}
\end{equation}
The face difference on the right still uses the two endpoint axis restrictions from the first round. Work in the common coordinates
\begin{equation}\label{body:even-common-coordinates}
 \begin{gathered}
 A=L_i,\quad W=L_{i-1},\qquad L_{i+1}=s_iW+q_iA+\xi_iU,\\
 s_i=-h_i/h_{i-1},\quad q_i=1+h_i/h_{i-1},\quad
 \xi_i=h_i\delta_i.
 \end{gathered}
\end{equation}
Because the neighbor forms have been rescaled, $s_i$ need not be negative; at the two exterior nodes, $s_i=1$. The radial coefficient must be extracted in these three variables: each $L_j$ also contains the center form $U$, so the already known $M_i,\Pi_i$ contribute to the next layer.

Divide \eqref{body:even-second-difference} by $U^{r+1}$, extract the coefficient of $U^2$ in the coordinates \eqref{body:even-common-coordinates}, and divide by $A^{r-2}$. Comparing the constant terms then shows that the endpoint coefficients at the shared node agree. We may therefore write the three coefficients as
\begin{equation}\label{body:even-shared-gamma}
 G_i\equiv h_i^{-r}L_i^{r-2}L_{i+1}^{r-2}
 (\gamma_{i+1}L_i^2+\eta_iL_iL_{i+1}+\gamma_iL_{i+1}^2)
 \pmod{UA_{2r-3}}.
\end{equation}
Here $\gamma_i$ is shared node data and $\eta_i$ is the middle coefficient of face $i$. The coefficients of $A,A^2$ give two transfers for $\eta,\gamma$; that of $A^3$ supplies an additional compatibility condition. The unknown shared-edge remainder starts at $A^4$ and affects none of these comparisons. The third raw coefficient contains $q_i(r-2)$, which is nonzero at block nodes $0\le i\le n$. At the two exterior nodes $q_i=0$ and only the first two equations remain. The expansions and raw coefficients are given in \eqref{eq:second-0}--\eqref{eq:rawjet3}; the resulting transfer equations are \eqref{eq:eta}--\eqref{eq:thirdjet}.

\subsection{Positivity, reflection symmetry, and parity separation}

At block nodes $0\le i\le n$, the first-round equations show that $\gamma_i^0=\beta_i^2$, $\eta_i^0=\beta_i\beta_{i+1}$ solve the first two transfers. Subtract this particular solution and define
\begin{equation}\label{body:even-states}
 \mathsf u_i=\frac{\gamma_i-\beta_i^2}{\beta_i^2},\qquad
 \mathsf v_i=\frac{\eta_i-\beta_i\beta_{i+1}}{\beta_i\beta_{i+1}}.
\end{equation}
The scalar $\mathsf u_i$ is attached to node $i$, and $\mathsf v_i$ to interval $i$ between neighboring directions. Define $d_i=-q_i\beta_i/\beta_{i+1}$. Using \eqref{body:even-product} and the exterior recurrence expresses these parameters, determined by the direction intervals and node data, as
\begin{equation}\label{body:even-parameters}
 \begin{aligned}
 d_0&=(1/h_0-1)\beta_0^2,&d_n&=(1/h_{n-1}-1)\beta_n^2,\\
 d_i&=(1/h_{i-1}+1/h_i)\beta_i^2&&(1\le i<n).
 \end{aligned}
\end{equation}
In particular $d_i>0$. At $0\le i\le n$, the three second-round relations become
\begin{equation}\label{body:even-three-relations}
 \begin{aligned}
 \mathsf v_i+\mathsf v_{i-1}&=rd_i\mathsf u_i,\\
 \mathsf u_{i+1}-\mathsf u_{i-1}
 &=\frac{r-1}{2}d_i(\mathsf v_i-\mathsf v_{i-1}),\\
 \mathsf u_{i-1}+\mathsf u_{i+1}
 -\frac{r(r-1)}6d_i^2\mathsf u_i
 &=-\frac{2(r+2)}{3r}-\frac{r-1}{3}d_i.
 \end{aligned}
\end{equation}
Applying the first two transfer relations at the exterior nodes $n+1,n+2$ gives the endpoint conditions \eqref{eq:exterior-u}, where $\mathsf R=(\beta_0\beta_n)^2$. These two nodes have $q_i=0$; they are distinct from the three exterior intervals included in cyclic sums.

First, the third equation and endpoint conditions form a strictly diagonally dominant linear system with positive right-hand side, giving $\mathsf u_i>0$. The first two equations imply that
\[
 I_i=\mathsf u_i\mathsf u_{i+1}-\frac{r-1}{2r}\mathsf v_i^2
\]
is independent of $i$. Together with the exterior relations, this gives $\mathsf u_0=\mathsf u_n$. After reversing the index order, the endpoint data are the same. Uniqueness of positive continuation therefore forces the two inward sequences to agree, so
\[
 \mathsf u_i=\mathsf u_{n-i},\qquad d_i=d_{n-i}.
\]
Combining this symmetry with the adjacent-difference inequalities compares the global extrema of the two groups: either every even-indexed $d_i$ exceeds every odd-indexed $d_i$, or the reverse holds. The details are given in Lemmas~\ref{lem:positive}, \ref{lem:branch}, and Propositions~\ref{prop:reflection}, \ref{prop:parity}.

\begin{proposition}[Reflection symmetry and parity separation]\label{body:even-geometry}
An actual nonzero even-order leading configuration satisfies $d_i=d_{n-i}$ and $\mathsf u_i=\mathsf u_{n-i}$.
Moreover, some $\varepsilon\in\{1,-1\}$ satisfies
\begin{equation}\label{body:even-separation}
 \varepsilon(d_{2j}-d_{2k+1})>0
 \qquad(0\le2j,2k+1\le n).
\end{equation}
Thus the closed intervals spanned by all even-indexed and all odd-indexed direction parameters are disjoint.
\end{proposition}

\subsection{Moment identities and a balanced path}

Using the previously defined secant-slope difference $\delta_i=\kappa_i-\kappa_{i-1}$, multiply the first-round Hermite equation \eqref{body:first-hermite1} by $1$ and by $t_i$, respectively, and sum around the complete cycle, including all three exterior intervals. Since
\[
 \sum_i\delta_i=0,\qquad
 \sum_it_i\delta_i=-\sum_ih_i\kappa_i
 =-\sum_i(\zeta_{i+1}-\zeta_i)=0,
\]
the signed exterior lengths give
\begin{equation}\label{body:even-moments}
 \sum_{i=0}^n(-1)^i\frac{d_i}{|\beta_i|}=0,\qquad
 \sum_{i=0}^n(-1)^i\frac{d_it_i}{|\beta_i|}=0.
\end{equation}
Both moments come from the actual first-round equations. Cyclic reindexing and exterior cancellation are given in \eqref{eq:moment0}--\eqref{eq:moment1}.

\begin{proposition}[A balanced even path from the same configuration]\label{body:even-path}
Put
\[
 g_i=\frac{1-t_i}{|\beta_i|}\quad(0\le i\le n),\qquad g_{-1}=0.
\]
Then $g_0,\ldots,g_{n-1}>0$, $g_n=0$, and
\begin{equation}\label{body:even-path-balance}
 \frac{g_{i-1}+g_{i+1}}{g_i}=d_i\quad(0\le i<n),\qquad
 \sum_{i=0}^{n-1}(-1)^i\nu_i g_i=0,
\end{equation}
where $\nu_0=\nu_{n-1}=1$ and $\nu_i=2$ internally.
\end{proposition}
The ratio identity uses $h_i=-\beta_i\beta_{i+1}$; degree balance uses both zero moments. There are exactly the even number $n$ of positive values: geometric node $n$ has become the right zero endpoint. The calculation follows \eqref{eq:actualpath} in the appendix.

\begin{theorem}[Interval overlap on a balanced even path]\label{body:even-path-overlap}
Let $N\ge2$ be even, $g_0,\ldots,g_{N-1}>0$, and $g_{-1}=g_N=0$. Define
\[
 \mu_i=(g_{i-1}+g_{i+1})/g_i,\qquad
 \nu_i=\begin{cases}1,&i=0,N-1,\\2,&1\le i\le N-2.\end{cases}
\]
If $\sum_{i=0}^{N-1}(-1)^i\nu_i g_i=0$, then
\[
 [\min_{i\ \mathrm{even}}\mu_i,\max_{i\ \mathrm{even}}\mu_i]
 \cap[\min_{i\ \mathrm{odd}}\mu_i,\max_{i\ \mathrm{odd}}\mu_i]
 \ne\varnothing.
\]
\end{theorem}
Theorem~\ref{thm:path} proves this by a nonnegative dual certificate, whose finite sine blocks are constructed in Lemma~\ref{lem:pathdual}.

Removing the geometric node $n$ whose path value is zero preserves strict separation of the ratios $d_0,\ldots,d_{n-1}$. Interval overlap on the balanced even path contradicts this separation in Proposition~\ref{body:even-geometry}. A nonzero leading configuration is therefore impossible.

\begin{theorem}[Eliminating the leading coefficient at even orders]\label{body:even-leading}
For every even integer $r\ge4$, each representative $H_v$ of a dual defect class may be replaced by an equivalent representative modulo $\A_v$ such that
\[
 H_v-K_\tau\in P_v^{r+2}A_{2r-1}
\]
at every incident vertex and face, while auxiliary boundary representatives remain zero.
\end{theorem}
After the two-direction normalization, any nonempty remaining leading support has an extreme point giving the excluded configuration. The normalized leading coefficient therefore vanishes everywhere. Corollary~\ref{body:completion} completes the even orders; together with Theorem~\ref{body:oddcase}, this proves Theorem~\ref{thm:main}.

\section{Applying the framework at \texorpdfstring{$r=1,2$}{r=1,2}}\label{body:low}

The local vertex spaces, global endpoint map, vertex functionals, and unique common face potentials also exist for $r=1,2$. The distinction arises in leading-layer compatibility: the factor $r-1$ in the general Hermite elimination vanishes at $r=1$, and the factor $r-2$ in the second round vanishes at $r=2$. We give the low-order residuals and explain how the known quartic and septic dimension theorems fit this framework. Appendix~\ref{sec:septic-axes} derives the two septic identities used below. Let $n_j$ count interior vertices with exactly $j$ distinct unoriented supporting lines.

\subsection{Quartic splines with one continuous derivative: termination after the leading layer}
Take $r=1,d=4$. On a face $\tau=(u,v,w)$ put $U=P_u,V=P_v,W=P_w$. The two incident differential kernels give the residual at $u$ the form
\[
 H_u-K_\tau=\lambda_{u,\tau}U^2VW+U^3F_{u,\tau},
 \qquad F_{u,\tau}\in A_1.
\]
The unique face potential and endpoint coupling of leading coefficients still hold. Leading coefficients at two-direction vertices can still be removed by changes of representatives modulo the local annihilators. If the leading layer vanishes everywhere, then
$U^3V^3\mid K_\tau-K_{\tau'}$ along a shared edge $uv$. Since $6>4$, every face difference is zero. Boundary zero potentials and connected face adjacency give $K_\tau=0$, and the local annihilators then show that the original left kernel is zero. 

\begin{remark}[The low-order input for the leading layer]\label{body:quartic-scope}
At $r=1$, $F_i\bmod U$ is only linear. The two expressions can still be compared, but the factor $r-1$ vanishes when extracting the second Hermite relation. The two-step recurrence of Section~\ref{body:odd} does not follow. Here the quartic dimension theorem of \citet*{APS87} gives $M_4=0$. Thus the common dual framework admits zero representatives, which supply the required leading-layer elimination.
\end{remark}

The local formula gives $\sigma_v(1)=1$ if $t_v=2$ and $\sigma_v(1)=0$ if $t_v\ge3$. Substitution in the same global dimension identity gives
\begin{theorem}[The quartic formula, \citealp*{APS87}]\label{body:quartic}
For every disk triangulation considered here,
\[
 \dim S_4^1(\Delta)=15+6(E_I-2V_I)+n_2.
\]
\end{theorem}

\subsection{Septic splines with two continuous derivatives: the axes determine the full face difference}
Take $r=2,d=7$. The center residual is
\[
 H_u-K_\tau=\lambda_{u,\tau}U^3V^2W^2+U^4F_{u,\tau},
 \qquad F_{u,\tau}\in A_3.
\]
The signed-area normalization $\Lambda_u=\lambda_{u,\tau}\Delta_\tau^2$ agrees with positive-area normalization. The first round still supplies the shared $\beta_i$, both Hermite equations, and the cyclic sum identity.

There is a further septic simplification: the unspecified part of the full face-difference quotient belongs to $A_{r-3}=A_{-1}=0$. The two endpoint axes therefore determine the full face difference. In the standard configuration of Section~\ref{body:framework}, put
$A=L_i$, $W=L_{i-1}$, and $T=L_{i+1}=s_iW+q_iA+\xi_iU$.
Substituting $r=2$ in the axis expressions gives the exact identity
\begin{equation}\label{body:septic-exact-jump}
 \frac{K_i-K_{i-1}}{\Lambda}
 =-\frac{q_i}{h_i^2}U^3A^3(s_iW+T).
\end{equation}
Subtracting the two face leading terms and using $\xi_i=h_i\delta_i$ gives
\begin{equation}\label{body:septic-cubic-difference}
 F_i-F_{i-1}=-\delta_i L_i^2
 \left(\frac{L_{i+1}}{h_i}-\frac{L_{i-1}}{h_{i-1}}\right).
\end{equation}
This displays the connection between the two computations for the septic case: compute the same face difference from the two shared-edge endpoints, then compare with the allowed cubic remainder inside a face. The septic leading-layer elimination uses the separate compatibility argument of \citet{WD26}. Here we compare its scope with the general even-order conditions; Appendix~\ref{app:degree} derives the two septic identities directly.

\begin{remark}[Degeneration of the second round at degree seven]\label{body:septic-scope}
The third raw coefficient for general even orders contains $q_i(r-2)$ and is identically zero at $r=2$. The three-equation system of Section~\ref{body:even} cannot replace the septic closure argument. Once the septic leading layer vanishes, however, each residual contains $P_v^4$. The endpoint factors then have total degree $8>7$, and the framework immediately gives a zero left kernel without higher-layer induction.
\end{remark}

The local corrections are $3$ if $t_v=2$, $1$ if $t_v=3$, and $0$ if $t_v\ge4$. Hence
\begin{theorem}[The septic formula, \citealp{WD26}]\label{body:septic}
For every disk triangulation considered here,
\[
 \dim S_7^2(\Delta)=36+15(E_I-2V_I)+3n_2+n_3.
\]
\end{theorem}

Thus the local dimensions, functional encoding, common face potentials, and final degree comparison apply directly at both lower orders. Each requires its own leading-layer input. At $r\ge3$, a shared-edge quotient of degree $r-3$ can remain after the leading coefficient has been eliminated, which is why the uniform higher-layer induction is then needed.

\section{Conclusion}\label{body:conclusion}

We have proved that, for every integer $r\ge3$, every nondegenerate geometric realization of a triangulated polygonal disk attains Schumaker's dimension lower bound at degree $3r+1$, as stated in Theorem~\ref{thm:main}. The local correction uses the actual supporting-line counts of the given realization. Together with the known quartic $C^1$ and septic $C^2$ results and Hong's higher-degree theorem, Corollary~\ref{cor:all-orders} gives the uniform sufficient range $r\ge1$, $d\ge3r+1$.

The key mechanism is to express compatibility of vertex functionals as polynomial conditions, so that both endpoints control one face difference. Hermite compatibility eliminates the leading coefficient at odd orders. At even orders, the next radial layer forces direction separation, contradicting interval overlap on a balanced path. The remaining coefficients are handled by a common classification on complete stars and induction through the radial layers.

Degree $3r+1$ ensures that the three edge traces uniquely determine a face polynomial. At degree $3r$ they fail to detect $U^rV^rW^r$, leaving one free scalar on each face. Coordinating these scalars requires a new argument across neighboring faces, so the present method does not extend directly to that degree. In particular, the question posed by \citet{WD26} of whether $S_6^2$ attains Schumaker's lower bound on every nondegenerate geometric realization remains unresolved.

\appendix
\section{Proofs of the algebraic framework}\label{app:framework}

On a general triangle $(u,v,w)$ this appendix uses $U=P_u,V=P_v,W=P_w$. Thus $V$ is the central form when the vertex under discussion is $v$. In the normal form in the main text the center is named $u$, and its form is $U$; this only changes the vertex names.

\subsection{Local vertex spaces and their dimensions}\label{human:sec:local}
\subsubsection{Separate the homogeneous binary pieces}
Set $R=\R[x,y,z]$, and let $R_d$ denote its degree-$d$ homogeneous part. Homogenization identifies $R_d$ with bivariate polynomials of degree at most $d$. For an interior edge $e$ with supporting-line equation $\ell_e$, use the local ideals \citep{SS97,MV13}
\[
 J(e)=(\ell_e^{r+1}),\qquad J(v)=\sum_{e\ni v}J(e).
\]
Only the supporting line matters: multiplying $\ell_e$ by a nonzero scalar leaves these spaces unchanged. At an interior vertex $v$, $J(v)_d$ is the space spanned by all allowed edge differences incident to $v$. Its dimension is therefore the number of independent coefficients in the vertex closure equation.

Translate $v$ to the affine origin. Every incident line form is now binary, $\ell_i=a_ix+b_iy$. Repeated supporting lines have the same contribution, so retain $t=t_v$ distinct forms and let
\[
 B=\R[x,y],\qquad I=(\ell_1^{r+1},\ldots,\ell_t^{r+1})\subset B.
\]
The homogenizing variable is free. Grouping by binary degree gives the direct sum
\begin{equation}\label{human:binary-pieces}
 J(v)_d=\bigoplus_{k=r+1}^d z^{d-k}I_k,
 \qquad I_k=\sum_{i=1}^t\ell_i^{r+1}B_{k-r-1}.
\end{equation}
This makes the local dimension problem finite and explicit. The ambient space $B_k$ has dimension $k+1$, and each summand has dimension $k-r$. Their sum has dimension at most $\min\{k+1,t(k-r)\}$. We next prove that this bound is attained for every set of distinct real supporting lines.

\subsubsection{First proof: the rank of an explicit coefficient matrix}
Choose a coordinate direction not parallel to any incident line and rescale the line forms to $\ell_i=x+c_i y$, with distinct $c_i$. Fix $k\ge r+1$ and write $m=k-r$. Since $(\ell_i,y)$ are binary coordinates,
\[
 \ell_i^{r+1}B_{m-1}
 =\Span\{\ell_i^k,y\ell_i^{k-1},\ldots,y^{m-1}\ell_i^{r+1}\}.
\]
It remains to place these $tm$ polynomials in the monomial basis of $B_k$. To simplify factorials, use the equivalent basis
\[
 Q_{i,a}:=\left.\frac{\dd^a}{\dd c^a}(x+cy)^k\right|_{c=c_i}
 =\frac{k!}{(k-a)!}y^a\ell_i^{k-a},\qquad 0\le a<m.
\]
Index the rows of $M_k$ by $x^{k-j}y^j$ ($0\le j\le k$), and its columns by $(i,a)$. The binomial expansion gives
\begin{equation}\label{human:matrix-entries}
 (M_k)_{j,(i,a)}=\binom kj
 \begin{cases}\dfrac{j!}{(j-a)!}c_i^{j-a},&j\ge a,\\[2pt]0,&j<a.\end{cases}
\end{equation}
An invertible diagonal row scaling gives the same rank as $V_k=(\left.\frac{\dd^a}{\dd T^a}T^j\right|_{T=c_i})$, the confluent Vandermonde matrix in question.

Read its transpose directly. For a vector $(b_0,\ldots,b_k)$, put $q(T)=\sum_{j=0}^k b_jT^j$. Then
\[
 V_k^{\mathsf T}(b_0,\ldots,b_k)^{\mathsf T}
 =\bigl(q^{(a)}(c_i)\bigr)_{i,\,0\le a<m}.
\]
All entries vanish exactly when each $c_i$ is a root of $q$ of multiplicity at least $m$. Since the nodes are distinct, this is equivalent to $\prod_{i=1}^t(T-c_i)^m\mid q$. The transpose kernel therefore has dimension $\pos{k+1-tm}$, proving
\begin{equation}\label{human:matrix-rank}
 \dim I_k=\rank M_k=(k+1)-\pos{k+1-tm}
 =\min\{k+1,t(k-r)\}.
\end{equation}
This gives a direct matrix proof of the local rank formula.

\subsubsection{Second proof: binary apolarity}

Use the binary differential pairing
\[
 \langle b,H\rangle=b(\pd_\xi,\pd_\eta)H\big|_0,
 \qquad b\in B_k,\quad H\in B_k^\vee:=\R[\xi,\eta]_k.
\]
Monomials of the same exponents pair to their factorials, so the pairing is perfect. For $D_i=a_i\pd_\xi+b_i\pd_\eta$,
\[
 \langle\ell_i^{r+1}q,H\rangle=\langle q,D_i^{r+1}H\rangle.
\]
Consequently, $H$ annihilates $I_k$ exactly when $D_i^{r+1}H=0$ for every $i$. Put $\ell_i^\perp=b_i\xi-a_i\eta$, so $D_i\ell_i^\perp=0$. Set $T_i=(a_i\xi+b_i\eta)/(a_i^2+b_i^2)$ and $N_i=\ell_i^\perp$. These are independent linear forms with $D_iT_i=1$ and $D_iN_i=0$. Expand
\[
 H=\sum_{j=0}^k c_jT_i^jN_i^{k-j},\qquad
 D_i^{r+1}H=\sum_{j=r+1}^k c_j\frac{j!}{(j-r-1)!}T_i^{j-r-1}N_i^{k-j}.
\]
The monomials on the right are distinct. Thus the derivative vanishes exactly when $c_j=0$ for every $j\ge r+1$, or equivalently when $N_i^{k-r}$ divides $H$. Distinct $\ell_i^\perp$ are relatively prime. Hence
\begin{equation}\label{human:binary-annihilator}
 I_k^\perp=
 \left(\prod_{i=1}^t(\ell_i^\perp)^{k-r}\right)
 B^\vee_{k-t(k-r)}\qquad(k\ge r+1),
\end{equation}
where a negative-degree space is zero. Taking dimensions proves
\begin{equation}\label{human:binary-rank}
 \dim I_k=\min\{k+1,t(k-r)\},\qquad
 \dim(B/I)_k=\pos{k+1-t(k-r)}.
\end{equation}
For $k\le r$, $I_k=0$ and $\dim(B/I)_k=k+1$. This argument computes the dimension directly from the supporting directions, including configurations with opposite incident rays.

The two calculations independently give the same binary rank. Summing the pieces, for every $d\ge r$ we obtain
\begin{equation}\label{human:local-general}
 \begin{aligned}
 \dim J(v)_d
 &=\sum_{k=r+1}^d\min\{k+1,t(k-r)\}\\
 &=\binom{d+2}{2}-\binom{r+2}{2}
 -\sum_{s=1}^{d-r}\pos{r+1-(t-1)s}.
 \end{aligned}
\end{equation}
The second line uses only $\min\{a,b\}=a-\pos{a-b}$ and $k=r+s$. We now apply this general formula at the degree of the paper.

\begin{lemma}[Local direction contribution]\label{lem:localdim}
At $d=3r+1$,
\begin{equation}\label{eq:localdim}
 \dim(R/J(v))_d=\binom{r+2}{2}+\sigma_v(r).
\end{equation}
\end{lemma}
\begin{proof}
The decomposition \eqref{human:binary-pieces} gives
\[
 \dim(R/J(v))_d=\sum_{k=0}^d\dim(B/I)_k.
\]
The terms with $k\le r$ sum to $\binom{r+2}{2}$. For $k=r+j$, the remaining summand in \eqref{human:binary-rank} is $\pos{r+1-j(t_v-1)}$. Since $t_v\ge2$, it vanishes for $k>2r$, so the finite sum through $d=3r+1$ equals the sum in \eqref{eq:sigma}.
\end{proof}
The binomial identity
\begin{equation}\label{eq:binomialidentity}
 \binom{3r+3}{2}-\binom{r+2}{2}=2\binom{2r+2}{2}
\end{equation}
therefore yields
\begin{equation}\label{human:vertex-row-count}
 \dim J(v)_{3r+1}=2\binom{2r+2}{2}-\sigma_v(r).
\end{equation}
These are the actual local row counts to use when assembling the global matrix.

\subsection{The global endpoint matrix and its row relations}\label{sec:duality}
\subsubsection{From edge cofactors to vertex closure}
Let $s_\tau\in R_d$ be the homogeneous polynomial on a triangle. Across $e$, choose its line form as a transverse coordinate. Equality of all derivatives through order $r$ says that the coefficients of transverse powers $0,\ldots,r$ in the difference vanish. Thus the smoothness condition is exactly
\begin{equation}\label{human:cofactor}
 s_{\tau_+}-s_{\tau_-}=\ell_e^{r+1}q_e,
 \qquad q_e\in R_{d-r-1}.
\end{equation}
At $d=3r+1$, each cofactor belongs to $R_{2r}$ and has $\binom{2r+2}{2}$ coefficients. Multiplication by $\ell_e^{r+1}$ is injective, identifying its input space with $J(e)_d$.

Orient the triangles counterclockwise and choose orientations of the interior edges. The differences \eqref{human:cofactor} telescope around every interior vertex. This gives the signed endpoint map \citep{Billera88}
\begin{equation}\label{eq:idealmap}
 \pd_J:\bigoplus_{e\in\Delta_1^\circ}J(e)_d
 \longrightarrow\bigoplus_{v\in\Delta_0^\circ}J(v)_d,
 \qquad (g_e)_e\longmapsto
       \left(\sum_{e\ni v}\operatorname{sgn}(v,e)g_e\right)_v.
\end{equation}
Here $\operatorname{sgn}(v,e)$ is $+1$ at the head and $-1$ at the tail. Boundary endpoints have no target component. A cofactor family coming from spline pieces lies in $\ker\pd_J$.

Choose an independent basis of each $J(v)_d$ and a basis of each edge input. Denote the resulting matrix by $T_d$. Its number of rows is $m_d=\sum_v\dim J(v)_d$, and at the degree in question it has $E_I\binom{2r+2}{2}$ columns. To distinguish these rows from raw coefficient equations, let $B_v$ embed a chosen basis of $J(v)_d$ into the monomial coefficients of $R_d$. The ambient matrix satisfies
\[
 T_d^{\rm amb}=\diag(B_v)_v\,T_d.
\]
Every $B_v$ has independent columns, so the two matrices have the same rank. The ambient matrix has additional local row relations whose dimensions were already counted in Lemma~\ref{lem:localdim}.

\begin{example}[Adjacent interior vertices with different local dimensions]\label{human:ex:mesh}
Take boundary vertices $A=(0,0)$, $B=(6,0)$, $C=(6,6)$, $D=(0,6)$ and four interior vertices
\[
 u=(2,2),\qquad v=(4,2),\qquad w=(4,4),\qquad p=(2,4).
\]
The ten counterclockwise triangles are
\[
 \begin{gathered}
 ABv,\quad Avu,\quad BCv,\quad Cwv,\quad CDw,\\
 Dpw,\quad DAp,\quad Aup,\quad uvw,\quad uwp.
 \end{gathered}
\]
Their doubled areas are $12,4,12,4,12,4,12,4,4,4$, all positive and summing to $72$, twice the area of the outer square. At $u$ the three supporting lines are $y=x$, $x=2$, and $y=2$. At its neighbor $v$ the five supporting lines are $y=2$, $x=4$, $x+y=6$, $x-2y=0$, and $2x-y=6$. Each of $w,p$ has four supporting lines. Thus, for $r=3,d=10$, Lemma~\ref{lem:localdim} gives
\[
 \dim J(u)_{10}=54,\qquad \dim J(v)_{10}=56,\qquad
 \dim J(w)_{10}=\dim J(p)_{10}=55.
\]
There are $13$ interior edges, so $T_{10}$ has size $220\times364$, while its ambient coefficient form has $264$ rows. On the interior edge $e=uv$, choose $\ell_e=y-2z$ and orient it from $v$ to $u$. The input $q_e\in R_6$ has $28$ coefficients. Then $g_e=(y-2z)^4q_e$ contributes $+g_e$ to the $54$-dimensional row block at $u$, $-g_e$ to the $56$-dimensional row block at $v$, and zero to the other two interior vertex blocks. This exhibits one edge connecting two local output spaces of different dimensions.
\end{example}

\begin{figure}[!b]
\centering
\begin{minipage}[c]{.47\linewidth}
\centering
\begin{tikzpicture}[x=.82cm,y=.82cm,every node/.style={font=\small},line join=round]
 \coordinate(A)at(0,0);\coordinate(B)at(6,0);
 \coordinate(C)at(6,6);\coordinate(D)at(0,6);
 \coordinate(u)at(2,2);\coordinate(v)at(4,2);
 \coordinate(w)at(4,4);\coordinate(p)at(2,4);
 \draw[line width=.7pt](A)--(B)--(C)--(D)--cycle;
 \draw[line width=.45pt](u)--(v)--(w)--(p)--cycle;
 \draw[line width=.45pt](u)--(w);
 \draw[line width=.45pt](A)--(v);\draw[line width=.45pt](A)--(u);
 \draw[line width=.45pt](A)--(p);\draw[line width=.45pt](B)--(v);
 \draw[line width=.45pt](C)--(v);\draw[line width=.45pt](C)--(w);
 \draw[line width=.45pt](D)--(w);\draw[line width=.45pt](D)--(p);
 \foreach\q in {u,v,w,p}{
  \draw[draw=black!65,line width=.55pt,
   preaction={draw=white,line width=2pt,arrows=-},
   -{Stealth[length=1.35mm,width=.95mm]}]
   (\q) ++(15:.31) arc[start angle=15,end angle=345,radius=.31];
 }
 \draw[line width=1.05pt,shorten <=3.8mm,shorten >=3.8mm,
  -{Stealth[length=1.8mm,width=1.2mm]}](v)--(u);
 \foreach\q/\pos in {A/below left,B/below right,C/above right,D/above left}{
  \fill(\q)circle(1.6pt);\node[\pos]at(\q){$\q$};
 }
 \foreach\q in {u,v,w,p}{\fill(\q)circle(1.6pt);}
 \node[fill=white,inner sep=.5pt]at(1.53,1.95){$u$};
 \node[fill=white,inner sep=.5pt]at(4.47,1.88){$v$};
 \node[fill=white,inner sep=.5pt]at(3.71,4.43){$w$};
 \node[fill=white,inner sep=.5pt]at(1.53,4.15){$p$};
\end{tikzpicture}

\smallskip
{\small (a) A mesh with four interior vertices}
\end{minipage}\hfill
\begin{minipage}[c]{.50\linewidth}
\centering
\begin{tikzpicture}[every node/.style={font=\small},>=Stealth]
 \node[draw,align=center,text width=22mm,inner sep=4pt] (input)at(0,0){$q_e\in R_6$\\$28$ coefficients};
 \node[draw,align=center,text width=30mm,inner sep=4pt] (outU)at(3.8,1.2){$+g_e\in J(u)_{10}$\\$54$ coordinates};
 \node[draw,align=center,text width=30mm,inner sep=4pt] (outV)at(3.8,-1.2){$-g_e\in J(v)_{10}$\\$56$ coordinates};
 \draw[->](input.east)--(outU.west);
 \draw[->](input.east)--(outV.west);
 \node[align=center,font=\small]at(1.75,-2.05){$g_e=(y-2z)^4q_e$};
\end{tikzpicture}

\smallskip
{\small (b) The edge input block for $uv$}
\end{minipage}
\caption{The four-interior-vertex mesh of Example~\ref{human:ex:mesh} and its $uv$ input block. Each small circular arrow indicates counterclockwise traversal of the closure condition around its interior vertex. The thick straight arrow only specifies the input orientation $e:v\to u$. The two output row blocks have dimensions $54$ and $56$; the other interior vertex blocks are zero on this edge input.}\label{human:fig:matrix}
\end{figure}
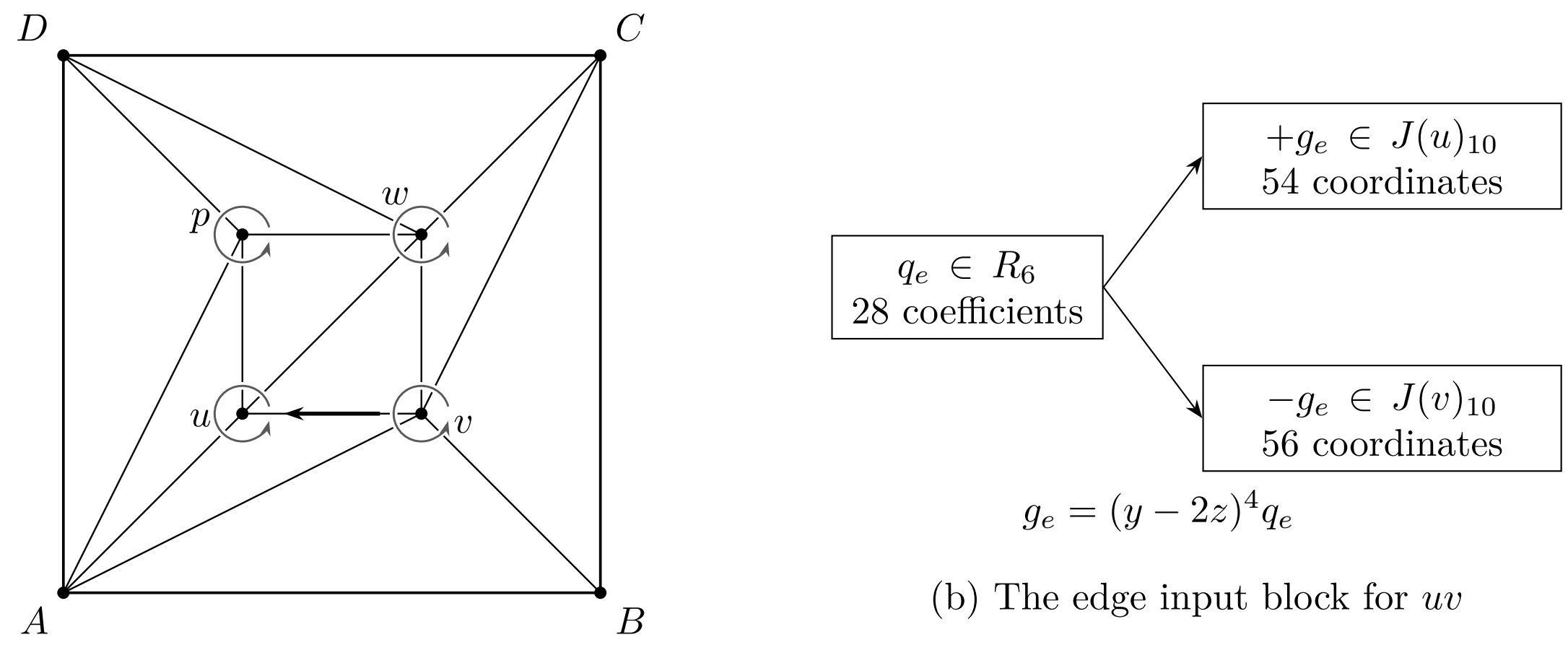

\subsubsection{The remaining dimension term is a cokernel}
To relate the matrix back to splines, we also need to integrate a family in its kernel. For a disk this follows from the relative chain complex, whose coefficients here are in $R_d$.
\begin{proposition}[Dimension identity]\label{prop:dimension}
For a triangulated disk and $d\ge r+1$,
\begin{equation}\label{eq:dimension-identity}
 \dim S_d^r=\binom{d+2}{2}+E_I\binom{d-r+1}{2}
             -\sum_v\dim J(v)_d+\dim M_d,
\end{equation}
where
\begin{equation}\label{eq:defect}
 M_d=\coker\pd_J,\qquad\eps_d^r=\dim M_d.
\end{equation}
\end{proposition}
\begin{proof}
In the relative simplicial chain complex of $(\Delta,\partial\Delta)$, the first homology is zero and the second homology is the common coefficient space $R_d$. The first boundary map restricts to \eqref{eq:idealmap}. Hence every $g\in\ker\pd_J$ is the boundary of a family $(s_\tau)_\tau$ of face polynomials. Membership of its edge differences in $J(e)_d$ gives the $C^r$ conditions \eqref{human:cofactor}. The kernel of this face-difference map consists of one common polynomial on all faces. We therefore have an exact sequence
\[
 0\longrightarrow R_d\longrightarrow S_d^r(\Delta)
       \longrightarrow\ker\pd_J\longrightarrow0.
\]
Rank--nullity and $\dim J(e)_d=\binom{d-r+1}{2}$ give \eqref{eq:dimension-identity}.
\end{proof}
Substituting \eqref{human:vertex-row-count} in \eqref{eq:dimension-identity} shows that Theorem~\ref{thm:main} is equivalent to $M_{3r+1}=0$, or to full row rank of $T_{3r+1}$.

\subsubsection{A row relation is a family of vertex functionals}
Let $\eta$ assign weights to the independent output coordinates of $T_d$. Write the corresponding linear functional on $J(v)_d$ as $\phi_v$. The equality $\eta^{\mathsf T}T_d=0$ says exactly
\begin{equation}\label{eq:functionals}
 \sum_v\phi_v\bigl((\pd_Jg)_v\bigr)=0
                  \qquad\text{for every }g=(g_e)_e.
\end{equation}
Choose an input supported on one edge. If both endpoints $u,v$ are interior, its two signed contributions give $\phi_u(g_e)=\phi_v(g_e)$ for all $g_e\in J(e)_d$. If only one endpoint is interior, its restriction to $J(e)_d$ is zero. This is the first useful consequence of the row relation: it identifies the restrictions of adjacent vertex functionals. The next step is to encode all of these restrictions by polynomial differentiation.

\subsection{Apolarity and point-form coordinates}\label{human:sec:apolar}
\subsubsection{Encode the row weights in polynomials}
The local calculation used a perfect differential pairing to turn multiplication by a line power into differentiation. We apply the same construction in three variables:
\[
 A=\R[X,Y,Z],\qquad
 \langle f,H\rangle=f(\pd_X,\pd_Y,\pd_Z)H\big|_0,
 \qquad f\in R_d,\ H\in A_d.
\]
The notation $X,Y,Z$ refers to the original dual variables throughout this section. For equal total degrees,
\[
 \langle x^iy^jz^k,X^aY^bZ^c\rangle
   =i!j!k!\,\delta_{ia}\delta_{jb}\delta_{kc}.
\]
Thus a functional assigning weights $\eta_{ijk}$ to monomial coefficients of $R_d$ is represented by $\sum\eta_{ijk}X^iY^jZ^k/(i!j!k!)$. This explains how apolarity stores a row relation: the factorials are fixed changes of coordinates, and the output remains the same scalar evaluation.

For $\ell=ax+by+cz$, put $D_\ell=a\pd_X+b\pd_Y+c\pd_Z$. The multiplication--differentiation identity is
\begin{equation}\label{eq:adjoint}
 \langle\ell^{r+1}q,H\rangle=\langle q,D_\ell^{r+1}H\rangle.
\end{equation}
Extend $\phi_v$ from $J(v)_d$ to $R_d$ and represent it by $H_v\in A_d$. Two extensions represent the same functional $\phi_v$ precisely when their polynomial representatives differ by an element of the local annihilator $\A_v$, where
\begin{equation}\label{eq:annihilator}
 \mathcal Z_e=\ker(D_{\ell_e}^{r+1}:A_d\longrightarrow A_{d-r-1}),\qquad \A_v=\bigcap_{e\ni v}\mathcal Z_e.
\end{equation}
Here $\mathcal Z_e$ is the zero-trace space on one edge, while $\A_v$ consists of the polynomials with zero dual trace on every edge incident to $v$. This construction allows the spaces $J(v)_d$ to have different dimensions. We therefore have $J(v)_d^*\simeq A_d/\A_v$. In particular, the vertex functional is zero exactly when its representative lies in $\A_v$.

\begin{lemma}[Edge compatibility]\label{lem:dual-defect}
The dual space $M_d^*$ consists of the classes $(H_v+\A_v)_v$ satisfying
\begin{equation}\label{eq:edge-trace}
 D_{\ell_e}^{r+1}H_u=D_{\ell_e}^{r+1}H_v=:G_e
                         \qquad(e=uv).
\end{equation}
At a boundary endpoint the auxiliary representative is $H_b=0$.
\end{lemma}
\begin{proof}
For $g_e=\ell_e^{r+1}q$, the restrictions obtained from \eqref{eq:functionals} are equal for every $q\in R_{d-r-1}$. By \eqref{eq:adjoint} and the perfect pairing in that degree, their equality is equivalent to \eqref{eq:edge-trace}. A single interior endpoint gives zero on the other side. These are the equations of $\ker\pd_J^*=(\coker\pd_J)^*$.
\end{proof}
We call $G_e=D_{\ell_e}^{r+1}H_v$ the dual edge trace. It is a degree-$2r$ auxiliary polynomial when $d=3r+1$. Boundary edges are assigned zero auxiliary traces. An interior chord with both endpoints on the boundary has a zero endpoint column; assign it the same zero auxiliary trace on its two adjacent faces. Every edge now has a $G_e$, which is zero if an endpoint is on the boundary. These conventions complete the face-edge data used next.

\subsubsection{Make the edge derivatives coordinate derivatives}
To organize three traces on one triangle, we need coordinates adapted to its three vertices. Define the point form
\begin{equation}\label{eq:pointform}
 P_v=x_vX+y_vY+Z.
\end{equation}
It records the actual vertex: $D_\ell P_v=\ell(x_v,y_v,1)$. For $\tau=(u,v,w)$ set $U=P_u,V=P_v,W=P_w$. Then
\[
 \begin{pmatrix}U\\V\\W\end{pmatrix}
 =C_\tau\begin{pmatrix}X\\Y\\Z\end{pmatrix},\qquad
 C_\tau=\begin{pmatrix}x_u&y_u&1\\x_v&y_v&1\\x_w&y_w&1\end{pmatrix}.
\]
The determinant is the nonzero oriented doubled area $[u,v,w]$. For $e=uv$, the line vanishes at $u,v$, whereas $\ell_e(w,1)\ne0$. The chain rule gives
\begin{equation}\label{human:point-derivative}
 D_{\ell_e}U=D_{\ell_e}V=0,\quad
 D_{\ell_e}W=\ell_e(w,1),\qquad
 D_{\ell_e}=\ell_e(w,1)\pd_W.
\end{equation}
The other two edges give nonzero multiples of $\pd_U,\pd_V$.

For completeness, the primal variables dual to these point forms are obtained by the inverse transpose:
\begin{equation}\label{eq:dualcoordinates}
 (\mathsf u,\mathsf v,\mathsf w)^{\mathsf T}
       =C_\tau^{-\mathsf T}(x,y,z)^{\mathsf T}.
\end{equation}
Their cross-pairings with $U,V,W$ vanish and the diagonal pairings equal one; monomials retain the factorial pairing rule. This coordinate relation also gives $\ell_{uv}=\ell_{uv}(w,1)\mathsf w$.

\subsection{A common face polynomial and its residuals}\label{sec:faces}
From now on fix $d=3r+1$. Equation~\eqref{eq:edge-trace} compares derivatives of representatives at neighboring vertices. To extract factors from both endpoints of a shared edge, we seek one polynomial $K_\tau$ on each face that has the same three traces. Then $H_v-K_\tau$ is killed by both derivatives incident to $v$, and subtracting residuals on two faces removes $H_v$ itself.

\subsubsection{Construct the face polynomial coefficient by coefficient}
For an edge $vw$ opposite $u$, divide its dual trace by $\ell_{vw}(u,1)^{r+1}$ and call the result $T_U\in A_{2r}$. Define $T_V,T_W$ cyclically. The required equations are
\[
 \pd_U^{r+1}K_\tau=T_U,\qquad
 \pd_V^{r+1}K_\tau=T_V,\qquad
 \pd_W^{r+1}K_\tau=T_W.
\]
The pair $T_U,T_V$ comes from the same endpoint representative $H_w$, after the same geometric normalizations. Hence
\begin{equation}\label{human:mixed-traces}
 \pd_V^{r+1}T_U=\pd_U^{r+1}T_V,
\end{equation}
and likewise for the other two pairs.
\begin{lemma}[Common face polynomial]\label{lem:facepotential}
Each face $\tau$ has a unique $K_\tau\in A_{3r+1}$ reproducing its three dual edge traces. If it contains a boundary edge or an interior chord with two boundary endpoints, then $K_\tau=0$.
\end{lemma}
\begin{proof}
For $i+j+k=3r+1$, at least one exponent is at least $r+1$. If $i\ge r+1$, recover the coefficient of $U^iV^jW^k$ by
\begin{equation}\label{human:recover-coefficient}
 [U^iV^jW^k]K_\tau
   =\frac{(i-r-1)!}{i!}
      [U^{i-r-1}V^jW^k]T_U.
\end{equation}
If $j$ is also at least $r+1$, \eqref{human:mixed-traces} equates this value with the coefficient recovered from $T_V$, including both factorial factors. Thus these rules consistently construct every coefficient and reproduce the three traces. They also prove uniqueness: a polynomial killed by all three $(r+1)$st derivatives would have all three exponents at most $r$, and therefore degree at most $3r$.

If one edge has two boundary endpoints, all three face traces are zero: each other edge also has a boundary endpoint. The boundary or chord convention supplies the trace on the given edge. Uniqueness then gives $K_\tau=0$.
\end{proof}

\begin{example}[A degree-ten coefficient calculation]\label{human:ex:K}
At $r=3$, consider a face polynomial with the three terms
\[
 K=aU^4V^3W^3+bU^4V^4W^2+cU^3V^3W^4.
\]
Its normalized traces are
\[
 T_U=24aV^3W^3+24bV^4W^2,\quad
 T_V=24bU^4W^2,\quad T_W=24cU^3V^3.
\]
The first and third coefficients are recovered from one trace each. The coefficient $b$ is recovered from both $T_U$ and $T_V$; their compatibility reads $\pd_V^4T_U=\pd_U^4T_V=576bW^2$. This is the overlap handled in \eqref{human:recover-coefficient}.
\end{example}

\subsubsection{The exponent bounds and the radial layers}
At the vertex $v$ of $\tau=(u,v,w)$, the two incident edge derivatives are nonzero multiples of $\pd_U$ and $\pd_W$. Their traces on $H_v$ and $K_\tau$ agree. Hence the $U,W$ exponents in the difference are at most $r$:
\begin{equation}\label{eq:radial-expansion}
 H_v-K_\tau=\sum_{0\le i,j\le r}c_{ij}U^iW^jV^{3r+1-i-j}.
\end{equation}
These $(r+1)^2$ monomials form the intersection of the two incident derivative kernels. The smallest possible $V$ exponent is $r+1$, attained only at $i=j=r$. Therefore
\begin{equation}\label{eq:leading}
 H_v-K_\tau\equiv\lambda_{v,\tau}V^{r+1}U^rW^r
                         \pmod{V^{r+2}A_{2r-1}}.
\end{equation}
We call the coefficient of $V^{r+1+p}$ the $p$th radial layer, with $p=0$ for the leading layer. For $0\le p\le r$, it has the form
\begin{equation}\label{human:radial-layer}
 U^{r-p}W^{r-p}\sum_{j=0}^p c_jU^jW^{p-j}.
\end{equation}
Thus the first layers have $1,2,3,\ldots$ coefficients per face. The term ``radial'' refers to this filtration by the center point form $V$.
In particular, the full remainder has the form
\[
 \begin{aligned}
 H_v-K_\tau&=\lambda_{v,\tau}V^{r+1}U^rW^r+V^{r+2}F_{v,\tau},\\
 F_{v,\tau}&=U^{r-1}W^{r-1}(aU+bW)+VC.
 \end{aligned}
\]
At the $V^{r+2}$ layer the transverse exponents sum to $2r-1$ and each is at most $r$, leaving exactly these two monomials. Every further layer contains one more factor $V$. This proves \eqref{body:F-factor} with the vertex names interchanged.

\subsubsection{Boundary face chains give zero potentials}
\begin{lemma}[Boundary propagation]\label{lem:boundaryfaces}
Every face containing a boundary vertex has $K_\tau=0$.
\end{lemma}
\begin{proof}
Fix a boundary vertex $b$. Its incident faces form an open chain beginning with a face containing a boundary edge. That face has zero potential by Lemma~\ref{lem:facepotential}. Suppose the next face is $(b,u,w)$ and the common edge with a zero face is $bu$. At $u$, both residuals in \eqref{eq:radial-expansion} contain $P_u^{r+1}$. Their difference shows $P_u^{r+1}\mid K_{(b,u,w)}$. The zero trace on $bw$ also gives $\pd_{P_u}^{r+1}K_{(b,u,w)}=0$ in the independent coordinates $(P_b,P_u,P_w)$. Thus its $P_u$ exponent is both at least $r+1$ and at most $r$, so the new potential is zero. This propagates along the chain; a double-boundary chord also has zero potentials on its adjacent faces by Lemma~\ref{lem:facepotential}.
\end{proof}

\subsubsection{Compare the two axes of a shared-edge jump}
Consider $\tau=(u,v,w)$ and $\tau'=(u,v,w')$ on opposite sides of $uv$. Write
\begin{equation}\label{eq:edgechange}
 W'=sW+pU+qV,\qquad s=\frac{[u,v,w']}{[u,v,w]}<0.
\end{equation}
Here $p,q$ are geometric scalars determined by the vertex positions. For actual unscaled point forms, $s+p+q=1$.
Subtract the two residuals at $u$. The representative $H_u$ cancels, so $U^{r+1}$ divides the face-potential difference. The same argument at $v$ gives the relatively prime factor $V^{r+1}$. Hence
\begin{equation}\label{eq:initialjump}
 J:=K_{\tau'}-K_\tau=U^{r+1}V^{r+1}Q,
                        \qquad Q\in A_{r-1}.
\end{equation}
This is the purpose of introducing the common face polynomial: each endpoint contributes a factor to the same jump.

At $u$, divide by $U^{r+1}$ and restrict to $U=0$. Equations~\eqref{eq:leading} and \eqref{eq:initialjump} give
\[
 VQ(0,V,W)=\lambda_{u,\tau}W^r
                  -\lambda_{u,\tau'}(sW+qV)^r.
\]
The constant term in $V$ yields $\lambda_{u,\tau'}=s^{-r}\lambda_{u,\tau}$. At the other endpoint the analogous equation gives the same transport for $\lambda_v$. Dividing the binomial numerators by their displayed variable, the two axis restrictions are
\begin{align}\label{human:jump-axes}
 Q(0,V,W)&=-\frac{q\lambda_{u,\tau}}{s^r}S_r(sW,sW+qV),\nonumber\\
 Q(U,0,W)&=-\frac{p\lambda_{v,\tau}}{s^r}S_r(sW,sW+pU),
\end{align}
where $S_r(a,b)=\sum_{j=0}^{r-1}a^{r-1-j}b^j$ satisfies $b^r-a^r=(b-a)S_r(a,b)$. At $U=V=0$, their equality gives $q\lambda_{u,\tau}=p\lambda_{v,\tau}$. We have proved
\begin{equation}\label{eq:fulltransport}
 \lambda_{u,\tau'}=s^{-r}\lambda_{u,\tau},\quad
 \lambda_{v,\tau'}=s^{-r}\lambda_{v,\tau},\quad
 q\lambda_{u,\tau}=p\lambda_{v,\tau}=:m.
\end{equation}
All these identities retain the branches $p=0$ or $q=0$.

The polynomial $-mS_r(sW,W')/s^r$ realizes both restrictions in \eqref{human:jump-axes}. Their difference from $Q$ vanishes on both axes and is therefore divisible by $UV$. Consequently the complete jump is
\begin{equation}\label{eq:explicitjump}
 J=-\frac{m}{s^r}U^{r+1}V^{r+1}S_r(sW,W')
              +U^{r+2}V^{r+2}R_e,\qquad R_e\in A_{r-3}.
\end{equation}
We set $A_k=0$ for $k<0$. At $r=3$, the displayed remainder is already a scalar, but it is still part of the full compatibility conditions. More generally, after the leading layer vanishes, the jump may still have the form $U^{r+2}V^{r+2}R_e$ of residual degree $r-3$. The higher-layer induction will raise the endpoint powers further.

It remains to eliminate the leading and higher residual coefficients. Theorem~\ref{thm:induction} gives the uniform induction for the higher radial coefficients, and Theorem~\ref{thm:from-leading} completes the defect argument by comparing degrees.

\section{Proofs of leading-layer transport and odd-order leading-layer elimination}\label{app:odd}\label{sec:corner}

We give the calculations supporting Section~\ref{body:odd} in their dependency order. The geometric figure and the logical outline remain in the main text.

Our goal is to choose representatives for which the coefficient of $P_v^{r+1}$ in every residual $H_v-K_\tau$ vanishes. Leading terms at two-direction vertices can first be removed modulo the local annihilators; the transport and compatibility relations below deal with the remaining coefficients. The local calculations apply to all $r\ge3$. They eliminate the leading layer at odd orders; even orders also require Appendix~\ref{app:even}.

\subsection{Two computations of the same edge jump}
Let faces $(u,v,w)$ and $(u,v,z)$ share $uv$. Put $U=P_u,V=P_v,W=P_w$, and write
\[
 P_z=sW+pU+qV,\qquad s<0.
\]
The negative sign comes from the third vertices lying on opposite sides of the edge. The endpoint residuals make the face jump divisible by
$U^{r+1}V^{r+1}$. Comparing the leading terms at $v$ on $V=0$ gives
\begin{equation}\label{eq:lambda-transport}
 \lambda_{v,(u,v,z)}=s^{-r}\lambda_{v,(u,v,w)}.
\end{equation}
After division by $U^{r+1}V^{r+1}$, the restriction of this same quotient to $V=0$ is
\[
 \lambda_v\frac{W^r-(W+(p/s)U)^r}{U}.
\]
At $U=0$ its value is $-r\lambda_v(p/s)W^{r-1}$.
The computation at $u$ gives $-r\lambda_u(q/s)W^{r-1}$. Consequently,
\begin{equation}\label{eq:lambda-coupling}
 p\lambda_{v,(u,v,w)}=q\lambda_{u,(u,v,w)}.
\end{equation}
These are restrictions of polynomials, so the resulting identity also includes $p=0$ or $q=0$.

For an ordered face $\tau=(a,b,c)$, write
$\Delta_\tau=\det(b-a,c-a)$ for its signed doubled area in the fixed physical chart.
Fix the representative at the center $v$. First choose an initial ordered face $\tau_0=(v,u_0,u_1)$, and define
\begin{equation}\label{eq:vertex-normalization}
 \Lambda_v:=\lambda_{v,\tau_0}\Delta_{\tau_0}^r.
\end{equation}

Next prove invariance across adjacent faces. For the two faces above, $\tau=(u,v,w)$ and $\tau'=(u,v,z)$, the common order of $u,v$ gives
\[
 \Delta_{\tau'}=s\Delta_\tau,\qquad
 s=\frac{\Delta_{(u,v,z)}}{\Delta_{(u,v,w)}}<0.
\]
The opposite signs of these areas already account for the sign of $s$.
Applying \eqref{eq:lambda-transport} at either endpoint therefore gives
\begin{equation}\label{eq:area-normalized-transport}
 \lambda_{a,\tau'}\Delta_{\tau'}^r=\lambda_{a,\tau}\Delta_\tau^r,\qquad a\in\{u,v\}.
\end{equation}
The same simultaneous permutation of the two ordered faces preserves the area ratio, so the law applies to any matching slots for the shared vertices. Starting from $\tau_0$, continue around $v$ by keeping the center in the first slot and the shared neighbor in its current slot, replacing only the other neighbor. For example, the orders begin
\[
 (v,u_0,u_1),\quad (v,u_2,u_1),\quad (v,u_2,u_3),\quad\ldots.
\]
Equation~\eqref{eq:area-normalized-transport} then gives at each step
\[
 \lambda_{v,\tau'}\Delta_{\tau'}^r
 =\lambda_{v,\tau}\Delta_\tau^r=\Lambda_v.
\]
This proves invariance along the face chain continued from the initial face. Every signed area uses that continued order, with no additional sign factor in the definition.

Finally check closure around the entire star. Let $N$ be the valence of $v$. Each crossing reverses the signed area. After $N$ crossings the underlying face returns with continued area $(-1)^N\Delta_{\tau_0}$; its neighbor slots are exchanged if $N$ is odd. Both neighbor forms have exponent $r$ in the leading term, so the face coefficient remains $\lambda_{v,\tau_0}$. Closure consequently requires
\[
 \Lambda_v=(-1)^{rN}\Lambda_v.
\]
For even $r$, including $r=2$, signed-area normalization agrees with ordinary positive-area normalization and closure has no sign obstruction. For odd $r$, a nonzero leading layer requires even $N$, so odd-valence vertices have $\Lambda_v=0$. The local normalization is uniform, and the parity restriction is read from closure. Reversing the initial order multiplies $\Lambda_v$ by $(-1)^r$; its vanishing is independent of this choice. The local coefficient $s$ remains in the calculations below.

Number the neighbors counterclockwise as $u_0,\ldots,u_{N-1}$ and write $\tau_i=\{v,u_i,u_{i+1}\}$ and $K_i=K_{\tau_i}$ for the underlying faces and their potentials, with indices modulo $N$. For binary restrictions put $L_i=P_{u_i}\bmod(P_v)$ in the coordinates induced by the fixed physical chart. The determinants of consecutive coefficient vectors are the positive doubled areas in this counterclockwise order.

\subsection{Normalize the removable terms and choose a support extreme point}
A complete interior star on two supporting lines has four rays. Choose actual neighbor forms $U,W$ on the two lines. Then
\[
 G=V^{r+1}U^rW^r\in\A_v.
\]
The first line derivative kills $V,U$ and only differentiates $W$, whose exponent is $r$; the other line is treated in the same way.
The leading terms of $G$ on all four faces are nonzero and satisfy \eqref{eq:area-normalized-transport}.
Transport determines all coefficients from one, so a unique multiple of $G$ removes this vertex's leading layer.

If a third supporting direction is present, its derivative is $a\pd_U+b\pd_W$ with $ab\ne0$. Expanding
the expression
\[
 (a\pd_U+b\pd_W)^{r+1}(U^rW^r)
\]
gives nonzero terms with distinct monomials. Every incident derivative kills $V$, so different powers of $V$ remain independent.
Hence elements of $\A_v$ have zero leading coefficient. We now perform the preceding normalization at every two-direction vertex.
At the boundary, Lemma~\ref{lem:boundaryfaces} and $H_b=0$ already give zero leading terms.
Define the normalized leading support by
\[
 \Sigma=\{v:\lambda_{v,\tau}\ne0\text{ for some incident face }\tau\}.
\]
Transport makes ``some face'' equivalent to ``every incident face.''
If the finite set $\Sigma$ is nonempty, choose a linear functional taking distinct values on its vertices
and let $v$ be its unique maximizer. All other support points lie in a strict open half-plane through $v$.
This constructs the support extreme point used in the next lemma.

\begin{lemma}[The standard form at a leading-support extreme point]\label{lem:corner-shape}
If the normalized leading support is nonempty, the neighbor cycle at a support extreme point consists of $n+1$ consecutive nonzero-support neighbors and
exactly two zero-support neighbors, with $n\ge2$. Orientation-preserving affine coordinates may be chosen so that the rays are
\begin{equation}\label{eq:corner-rays}
 (1,0)=\mathbf r_0,\mathbf r_1,\ldots,\mathbf r_n=(0,1),
 \quad(-b,0),(0,-a),\qquad a,b>0,
\end{equation}
where $\mathbf r_i=\rho_i(1-t_i,t_i)$, $\rho_i>0$,
$0=t_0<\cdots<t_n=1$, and $\rho_0=\rho_n=1$.
A configuration of this form is called a leading corner, or simply a corner.
\end{lemma}
\begin{proof}
All nonzero-support neighbors lie in one strict open half-plane through the center.
By \eqref{eq:lambda-coupling}, a zero neighbor has opposite predecessor and successor rays.
If all neighbors were zero, repeated use of $L_{i+1}=-c_iL_{i-1}$, $c_i>0$,
would repeat a directed ray or give a four-ray two-direction star, already normalized to zero.
Thus a nonzero neighbor exists. A zero neighbor also exists, since a complete interior star fills a neighborhood of its center.

An isolated zero neighbor would have two opposite nonzero neighbors, contrary to their strict half-plane.
Three consecutive zero neighbors force rays four positions apart to have the same direction; embeddedness excludes this except for the normalized two-direction four-star.
Every maximal zero block therefore consists of two vertices. Write its four consecutive bordering rays as $u,z_1,z_2,w$, where $z_1,z_2$ are zero-support neighbors and $u,w$ are nonzero-support neighbors. Let their three consecutive angular gaps be $\theta_1,\theta_2,\theta_3\in(0,\pi)$. At $z_1$, the predecessor $u$ and successor $z_2$ are opposite; at $z_2$, the predecessor $z_1$ and successor $w$ are opposite. Thus
\[
 \theta_1+\theta_2=\pi,\qquad \theta_2+\theta_3=\pi.
\]
The arc through this zero block and its two bordering rays has angular length
\begin{equation}\label{eq:zero-block-span}
 \theta_1+\theta_2+\theta_3=2\pi-\theta_2>\pi.
\end{equation}
Distinct maximal zero blocks are separated by nonzero neighbors. Their corresponding arcs have disjoint interiors, even if they share a bordering nonzero ray. Two such arcs would have total length greater than $2\pi$, which is impossible in one complete neighbor cycle. There is therefore exactly one zero block. A nonzero block of length one gives consecutive opposite rays; length two again gives a two-direction four-star.
The nonzero block thus has at least three vertices. Its extreme rays are independent; using their actual lengths as coordinate units proves the claim.
\end{proof}

\subsection{The leading coefficient in standard coordinates}
In the configuration of Lemma~\ref{lem:corner-shape}, denote the center by $u$ and its neighbors by $p_i$. Let $E=(p_0-u\ \ p_n-u)$ and $D=\det E>0$. Then
\begin{equation}\label{eq:physicalchange}
 \binom\alpha\beta=E^{-1}\left(\binom xy-\binom{x_u}{y_u}\right),
 \qquad \binom xy=\binom{x_u}{y_u}+E\binom\alpha\beta.
\end{equation}
Writing the columns of $E$ as $e_0,e_n$, set
\[
 X_0=e_{0x}X_{\rm old}+e_{0y}Y_{\rm old},\qquad
 Y_0=e_{nx}X_{\rm old}+e_{ny}Y_{\rm old},\qquad U=P_u.
\]
The point form is $P_a=U+\alpha_aX_0+\beta_aY_0$. With
$\widehat X=X_0-U/b$, $\widehat Y=Y_0-U/a$, the exterior forms become $-b\widehat X,-a\widehat Y$. The remaining normalized forms are
\begin{equation}\label{eq:point-normalized}
 \begin{aligned}
 L_i&=P_{p_i}/\rho_i=(1-t_i)\widehat X+t_i\widehat Y+\zeta_iU,\\
 \zeta_i&=\rho_i^{-1}+(1-t_i)/b+t_i/a>0.
 \end{aligned}
\end{equation}
For the exterior vertices take $\rho_{n+1}=-b,\rho_{n+2}=-a$, so that $L_{n+1}=\widehat X,L_{n+2}=\widehat Y$. The complete cyclic data are
\begin{gather}\label{eq:exterior}
 (t_{n+1},\zeta_{n+1})=(0,0),\quad(t_{n+2},\zeta_{n+2})=(1,0),\\
 h_i=t_{i+1}-t_i>0\ (0\le i<n),\quad\sum_{i=0}^{n-1}h_i=1,
 \quad h_n=-1,\ h_{n+1}=1,\ h_{n+2}=-1.\nonumber
\end{gather}
The positive doubled area is $\Delta_i^*=D\rho_i\rho_{i+1}h_i$. Indices are modulo $N=n+3$; the negative exterior scales and the signs of all three exterior intervals are retained.

Choose the initial ordered face $(u,p_0,p_1)$ counterclockwise. Counterclockwise face orders differ from the continued orders by $i$ reversals; put $\sigma_i=(-1)^{ri}$. If $\lambda_i^*$ denotes the coefficient in the actual neighbor forms, \eqref{eq:vertex-normalization} gives
\[
 (\Delta_i^*)^r\lambda_i^*=\Lambda_u\sigma_i=D^r\Lambda\sigma_i,
 \qquad \Lambda=\Lambda_u/D^r\ne0.
\]
In the normalized forms $L_i$, the residual is
\begin{equation}\label{eq:leading-corner}
 \frac{H_u-K_i}{\Lambda}
 =\frac{\sigma_i}{h_i^r}U^{r+1}L_i^rL_{i+1}^r+U^{r+2}F_i,
 \qquad F_i\in A_{2r-1}.
\end{equation}
This equation divides the residual by the common scalar $\Lambda$; $H_u$ and $K_i$ remain fixed.

\subsection{Hermite endpoint data and face-difference coefficients}
The two transverse differential conditions give
\[
 F_i=L_i^{r-1}L_{i+1}^{r-1}(a_iL_i+b_iL_{i+1})+UG_i.
\]
Indeed, the transverse exponents of $F_i\bmod U$ sum to $2r-1$ and are individually at most $r$, leaving only $(r,r-1)$ and $(r-1,r)$. In $(\widehat X,\widehat Y,U)$, define
$f_i(t)=F_i(-t,1-t,0)$. Since $L_j(-t,1-t,0)=t_j-t$, the two endpoints have zeros of multiplicity at least $r-1$. The endpoint coefficients on the two sides of $t_i$ are independently defined by
\begin{equation}\label{human:separate-endpoint-data}
 \begin{aligned}
 A_i^+&:=[x^{r-1}]f_i(t_i+x)=(-1)^{r-1}b_ih_i^r,\\
 A_i^-&:=[x^{r-1}]f_{i-1}(t_i+x)=-a_{i-1}h_{i-1}^r.
 \end{aligned}
\end{equation}
Put $M_i=\sigma_iL_i^rL_{i+1}^r/h_i^r$ and $J_i=(K_i-K_{i-1})/\Lambda$. Subtraction gives
\begin{equation}\label{human:exact-first-subtraction}
 U^{r+2}(F_i-F_{i-1})=-J_i-U^{r+1}(M_i-M_{i-1}).
\end{equation}
Insert the two endpoint axis restrictions \eqref{eq:axis1}--\eqref{eq:axis2} of $J_i$ to obtain the three coefficients in \eqref{eq:fdiff1}--\eqref{eq:fdiff3}. Appendix~\ref{app:coefficients} gives the binomial expansion and the remainder-order estimate.

Equation~\eqref{eq:fdiff1} gives $A_i^+=A_i^-=:A_i$, so we may define $\beta_i=A_i/\sigma_{i-1}$. Using $(-1)^{r-1}\sigma_i=-\sigma_{i-1}$ in \eqref{human:separate-endpoint-data} gives
\begin{equation}\label{human:ab-beta}
 a_i=-\frac{\sigma_i}{h_i^r}\beta_{i+1},\qquad
 b_i=-\frac{\sigma_i}{h_i^r}\beta_i.
\end{equation}
Substitute \eqref{human:ab-beta} into the Hermite reconstruction \eqref{eq:hermite-polynomial} and compare with the next two face-difference coefficients. This yields the following two equations.

\begin{lemma}[Hermite identities]\label{lem:hermite}
For $r\ge2$, at every node including the two exterior nodes,
\begin{align}
 r\delta_i&=\frac{\beta_{i-1}+r\beta_i}{h_{i-1}}
             +\frac{r\beta_i+\beta_{i+1}}{h_i},\label{eq:hermite1}\\
 \frac r2\delta_i(1/h_i-1/h_{i-1})
 &=\frac r2\beta_i(1/h_i^2-1/h_{i-1}^2)
    +\frac{\beta_{i+1}}{h_i^2}-\frac{\beta_{i-1}}{h_{i-1}^2}.\label{eq:hermite2}
\end{align}
They imply
\begin{align}
 \frac{\beta_{i+1}}{h_i}&=\frac{\beta_{i-1}}{h_{i-1}}
                   &&(0\le i\le n),\label{eq:betarec}\\
 \beta_{n+1}&=\beta_0-r\zeta_0,\qquad
 \beta_{n+2}=\beta_n-r\zeta_n,\label{eq:betaend}\\
 \sum_{i=0}^{n-1}\frac{\beta_i+\beta_{i+1}}{h_i}
                   &=\beta_0+\beta_n.\label{eq:flux}
\end{align}
\end{lemma}
\begin{proof}
The first two identities follow from the coefficient comparison in Appendix~\ref{app:coefficients}; the second cancels $r-1\ne0$.
Subtract $(1/h_i-1/h_{i-1})/2$ times the first equation from the second and rearrange:
\[
 \frac{h_{i-1}+h_i}{2h_{i-1}h_i}
 \left(\frac{\beta_{i+1}}{h_i}-\frac{\beta_{i-1}}{h_{i-1}}\right)=0.
\]
The first factor is nonzero on the block, using $h_0,h_{n-1}<1$ at its ends. This gives the recurrence.
At the exterior nodes the interval sum is zero; substitution into the first equation, without division by that sum, gives \eqref{eq:betaend}.
Finally sum around the complete star: $\sum\delta_i=0$, and each interval contributes
$(r+1)(\beta_i+\beta_{i+1})/h_i$.
The three exterior contributions are
$-(\beta_n+\beta_{n+1})+(\beta_{n+1}+\beta_{n+2})-(\beta_{n+2}+\beta_0)
=-\beta_0-\beta_n$, proving \eqref{eq:flux}.
\end{proof}

\begin{proposition}[Excluding odd-length corners]\label{prop:oddcorner}
A nonzero leading corner cannot have odd $n$. In particular, odd $r\ge3$ admits no normalized nonzero leading support.
\end{proposition}
\begin{proof}
If $n$ is odd, the recurrence gives positive constants $c_1,c_2$ such that $\beta_n=c_1\beta_1$ and
$\beta_{n-1}=c_2\beta_0$. The endpoint recurrence and exterior equations give
\[
 \beta_1=\frac{rh_0\zeta_n}{1+h_0c_1}>0,\qquad
 \beta_0=\frac{r\zeta_0}{1+c_2/h_{n-1}}>0.
\]
Both positive-multiplier chains are therefore positive. Since $0<h_i<1$, the left side of the cyclic sum identity strictly exceeds
$\beta_0+\beta_n+2\sum_{i=1}^{n-1}\beta_i$, contradicting \eqref{eq:flux}.

For odd $r$, a nonzero leading layer requires even $N=n+3$, hence odd $n$, already excluded.
\end{proof}

\begin{lemma}[Eliminating the leading coefficient at odd smoothness orders]\label{lem:odd-leading}
For each odd integer $r\ge3$, we may replace the representatives $H_v$ of any dual defect class by equivalent representatives modulo $\A_v$ so that
\[
 H_v-K_\tau\in P_v^{r+2}A_{2r-1}
\]
at every incident vertex and face. Auxiliary boundary representatives remain zero.
\end{lemma}
\begin{proof}
Normalize the two-direction vertices as above. If the leading support remains nonempty, take a support extreme point. Odd $r$ forces its valence $N$ to be even. Its shape gives $N=n+3$, hence odd $n$, contradicting Proposition~\ref{prop:oddcorner}. Thus the normalized leading coefficient vanishes everywhere. The normalization used only elements of the local annihilators, so the functional classes and the actual face potentials have not changed.
\end{proof}

After the leading coefficient has been eliminated, a shared-edge difference may still have a quotient of degree $r-3$. The uniform induction in Theorem~\ref{thm:induction} raises the endpoint powers until the degree argument in Theorem~\ref{thm:from-leading} forces every face difference to vanish.

\section{Elimination of the higher radial coefficients}\label{app:higher}

\subsection{Successive elimination of radial coefficients}\label{sec:higher}

Let $r\ge3$ and suppose each residual is divisible by $P_v^{r+2}$. A shared-edge difference can still have a quotient of degree $r-3$. Without a parity assumption, we eliminate the lowest remaining residual coefficient successively and thereby increase both endpoint powers.

\subsubsection{The lowest remaining residual coefficient}
Fix integers $p\ge1$ and $r\ge2p+1$. Suppose representatives of the same functional classes have already been chosen so that, simultaneously at every vertex and incident face,
\begin{equation}\label{eq:layerhyp}
 H_v-K_\tau\in P_v^{r+1+p}A_{2r-p}.
\end{equation}
We prove that the representatives can be chosen so that $P_v^{r+2+p}$ divides every residual. We first extract the coefficient of the lowest remaining power of $P_v$.

In the coordinates $U=P_u,V=P_v,W=P_w$ of a face $(u,v,w)$, equation~\eqref{eq:radial-expansion} still bounds the $U,W$ exponents by $r$. The coefficient of $V^{r+1+p}$ has total degree $2r-p$, so both transverse exponents are also at least $r-p$. There are exactly $p+1$ possible coefficients:
\begin{align}\label{eq:higherlayer}
 H_v-K_{(u,v,w)}\equiv{}&
 V^{r+1+p}U^{r-p}W^{r-p}
       \sum_{j=0}^p c_{v,\tau,j}U^jW^{p-j}
 \pmod{V^{r+2+p}A_{2r-p-1}}.
\end{align}
These are precisely the coefficients of the lowest remaining power of $P_v$ in the original residual. The two transverse exponents sum to $2r-p$, and each is at most $r$.

On an interior edge $e=uv$, the induction hypothesis at its two endpoints makes the face-potential difference divisible by $P_u^{r+1+p}$ and $P_v^{r+1+p}$. These forms are coprime, hence
\begin{equation}\label{eq:fulljump}
 K_{\tau'}-K_\tau=P_u^{r+1+p}P_v^{r+1+p}Q_e,
 \qquad Q_e\in A_{r-2p-1}.
\end{equation}
We retain this three-variable quotient $Q_e$ to connect the computations at the two endpoints. If the edge has a boundary endpoint, both neighboring face potentials are zero by Lemma~\ref{lem:boundaryfaces}.

\subsubsection{From the current coefficient to transfer around a complete star}
To remove the current term in \eqref{eq:higherlayer}, we need not first solve for every higher term. Dividing the residual by $P_v^{r+1+p}$ and then reducing modulo $P_v$ retains exactly the coefficient in question. The quotient $A/(P_v)$ has two independent polynomial coordinates, which we use below.

List the complete neighbor cycle $u_0,\ldots,u_{N-1}$ counterclockwise and put
\[
 L_i=P_{u_i}\bmod(P_v),\qquad
 \Delta_i=\det(L_i,L_{i+1})>0,\qquad
 \Gamma_i=\det(L_{i-1},L_{i+1}).
\]
The determinants act on the coefficient vectors of binary linear forms. Thus
\begin{equation}\label{eq:starspq}
 L_{i+1}=s_iL_{i-1}+q_iL_i,\qquad
 s_i=-\frac{\Delta_i}{\Delta_{i-1}}<0,\qquad
 q_i=\frac{\Gamma_i}{\Delta_{i-1}}.
\end{equation}
Distinct incident edges have distinct directed rays, and each consecutive angular gap belongs to $(0,\pi)$. Consequently, $q_i=0$ means exactly that the predecessor and successor rays are opposite.

Define the coefficient polynomial
\[
 F_i=\left[(H_v-K_{(v,u_i,u_{i+1})})/P_v^{r+1+p}\right]\bmod(P_v).
\]
Equation~\eqref{eq:higherlayer} gives its factorization within a face. Dividing \eqref{eq:fulljump} by the central factor and reducing modulo that point form gives the difference between neighboring coefficients. They therefore satisfy the system
\begin{equation}\label{eq:starsystem}
 \begin{aligned}
 F_i&=L_i^{r-p}L_{i+1}^{r-p}
       \sum_{j=0}^p a_{i,j}L_i^jL_{i+1}^{p-j},\\
 F_i-F_{i-1}&\in L_i^{r+1+p}\R[X,Y]_{r-2p-1}.
 \end{aligned}
\end{equation}
We call the collection of coefficient polynomials $(F_i)_i$ a local state on the star; equivalently, it is given by the coefficient vectors $(a_{i,0},\ldots,a_{i,p})$. We first classify the local solutions of \eqref{eq:starsystem} and then determine which can occur on the whole triangulation. A mode refers to one of the solution-space types in that classification.

\begin{lemma}[Complete coefficient transfer]\label{lem:coefficienttransfer}
At node $i$, write $U=L_i,W=L_{i-1},L_{i+1}=sW+qU$. The system \eqref{eq:starsystem} at this node is equivalent to
\begin{equation}\label{eq:fullcoefficient}
 s^{r-k}\sum_{j=0}^{\min(p,k)}\binom{r-j}{k-j}q^{k-j}a_{i,j}
 =\begin{cases}a_{i-1,p-k},&0\le k\le p,\\0,&p<k\le2p.\end{cases}
\end{equation}
The first $p+1$ equations define an invertible transfer. If $q=0$, they further imply $F_i=F_{i-1}$.
\end{lemma}
\begin{proof}
Expand the two faces in the same $U,W$ coordinates. After canceling their common factor $U^{r-p}$, the new and old polynomials are respectively
\begin{equation}\label{human:higher-reduced-pair}
 \sum_{j=0}^p a_{i,j}U^j(sW+qU)^{r-j},\qquad
 \sum_{j=0}^p a_{i-1,j}U^{p-j}W^{r-p+j}.
\end{equation}
The original difference must contain $U^{r+1+p}$, so the reduced difference must contain $U^{2p+1}$. Every coefficient from $U^0$ through $U^{2p}$ must therefore match. The first two equations, for example, are
\[
 s^ra_{i,0}=a_{i-1,p},\qquad
 s^{r-1}(rq\,a_{i,0}+a_{i,1})=a_{i-1,p-1}.
\]
The binomial expansion gives \eqref{eq:fullcoefficient} in general. Its first $p+1$ equations are triangular in the new coefficients, with nonzero diagonal entries $s^{r-k}$. Including the reversal of the old vector, the old-to-new transfer has determinant
\begin{equation}\label{eq:transferdet}
 (-1)^{p(p+1)/2}s^{-(p+1)r+p(p+1)/2}\ne0.
\end{equation}
The last $p$ equations are additional compatibility conditions. When $q=0$, both reduced polynomials in \eqref{human:higher-reduced-pair} have $U$-degree at most $p$, so matching the first $p+1$ coefficients already makes them identical.
\end{proof}

Invertible transfer has an immediate consequence: a nonzero solution on a complete star is nonzero on every face. This will justify the nonzero scales needed in subsequent coefficient ratios.

\begin{lemma}[Boundary states]\label{lem:boundary-chain}
At each stage of the induction, the local state at every boundary vertex is zero.
\end{lemma}
\begin{proof}
At a boundary vertex $b$, we have $H_b=0$ and $K_\tau=0$ on every incident face by Lemma~\ref{lem:boundaryfaces}. The residual itself is therefore zero. All changes of representatives below are made at interior vertices; the boundary representatives remain $H_b=0$.
\end{proof}

\subsubsection{Turning the additional coefficient equations into interpolation}
The first $p+1$ equations already govern transfer; what remains is to solve the last $p$ zero conditions. For $q\ne0$, their convolution coefficients satisfy
\[
 \binom{r-j}{k-j}=\binom rk\frac{\binom kj}{\binom rj}.
\]
This identity suggests removing the geometric powers and the fixed binomial factors first. Set
\begin{equation}\label{human:higher-newton-data}
 z_j=\frac{a_{i,j}}{q^j\binom rj},\qquad
 h(x)=\sum_{j=0}^p z_j\binom{x}{j}.
\end{equation}
Here $x$ is an interpolation variable and $\binom{x}{j}$ is a polynomial of degree $j$. The left side of \eqref{eq:fullcoefficient} becomes $s^{r-k}q^k\binom rk h(k)$. Since $r\ge2p+1$, all factors being removed are nonzero, and the last $p$ equations become exactly
\[
 h(p+1)=\cdots=h(2p)=0.
\]
A polynomial of degree at most $p$ with these $p$ distinct roots has only one scale left. Record it by $\lambda=h(0)=a_{i,0}$. Then
\begin{equation}\label{eq:Newtonpolynomial}
 h(x)=\lambda\prod_{\nu=p+1}^{2p}\frac{\nu-x}{\nu}
 =\lambda\frac{\binom{2p-x}p}{\binom{2p}p}.
\end{equation}
If $\lambda=0$, the additional root at $x=0$ forces $h=0$. A nonzero state thus has $\lambda\ne0$.

We can now expand this root polynomial back into the Newton basis with which we started. Its $j$th forward difference at zero is $(-1)^j\binom{2p-j}{p-j}/\binom{2p}p$, giving
\begin{equation}\label{eq:Newtonidentity}
 \frac{\binom{2p-x}p}{\binom{2p}p}
 =\sum_{j=0}^p(-1)^j\frac{\binom pj}{\binom{2p}j}\binom{x}{j}.
\end{equation}
Comparing the coefficients in \eqref{human:higher-newton-data} yields the following normal form.

\begin{proposition}[Binomial normal form and the first unmatched coefficient]\label{prop:normalform}
When $q\ne0$, the additional equations force
\begin{equation}\label{eq:bj}
 a_{i,j}=\lambda(-q)^jb_j(r,p),\qquad
 b_j(r,p)=\frac{\binom rj\binom pj}{\binom{2p}j}>0.
\end{equation}
The first $p+1$ transfer equations give
\begin{equation}\label{eq:oldcoeff}
 a_{i-1,p-k}=\lambda s^{r-k}q^k b_k(r,p)
 \qquad(0\le k\le p).
\end{equation}
The jump coefficient immediately beyond the matching range is
\begin{equation}\label{eq:firstremainder}
 \left.\frac{F_i-F_{i-1}}{U^{r+1+p}}\right|_{U=0}
 =\lambda(-1)^p\frac{\binom r{2p+1}}{\binom{2p}p}
        s^{r-2p-1}q^{2p+1}W^{r-2p-1}.
\end{equation}
\end{proposition}
\begin{proof}
Equation~\eqref{eq:bj} has been derived by interpolation. For $0\le k\le p$, we have $\binom rk h(k)=\lambda b_k$; substitution into the full coefficient equation gives \eqref{eq:oldcoeff}.
To detect an actual jump, read precisely the next coefficient, $k=2p+1$. Equation~\eqref{eq:Newtonpolynomial} gives
\[
 h(2p+1)=\lambda\frac{(-1)^p}{\binom{2p}p}.
\]
The old reduced polynomial has $U$-degree at most $p$ and makes no contribution here. The new coefficient is $s^{r-2p-1}q^{2p+1}\binom r{2p+1}h(2p+1)$. Restoring the common $U^{r-p}$ gives precisely the exponent $r+1+p$, proving \eqref{eq:firstremainder}. The condition $\binom r{2p+1}\ne0$ explains the range needed for the later nonzero test.
\end{proof}

For example, at $p=1$ the two original face coefficients are compressed by the additional equation to $\lambda(1,-rq/2)$. This is still a potentially nonzero face state. Its description at the other endpoint must agree, and transfer around the entire cycle must close.

\subsection{Comparing the two ends of a face and classifying complete stars}\label{sec:stars}

\subsubsection{Eliminating the face scale before testing geometric closure}
Suppose the star state is nonzero and two adjacent nodes satisfy $q_iq_{i+1}\ne0$. At node $i$, equation~\eqref{eq:bj} describes face $i$; at node $i+1$, this face is the old face, so \eqref{eq:oldcoeff} describes the same coefficients a second way:
\[
 a_{i,j}=\lambda_i(-q_i)^jb_j
 =\lambda_{i+1}s_{i+1}^{r-p+j}q_{i+1}^{p-j}b_{p-j}.
\]
All face states are nonzero, hence both scales are nonzero. Taking the first two coefficients as the ratio $a_{i,1}/a_{i,0}$ eliminates those scales and leaves the geometric parameters to the first power:
\[
 -q_i b_1=\frac{s_{i+1}}{q_{i+1}}\frac{b_{p-1}}{b_p}.
\]
Using $b_1=r/2$ and $b_{p-1}/b_p=p(p+1)/(r-p+1)$, we obtain
\begin{equation}\label{eq:uniform-star-ratio}
 \frac{\Gamma_i\Gamma_{i+1}}{\Delta_{i-1}\Delta_{i+1}}
 =\frac{q_iq_{i+1}}{-s_{i+1}}
 =C_p(r):=\frac{2p(p+1)}{r(r-p+1)}.
\end{equation}
This constant is rational and belongs to $(0,1)$: at the smallest allowed value $r=2p+1$, the denominator exceeds the numerator by $3p+2$, and thereafter the denominator increases. We next test whether a complete embedded star can satisfy this fixed ratio at all nodes where it applies.

The following geometric lemma applies to any rational constant $C\in(0,1)$ and converts the ratio condition into a restriction on the directions of a complete star.

\begin{lemma}[A geometric obstruction on complete ray stars]\label{rthree:geometric-obstruction}
Let $L_0,\ldots,L_{N-1}$ form a complete counterclockwise embedded ray star with consecutive angular gaps in $(0,\pi)$, and put $\Delta_i=\det(L_i,L_{i+1})$ and $\Gamma_i=\det(L_{i-1},L_{i+1})$, with indices modulo $N$.
Suppose a rational number $C\in(0,1)$ satisfies
\[
 \frac{\Gamma_i\Gamma_{i+1}}{\Delta_{i-1}\Delta_{i+1}}=C
\]
at every node with $\Gamma_i\Gamma_{i+1}\ne0$.
Then $N=4$, and at least one pair of rays is opposite.
\end{lemma}
\begin{proof}
We first exclude all other valences and direction configurations. The binary determinant identity gives
\begin{equation}\label{eq:plucker}
 \Gamma_i\Gamma_{i+1}
 =\Delta_{i-1}\Delta_{i+1}
   +\det(L_{i-1},L_{i+2})\Delta_i.
\end{equation}
If $N\ge6$, the average sum of three consecutive angular gaps is $6\pi/N\le\pi$.
Some four consecutive rays therefore span at most a half-turn.
Their two middle $\Gamma$ values are strictly positive, while the last term of \eqref{eq:plucker} is nonnegative. The specified ratio is consequently at least one, contradicting $C<1$.
If $N=3$, cyclic indices give $\Gamma_i=-\Delta_{i+1}$, and the same ratio equals one, again a contradiction.

For $N=5$, the existence of four consecutive rays spanning at most a half-turn has already been excluded.
Otherwise every three consecutive angular gaps have sum greater than $\pi$, so their complementary pair has sum less than $\pi$. All $\Gamma_i$ are therefore positive.
Take indices modulo five. Since $\det(L_{i-1},L_{i+2})=-\Gamma_{i+3}$, equation~\eqref{eq:plucker} and the specified ratio give
\[
 \Gamma_{i+3}\Delta_i=(1-C)\Delta_{i-1}\Delta_{i+1}.
\]
Multiply these five equations, and separately multiply the five ratio equations, to obtain
\begin{equation}\label{human:higher-five-products}
 \prod_i\Gamma_i=(1-C)^5\prod_i\Delta_i,\qquad
 (\prod_i\Gamma_i)^2=C^5(\prod_i\Delta_i)^2.
\end{equation}
All products are positive. Cancellation gives $C^5=(1-C)^{10}$. Taking the positive fifth root gives $C=(1-C)^2$, or
\begin{equation}\label{eq:pentagon}
 C^2-3C+1=0.
\end{equation}
Its root in $(0,1)$ is $(3-\sqrt5)/2$, contradicting the rationality of $C$.

It remains to consider $N=4$. If the four rays lie on four distinct supporting lines, all $\Gamma_i$ are nonzero and $\Gamma_{i+2}=-\Gamma_i$.
The specified ratio requires every adjacent product to be positive, forcing all $\Gamma_i$ to have the same sign, a contradiction.
An opposite pair must therefore exist.
\end{proof}

\begin{samepage}
\begin{theorem}[Classification of complete higher-layer stars]\label{thm:starclassification}
For $p\ge1$ and $r\ge2p+1$, the solution spaces of \eqref{eq:starsystem} on a complete counterclockwise embedded ray star are as follows.
\begin{center}
\begin{tabular}{@{}p{.56\textwidth}p{.35\textwidth}@{}}
\toprule
Complete star & Solution space\\\midrule
Three rays, at least five rays, or four rays on four supporting lines & Zero\\[2pt]
Four rays on three supporting lines, with one opposite pair & One-dimensional; jumps only on that pair\\[2pt]
Four rays on two supporting lines, with two opposite pairs & $(p+1)$-dimensional; one polynomial on all faces\\\bottomrule
\end{tabular}
\end{center}
\end{theorem}
\end{samepage}
\begin{proof}
For every adjacent pair with nonzero $q$, equation~\eqref{eq:uniform-star-ratio} has the same rational constant $C_p(r)\in(0,1)$. Lemma~\ref{rthree:geometric-obstruction} therefore leaves only four rays with at least one opposite pair. This excludes the first line of the table.

Two four-ray cases remain. With two opposite pairs, every $q_i=0$, so full transfer makes all face states one and the same polynomial. In coordinates on the two directions, the common polynomial belongs exactly to
\begin{equation}\label{eq:twodirectionkernel}
 U^{r-p}W^{r-p}\R[U,W]_p.
\end{equation}
Every polynomial in this space has the required face factorizations and zero jumps; the dimension is exactly $p+1$.

With only one opposite pair, the nonzero $q$ values occur precisely on that pair. One nonzero $q$ bounds the solution dimension by one, and invertible transfer does not increase it. We must still exhibit a nonzero solution that actually closes. After an orientation-preserving linear change and positive rescaling of the rays, take
\[
 (L_0,L_1,L_2,L_3)=(X,Y,-X,c_0X-Y),\qquad c_0\ne0.
\]
Then $s_i=-1$ and $(q_0,q_1,q_2,q_3)=(c_0,0,-c_0,0)$. Set $b_j=b_j(r,p)$ and choose
\begin{equation}\label{eq:fourraymode}
 \begin{aligned}
 a_{0,j}&=(-c_0)^jb_j,&a_{1,j}&=(-1)^{r-j}a_{0,p-j},\\
 a_{2,j}&=(-1)^pc_0^jb_j,&a_{3,j}&=(-1)^{r-j}a_{2,p-j}.
 \end{aligned}
\end{equation}
At $q_1=q_3=0$, the second and fourth groups are exactly the triangular transfer equations, giving $F_1=F_0$ and $F_3=F_2$. At $q_2=-c_0$, the new coefficients have normal-form scale $(-1)^p$. The old low coefficients required there are
\[
 (-1)^p(-1)^{r-k}(-c_0)^kb_k
 =(-1)^{r-p}c_0^kb_k=a_{1,p-k}.
\]
The last $p$ zero equations follow from the normal form. At the final node $q_0=c_0$, the new coefficients have scale one, while the old low coefficients are indeed
\[
 a_{3,p-k}=(-1)^{r-p+k}(-1)^pc_0^kb_k
          =(-1)^{r-k}c_0^kb_k.
\]
These are exactly the equations required to return to the first group; its high zero conditions also hold. The entire cycle has now been checked. Since $a_{0,0}=b_0=1$, the solution is nonzero. Equation~\eqref{eq:firstremainder} makes the jumps at both nonzero $q$ nodes nonzero. The solution space is therefore exactly one-dimensional, with jumps precisely on the opposite pair.
\end{proof}

The two remaining spaces require different treatment. In the two-direction mode, the coefficient polynomial $F_i$ is the same on every incident face. After lifting this common polynomial to the local annihilator below, we will remove it by changing the vertex representative. The three-direction mode carries actual face differences and cannot be discarded as a freedom of representation. We first determine whether it can continue to transmit across the triangulation.

\subsubsection{Connecting the endpoints by the shared-edge quotient}
The tangential coordinates and local scales at distinct vertices differ. To compare transmission at the two ends, return to the quotient $Q_e$ fixed in \eqref{eq:fulljump}. If $\xi_e$ is the coefficient vector of the edge-line equation, its two endpoints satisfy $P_u(\xi_e)=P_v(\xi_e)=0$, while the third point of a nondegenerate neighboring triangle satisfies $P_w(\xi_e)\ne0$. Evaluation in this common direction reads the lowest jump at both endpoints.

\begin{lemma}[Shared-edge evaluation]\label{lem:shared-edge}
Fix the order of the two neighboring faces at an interior edge $e=uv$ and define
\begin{equation}\label{eq:edgescalar}
 c_e=Q_e(\xi_e).
\end{equation}
Its vanishing is an intrinsic edge property. A nonzero three-direction four-ray mode has exactly two nonzero incident $c_e$, on its opposite pair; a zero mode or a two-direction common mode has all incident $c_e$ equal to zero.
\end{lemma}
\begin{proof}
At endpoint $v$, use the face coordinates $U=P_u,V=P_v,W=P_w$, keeping the same symbol for $Q_e$ expressed in these coordinates. Divide the face difference by $V^{r+1+p}$ and reduce modulo $V$. This gives
\[
 F_i-F_{i-1}=\varepsilon_{v,e}U^{r+1+p}Q_e(U,0,W),
 \qquad \varepsilon_{v,e}\in\{1,-1\},
\]
where the sign records the difference between the fixed global side order and the local cyclic order. Therefore
\begin{equation}\label{eq:edge-evaluation-interface}
 c_e=\varepsilon_{v,e}
 \left.\frac{F_i-F_{i-1}}{U^{r+1+p}}
 \right|_{U=0,\,W=P_w(\xi_e)}.
\end{equation}
Repeating this operation at the other endpoint evaluates the same three-variable $Q_e$. Reversing the side order changes its sign; rescaling the line coefficient vector multiplies its value by a nonzero factor. Vanishing is unchanged.

If $q_i\ne0$ and the mode is nonzero, equation~\eqref{eq:firstremainder} gives
\[
 c_e=\varepsilon_{v,e}\lambda(-1)^p
 \frac{\binom r{2p+1}}{\binom{2p}p}
 s^{r-2p-1}q^{2p+1}P_w(\xi_e)^{r-2p-1}\ne0.
\]
Every nonzero factor has been justified: the state scale is nonzero, $s\ne0$, $q\ne0$, $r\ge2p+1$, and the third vertex is off the shared edge line. If $q_i=0$, full polynomial equality $F_i=F_{i-1}$ gives $c_e=0$. The classification completes the proof. This also includes the limiting degree $r=2p+1$, when $Q_e$ is constant.
\end{proof}

We need not infer from $c_e=0$ that an arbitrary $Q_e$ is identically zero: the evaluation already detects every possible remaining one-dimensional transmission mode. An interior edge with a boundary endpoint has both neighboring $K$ equal to zero, hence $Q_e=c_e=0$.

\begin{lemma}[Vanishing of finite straight transmission]\label{lem:straight}
Every $c_e$ vanishes, and all three-direction four-ray modes vanish.
\end{lemma}
\begin{proof}
Suppose a nonzero edge exists and let $G$ be the finite graph of actual edges with $c_e\ne0$ and their endpoints. No endpoint is on the boundary. Classification and shared-edge evaluation show that each graph vertex has exactly two graph edges, with opposite rays: transmission can neither stop nor turn at an interior vertex.

Choose a linear function that is nonconstant on each graph edge; only finitely many directions must be avoided. At its maximal graph vertex, one of the two opposite edges must lead to a larger value, a contradiction. Thus $G$ is empty, and the detection formula makes every three-direction mode zero.
\end{proof}

\subsubsection{Lifting the common state to the local annihilator}
After transmission is excluded, a nonzero current state can only be a common polynomial from \eqref{eq:twodirectionkernel}. This is a statement in the binary quotient modulo $P_v$; it does not yet permit a change of the original $H_v$. At a two-direction vertex, choose actual neighbor forms $U=P_u,W=P_w$ on the two supporting lines and lift the common term to
\begin{equation}\label{eq:gaugehigher}
 G_v=P_v^{r+1+p}U^{r-p}W^{r-p}
       \sum_{j=0}^p g_jU^jW^{p-j}.
\end{equation}
Each term has total degree $3r+1$, and its $U,W$ exponents are respectively $r-p+j\le r$ and $r-j\le r$. The edge derivative on the first supporting line kills $P_v,U$ and acts only on $W$; a derivative on the other line acts only on $U$. Every incident $(r+1)$st derivative thus kills $G_v$, proving $G_v\in\A_v$.

We may now replace $H_v$ by $H_v-G_v$. This removes exactly the current tangential coefficient. Since $G_v$ contains $P_v^{r+1+p}$, no lower radial layer is restored. All actual edge traces remain unchanged, and uniqueness then preserves all face potentials. These replacements are made independently at interior vertices. Boundary representatives are not changed and remain $H_b=0$.

\begin{theorem}[Uniform higher-layer induction]\label{thm:induction}
For $p\ge1$ and $r\ge2p+1$, each representative $H_v$ satisfying \eqref{eq:layerhyp} may be replaced by an equivalent representative modulo $\A_v$ so that
\begin{equation}\label{eq:inductionstep}
 H_v-K_\tau\in P_v^{r+2+p}A_{2r-p-1}
\end{equation}
holds simultaneously at all vertices and incident faces.
\end{theorem}
\begin{proof}
The classification exhausts the local states. Lemma~\ref{lem:straight} eliminates the three-direction transmission modes. At each interior two-direction vertex, subtracting the local annihilator element in \eqref{eq:gaugehigher} from its representative removes the remaining common coefficient. Lemma~\ref{lem:boundary-chain} gives zero boundary states. Thus the current coefficient vanishes everywhere, increasing the point-form exponent by one without changing the functional classes.
\end{proof}

\begin{remark}[Induction range]\label{rem:layer-range}
The step uses $r>2p$. At the adjacent value $r=2p$, we have $C_p(r)=1$, and the common polynomial $(L_0L_1L_2)^p$ on three rays has the required face factorizations and zero jumps; also $\binom r{2p+1}=0$. The degree termination below occurs before that layer needs to be treated.
\end{remark}

\subsection{From leading-term elimination to vanishing of the defect}\label{sec:completion}
\begin{theorem}[From leading-coefficient elimination to zero defect]\label{thm:from-leading}
Fix $r\ge3$. If a dual defect class admits representatives with $H_v-K_\tau\in P_v^{r+2}A_{2r-1}$ at every incident vertex and face, that class is zero.
\end{theorem}
\begin{proof}
The assumed divisibility is exactly the induction hypothesis for $p=1$. Apply the step successively for
\[
 p=1,2,\ldots,\left\lfloor\frac{r-1}{2}\right\rfloor.
\]
Every application satisfies $r\ge2p+1$ and preserves the same face potentials. Finally,
\begin{equation}\label{eq:final-depth}
 H_v-K_\tau\in P_v^{r+1+s_*}A_{2r-s_*},\qquad
 s_*=\left\lfloor\frac{r+1}{2}\right\rfloor
     =\left\lceil\frac r2\right\rceil.
\end{equation}
We do not need to remove every residual coefficient: this divisibility depth already suffices. On an edge with two interior endpoints, the face jump contains $P_u^{r+1+s_*}P_v^{r+1+s_*}$, so its quotient would have degree
\begin{equation}\label{eq:negative-degree}
 (3r+1)-2(r+1+s_*)=
 \begin{cases}-1,&r\text{ even},\\-2,&r\text{ odd}.
 \end{cases}
\end{equation}
The jump is therefore zero. An edge with a boundary endpoint already has zero neighboring potentials. Hence all adjacent face potentials agree; face adjacency is connected and there is a zero boundary face, so every $K_\tau=0$.

Since each potential reproduces the actual edge traces, $D_{\ell_e}^{r+1}H_v=0$ for every $e\ni v$. Thus $H_v\in\A_v$ and the original vertex functional class is zero. 
\end{proof}

\begin{theorem}[The dimension formula for every odd order]\label{thm:oddcase}
Theorem~\ref{thm:main} holds for every odd integer $r\ge3$.
\end{theorem}
\begin{proof}
Lemma~\ref{lem:odd-leading} supplies the hypothesis of Theorem~\ref{thm:from-leading} for every dual defect class. Thus every class is zero and $M_{3r+1}=0$. Proposition~\ref{prop:dimension}, Lemma~\ref{lem:localdim}, and \eqref{eq:binomialidentity} give \eqref{eq:main}.
\end{proof}

This completes the odd orders. After eliminating the leading coefficient at even orders, we use the same Theorem~\ref{thm:from-leading}.

\section{The even-order leading-layer proof}\label{app:even}

\subsection{Second radial coefficients and a positive cyclic system for even orders}\label{sec:positivity}
The proof for all odd orders is now complete. It remains to eliminate the leading coefficient for even $r\ge4$; the higher-layer argument is already available. Proposition~\ref{prop:oddcorner} leaves only even $n\ge2$. We first remove its shortest possibility.

\begin{lemma}[The shortest even corner]\label{lem:short-corner}
A nonzero leading corner, as defined in Lemma~\ref{lem:corner-shape}, cannot have $n=2$.
\end{lemma}
\begin{proof}
Suppose $n=2$, and write $a=h_0$, $b=h_1$, so that $a,b>0$ and $a+b=1$. The recurrence \eqref{eq:betarec} gives
\[
 \beta_2=\frac ba\beta_0,\qquad \beta_3=-\frac{\beta_1}{b}.
\]
Substituting the first identity into the cyclic sum identity \eqref{eq:flux} yields
\[
 \frac{\beta_0}{a}+\frac{\beta_1}{ab}=0.
\]
Hence $\beta_1=-b\beta_0$ and $\beta_3=\beta_0$. Equation~\eqref{eq:betaend} then implies $r\zeta_0=0$, contradicting $\zeta_0>0$.
\end{proof}

Henceforth $r\ge4$ and $n\ge4$ are both even.
The recurrence and cyclic sum identity force alternating signs in the positive block, but we still need the next layer of the same residual.
That layer will supply positive scalar states and then force scalar reflection.

\subsubsection{Compatibility of the next radial coefficients}
Use the center form $U=P_u$, the shared neighbor form $A=L_i$, and the previous and next neighbor forms $W=L_{i-1}$, $T=L_{i+1}$. The first comparison gives
\[
 \frac{H_u-K_i}{\Lambda}=U^{r+1}M_i+U^{r+2}\Pi_i+U^{r+3}G_i,
\]
where $M_i=h_i^{-r}L_i^rL_{i+1}^r$ and
$\Pi_i=-h_i^{-r}L_i^{r-1}L_{i+1}^{r-1}(\beta_{i+1}L_i+\beta_iL_{i+1})$.
The transverse exponents in each face are at most $r$. Initially distinguish its endpoint coefficients by writing
\[
 G_i\equiv h_i^{-r}L_i^{r-2}L_{i+1}^{r-2}
 (\gamma_{i+1}^-L_i^2+\eta_iL_iL_{i+1}+\gamma_i^+L_{i+1}^2)
 \pmod{UA_{2r-3}}.
\]
Subtracting the two neighboring expressions gives
\begin{equation}\label{human:second-exact-difference}
 \begin{aligned}
 U^{r+3}(G_i-G_{i-1})={}&-\frac{K_i-K_{i-1}}{\Lambda}
 -U^{r+1}(M_i-M_{i-1})\\
 &-U^{r+2}(\Pi_i-\Pi_{i-1}).
 \end{aligned}
\end{equation}
First divide \eqref{human:second-exact-difference} by $U^{r+1}$. Insert the endpoint axis restrictions \eqref{eq:axis1}, \eqref{eq:axis2}, extract the coefficient of $U^2$, and remove $A^{r-2}$ as in \eqref{eq:second-0}--\eqref{eq:second-2} of Appendix~\ref{app:coefficients}. Constant-coefficient matching gives $\gamma_i^-=\gamma_i^+$; denote this value by $\gamma_i$.
Put $C_i=h_i^{-r}$ and $E_i=(H_u-K_i)/(\Lambda U^{r+1})$. Then
\begin{align}\label{eq:jet2}
 E_i=C_i\{&L_i^rL_{i+1}^r
 -UL_i^{r-1}L_{i+1}^{r-1}(\beta_{i+1}L_i+\beta_iL_{i+1})\nonumber\\
 &+U^2L_i^{r-2}L_{i+1}^{r-2}
 (\gamma_{i+1}L_i^2+\eta_iL_iL_{i+1}+\gamma_iL_{i+1}^2)\}
 \pmod{U^3A_{2r-3}}.
\end{align}
The common coordinates satisfy
\begin{equation}\label{eq:spq}
 \begin{gathered}
 s_i=-h_i/h_{i-1},\qquad q_i=1+h_i/h_{i-1},\qquad \xi_i=h_i\delta_i,\\
 T=s_iW+q_iA+\xi_iU.
 \end{gathered}
\end{equation}
Because the neighbor forms have been rescaled, $s_i$ need not be negative; at the two exterior nodes, $s_i=1$. In the unindexed local formulas, $s,q,p$ denote $s_i,q_i,\xi_i$, respectively.
The linear, quadratic, and cubic coefficients are expanded in \eqref{eq:rawjet1}--\eqref{eq:rawjet3}. Appendix~\ref{app:coefficients} also shows why the undetermined remainder starts at $A^4$ and affects none of them. We record the resulting transfers here.

\begin{lemma}[Necessary equations of the second radial layer]\label{lem:jet2}
At every node,
\begin{align}
 \eta_i&=s_i\eta_{i-1}-rq_i\gamma_i+r\beta_i\xi_i,\label{eq:eta}\\
 \gamma_{i+1}&=s_i^2\gamma_{i-1}-(r-1)s_iq_i\eta_{i-1}
 +\frac{r(r-1)}2q_i^2\gamma_i\nonumber\\
 &\quad +(r-1)\beta_{i+1}\xi_i-\frac{r(r-1)}2\xi_i^2.\label{eq:gamma}
\end{align}
For $0\le i\le n$, also
\begin{align}\label{eq:thirdjet}
 0={}&6\gamma_{i+1}+3(r-1)q_i\eta_i+r(r-1)q_i^2\gamma_i
 -3r(r-1)\beta_i\xi_iq_i\nonumber\\
 &-6(r-1)\beta_{i+1}\xi_i+2r(r-1)\xi_i^2.
\end{align}
\end{lemma}
\begin{proof}
Set \eqref{eq:rawjet1} and \eqref{eq:rawjet2} to zero and eliminate $\eta_i$ to obtain the first two equations. The third raw coefficient contains $q_i(r-2)$.
On the block $q_i\ne0$ and $r\ge4$, so division gives \eqref{eq:thirdjet}.
At the two exterior nodes $q_i=0$, and the first two equations remain. The same factor vanishes at $r=2$.
Thus the actual second radial layer supplies an additional condition for the present even-order argument.
\end{proof}

\subsubsection{Extracting linear states from a quadratic particular solution}
First establish that the subsequent divisions are legitimate. The recurrence expresses even coefficients as positive multiples of $\beta_0$ and odd coefficients as positive multiples of $\beta_1$.
Move the right-hand side of the cyclic sum identity \eqref{eq:flux} to the left: the endpoint weights are $1/h-1>0$, and the internal weights are positive.
Hence $\beta_0,\beta_1$ have opposite signs. If one were zero, all would vanish, and the exterior equations would give $\zeta_0=0$.
All block coefficients are consequently nonzero with alternating signs.
Multiplying the recurrence by $\beta_i$ also gives
\[
 \frac{\beta_i\beta_{i+1}}{h_i}
 =\frac{\beta_{i-1}\beta_i}{h_{i-1}}=-c_*^2,
 \qquad c_*>0.
\]
To absorb this common product constant into the radial unit, rewrite the coordinate and coefficients as
\[
 \widehat U=c_*U,\quad\widehat\Lambda=\Lambda/c_*^{r+1},\quad
 (\widehat\zeta_i,\widehat\beta_i)=(\zeta_i,\beta_i)/c_*,\quad
 (\widehat\gamma_i,\widehat\eta_i)=(\gamma_i,\eta_i)/c_*^2.
\]
The actual $H_u,K_i,L_i$ remain fixed, and \eqref{eq:leading-corner} and \eqref{eq:jet2} retain their form.
Hereafter $U$ denotes a positive multiple of the original central point form. Dropping hats gives
\begin{equation}\label{eq:product-normalization}
 h_i=-\beta_i\beta_{i+1}\quad(0\le i\le n\text{ or }i=n+2),\qquad
 \beta_{n+1}=\beta_n^{-1},\quad\beta_{n+2}=\beta_0^{-1}.
\end{equation}
The block identities follow from the product constant; the three exterior identities follow from the endpoint recurrence and $h_n=h_{n+2}=-1$.

The raw second-radial-layer equations contain known quadratic terms. First the Hermite equation gives
\[
 \beta_{i+1}=-s_i\beta_{i-1},\qquad
 \xi_i=h_i\delta_i=q_i\beta_i+\frac2r\beta_{i+1}\quad(0\le i\le n).
\]
This suggests testing $\gamma_i^0=\beta_i^2$ and $\eta_i^0=\beta_i\beta_{i+1}$.
The residual in the first transfer is
\[
 \beta_i[\beta_{i+1}-s_i\beta_{i-1}+rq_i\beta_i-r\xi_i]=0.
\]
In the second, cancel $\beta_{i+1}^2=s_i^2\beta_{i-1}^2$; what remains is
\[
 (r-1)\left[q_i\beta_i\beta_{i+1}+\frac r2q_i^2\beta_i^2
 +\beta_{i+1}\xi_i-\frac r2\xi_i^2\right]=0,
\]
which also cancels term by term after inserting $\xi_i$.
Subtract this known particular solution and define
\begin{equation}\label{human:shifted-states}
 \widetilde\gamma_i=\gamma_i-\beta_i^2,\quad e_i=\eta_i-\beta_i\beta_{i+1},\qquad
 \mathsf u_i=\widetilde\gamma_i/\beta_i^2,\quad
 \mathsf v_i=e_i/(\beta_i\beta_{i+1}).
\end{equation}
The scalar $\mathsf u_i$ is attached to node $i$, and $\mathsf v_i$ to interval $i$ between the directions $t_i,t_{i+1}$. They are distinct from the geometric vertices and point forms.
The first two transfers now read
\[
 e_i-s_ie_{i-1}=-rq_i\widetilde\gamma_i,\qquad
 \widetilde\gamma_{i+1}=s_i^2\widetilde\gamma_{i-1}-(r-1)s_iq_ie_{i-1}+\frac{r(r-1)}2q_i^2\widetilde\gamma_i.
\]
Divide the first by $\beta_i\beta_{i+1}$; its left side becomes $\mathsf v_i+\mathsf v_{i-1}$.
Define the positive parameters determined by the direction intervals and normalized node data:
\begin{equation}\label{eq:directions}
 d_i=-q_i\beta_i/\beta_{i+1}>0\qquad(0\le i\le n).
\end{equation}
Its internal positivity follows from alternating signs. At the ends use $q_0=1-h_0>0$,
$q_n=1-1/h_{n-1}<0$, and the exterior reciprocal relations.
Division of the second equation by $\beta_{i+1}^2$ gives
$\mathsf u_{i+1}-\mathsf u_{i-1}=-(r-1)d_i\mathsf v_{i-1}
+r(r-1)d_i^2\mathsf u_i/2$. Combining this with the first equation yields
\begin{equation}\label{eq:twotransfers}
 \begin{aligned}
 \mathsf v_i+\mathsf v_{i-1}&=rd_i\mathsf u_i,\\
 \mathsf u_{i+1}-\mathsf u_{i-1}&=\frac{r-1}{2}d_i(\mathsf v_i-\mathsf v_{i-1}).
 \end{aligned}
\end{equation}

The third coefficient equation retains a nonzero forcing term. Insert \eqref{human:shifted-states} and divide by $\beta_{i+1}^2$.
The unknown and known parts become, respectively,
\[
 6\mathsf u_{i+1}-3(r-1)d_i\mathsf v_i+r(r-1)d_i^2\mathsf u_i,
 \qquad \frac{2(r+2)}r+(r-1)d_i.
\]
The two transfers reconstruct
\[
 \mathsf v_i=\frac r2d_i\mathsf u_i+
 \frac{\mathsf u_{i+1}-\mathsf u_{i-1}}{(r-1)d_i}.
\]
After substitution, the unknown part becomes $3(\mathsf u_{i-1}+\mathsf u_{i+1})
-r(r-1)d_i^2\mathsf u_i/2$. Abbreviate the coefficients by
\begin{equation}\label{eq:constants}
 \alpha=\frac{r(r-1)}6,\qquad c=\frac{2(r+2)}{3r},\qquad
 \omega=\frac{r-1}3,\qquad\chi=\frac{r+1}{2r}.
\end{equation}
The first three constants give the scalar equation
\begin{equation}\label{eq:scalareq}
 \mathsf u_{i-1}+\mathsf u_{i+1}-\alpha d_i^2\mathsf u_i=-c-\omega d_i
 \qquad(0\le i\le n).
\end{equation}
The constant $\chi$ comes from the exterior nodes. There $s=1,q=0$, and $p=(\beta_{i+1}-\beta_{i-1})/r$.
The first two raw equations give
\[
 e_i=e_{i-1},\qquad
 \widetilde\gamma_{i+1}-\widetilde\gamma_{i-1}=-\frac{r+1}{2r}(\beta_{i+1}^2-\beta_{i-1}^2).
\]
Denote the multiplier determined by the endpoints by
\begin{equation}\label{eq:period}
 \rho=\beta_0\beta_n>0,\qquad \mathsf R=\rho^2.
\end{equation}
For example, divide $\widetilde\gamma_{n+2}=\widetilde\gamma_n-\chi(\beta_{n+2}^2-\beta_n^2)$ by
$\beta_{n+2}^2=\beta_0^{-2}$ to obtain the first equation below. The other endpoint gives the second, and the exterior $e$ states agree:
\begin{equation}\label{eq:exterior-u}
 \begin{aligned}
 \mathsf u_{-1}&=\mathsf R\mathsf u_n+\chi(\mathsf R-1),\\
 \mathsf u_{n+1}&=\mathsf R\mathsf u_0+\chi(\mathsf R-1),\\
 \mathsf v_{-1}&=\mathsf v_n.
 \end{aligned}
\end{equation}
Also \eqref{eq:betaend} gives $r\zeta_0=(\rho-1)/\beta_n$, so $\rho\ne1$ and $\mathsf R\ne1$.

\subsubsection{Reading positivity from a cyclic linear system}
The inward extension will require $\mathsf u_i>0$. To prove it, first express the direction coefficients in interval lengths.
Using \eqref{eq:product-normalization} gives
\begin{align}
 d_id_{i+1}&=\frac{(h_{i-1}+h_i)(h_i+h_{i+1})}{h_{i-1}h_{i+1}}>1
                      &&(1\le i\le n-2),\label{eq:di-products}\\
 d_0d_1&=\frac{(1-h_0)(h_0+h_1)}{h_1}>1,\label{eq:end-products}\\
 d_{n-1}d_n&=\frac{(1-h_{n-1})(h_{n-2}+h_{n-1})}{h_{n-2}}>1,\nonumber\\
 \frac{\mathsf R}{d_0d_n}&=\frac{h_0h_{n-1}}{(1-h_0)(1-h_{n-1})}<1.\label{eq:wrap-weight}
\end{align}
For an internal product minus one, the numerator is $h_i(h_{i-1}+h_i+h_{i+1})>0$.
At the left endpoint the difference is $h_0(1-h_0-h_1)/h_1>0$, and the right endpoint is analogous.
The last bound uses $h_0+h_{n-1}<1$. These inequalities are strict because $n\ge4$.

\begin{lemma}[Positivity of actual extensions]\label{lem:positive}
The actual states above satisfy
\begin{equation}\label{eq:positive-u}
 \mathsf u_i>0\ (0\le i\le n),\qquad
 \mathsf u_{-1}>-\chi,\quad\mathsf u_{n+1}>-\chi.
\end{equation}
\end{lemma}
\begin{proof}
Put $Y_i=d_i\mathsf u_i$, insert the exterior conditions into the scalar equation, and divide row $i$ by $-d_i$.
The diagonal becomes the constant $\alpha$, and adjacent coefficients are symmetric negative edge weights:
\[
 MY=f,\quad M_{ii}=\alpha,\quad
 M_{i,i+1}=M_{i+1,i}=-\frac1{d_id_{i+1}},\quad
 M_{0n}=M_{n0}=-\frac{\mathsf R}{d_0d_n}.
\]
All remaining entries vanish. The right side is
\[
 f_i=\frac c{d_i}+\omega+
 \frac{\chi(\mathsf R-1)}{d_i}\mathbf1_{\{0,n\}}(i).
\]
The internal right sides are positive. At the endpoints they equal $(c-\chi+\chi\mathsf R)/d_i+\omega>0$, since
$c-\chi=(r+5)/(6r)>0$.
Each row has two off-diagonal weights less than one, while $\alpha\ge2$.
Write $M=\alpha I-W$. Then $W\ge0$ and $\|W\|_\infty<\alpha$, so
\[
 M^{-1}=\frac1\alpha\sum_{j=0}^{\infty}(W/\alpha)^j
\]
converges entrywise nonnegatively and has positive diagonal. Since $f>0$, we obtain $Y_i>0$ and $\mathsf u_i>0$.
The two exterior bounds follow from \eqref{eq:exterior-u} and $\mathsf R>0$.
\end{proof}

\subsubsection{From cyclic endpoint data to forced reflection}
Positivity is now available. To compare inward extensions from the two sides, we first need equal current states at the endpoints.
Multiply the two exterior relations by $\mathsf u_0$ and $\mathsf u_n$, respectively. Their known common product cancels, giving
\[
 \mathsf u_{-1}\mathsf u_0-\mathsf u_n\mathsf u_{n+1}
 =\chi(\mathsf R-1)(\mathsf u_0-\mathsf u_n).
\]
Thus equality of the adjacent products on the left will imply equal endpoint states because $\mathsf R\ne1$.
The exterior interval states also satisfy $\mathsf v_{-1}=\mathsf v_n$. This suggests comparing the change of an adjacent vertex product with the change of the interval state.

Compute the same product difference using the two transfer equations:
\begin{align*}
 \mathsf u_i\mathsf u_{i+1}-\mathsf u_{i-1}\mathsf u_i
 &=\frac{r-1}{2}d_i\mathsf u_i(\mathsf v_i-\mathsf v_{i-1})\\
 &=\frac{r-1}{2r}(\mathsf v_i+\mathsf v_{i-1})
                         (\mathsf v_i-\mathsf v_{i-1})\\
 &=\frac{r-1}{2r}(\mathsf v_i^2-\mathsf v_{i-1}^2).
\end{align*}
The right side is a difference of squares. Moving it to the left produces the invariant
\begin{equation}\label{eq:invariant}
 I_i=\mathsf u_i\mathsf u_{i+1}-\frac{r-1}{2r}\mathsf v_i^2,
 \qquad I_i=I_{i-1}.
\end{equation}
Hence $I_{-1}=I_n$. The equal exterior interval states cancel their square terms and give
$\mathsf u_{-1}\mathsf u_0=\mathsf u_n\mathsf u_{n+1}$, as required.
Together with the exterior relation above, this yields
$\chi(\mathsf R-1)(\mathsf u_0-\mathsf u_n)=0$, and therefore $\mathsf u_0=\mathsf u_n$.

The current endpoint states are now equal. We next prove that the same exterior data admit at most one positive inward branch.

\begin{lemma}[Uniqueness of a positive inward branch]\label{lem:branch}
Fix an outward vertex state ${a_*}\in\R$, a current state ${u_*}>0$, and an outward interval state ${b_*}\in\R$.
At most one positive direction parameter satisfies the two transfers and scalar equation and produces a positive inward vertex state.
\end{lemma}
For continuation from left to right, these quantities are $a_*=\mathsf u_{i-1}$, $u_*=\mathsf u_i$, $b_*=\mathsf v_{i-1}$, and the inward state is $x_*=\mathsf u_{i+1}$.

\begin{proof}
For a parameter $d>0$, the transfers first determine
$v_{\rm in}=rd{u_*}-{b_*}$, ${x_*}={a_*}-(r-1)d{b_*}+3\alpha d^2{u_*}$.
The scalar equation also requires ${a_*}+{x_*}-\alpha d^2{u_*}=-c-\omega d$.
Comparing the two expressions for ${x_*}$ gives
\begin{equation}\label{eq:branchquadratic}
 2\alpha {u_*}d^2+[\omega-(r-1){b_*}]d+2{a_*}+c=0.
\end{equation}
If distinct positive roots $d,e$ produce positive inward states ${x_*},{z_*}$, their product gives
$2\alpha {u_*}de=2{a_*}+c$. The scalar equation also gives
${x_*}=\alpha d^2{u_*}-{a_*}-c-\omega d$ and its counterpart for ${z_*}$. Hence
\[
 \frac {x_*}d+\frac {z_*}e
 =\left(\alpha {u_*}-\frac{{a_*}+c}{de}\right)(d+e)-2\omega
 =-\frac c2(1/d+1/e)-2\omega<0,
\]
This contradicts ${x_*},{z_*}>0$. Once the parameter is fixed, the two transfers determine both inward states uniquely.
\end{proof}

\begin{proposition}[Reflection forced by actual extension]\label{prop:reflection}
An actual nonzero even-order corner satisfies
\begin{equation}\label{eq:reflection}
 d_i=d_{n-i},\qquad \mathsf u_i=\mathsf u_{n-i},\qquad
 \mathsf v_i=\mathsf v_{n-1-i}.
\end{equation}

\end{proposition}
\begin{proof}
We know $\mathsf u_0=\mathsf u_n$. The exterior relations give equal outward vertex states,
and $\mathsf v_{-1}=\mathsf v_n$. Reading the second transfer backward from the right endpoint
reverses both differences, so the same inward-branch equation applies.
Both inward states are positive. Branch uniqueness therefore gives
$d_0=d_n,\mathsf u_1=\mathsf u_{n-1},\mathsf v_0=\mathsf v_{n-1}$.
Repeat until reaching the center to obtain \eqref{eq:reflection}.

\end{proof}

For the later moment and path calculations, retain the original direction data and put $a_i=|\beta_i|>0$. The product normalization gives
\begin{equation}\label{eq:ai}
 \begin{gathered}
 h_i=a_ia_{i+1},\qquad
 d_i=a_i/a_{i-1}+a_i/a_{i+1}\quad(1\le i<n),\\
 d_0=a_0/a_1-a_0^2,\qquad d_n=a_n/a_{n-1}-a_n^2.
 \end{gathered}
\end{equation}

\begin{figure}[tb]
\centering
\begin{tikzpicture}[x=1.1cm,y=1cm,every node/.style={font=\small}]
\draw(-.25,0)--(6.25,0);
\foreach \i in {0,...,6}
 {\fill(\i,0)circle(1.4pt);\node[below]at(\i,0){$\i$};}
\draw[dashed](3,-.45)--(3,1.2);\node[above]at(3,1.2){$i=3$};
\draw[<->](0,.45)--(6,.45);
\draw[<->](1,.7)--(5,.7);
\draw[<->](2,.95)--(4,.95);
\node at(3,-.9){$d_i=d_{n-i},\quad\mathsf u_i=\mathsf u_{n-i}$};
\node at(3,-1.45){$\mathsf v_i=\mathsf v_{n-1-i}$};
\end{tikzpicture}
\caption{Scalar reflection for $n=6$: node states and parameters at indices $i$ and $n-i$ agree; interval states pair $i$ with $n-1-i$. Positions show indices only.}\label{human:reflection-figure}
\end{figure}
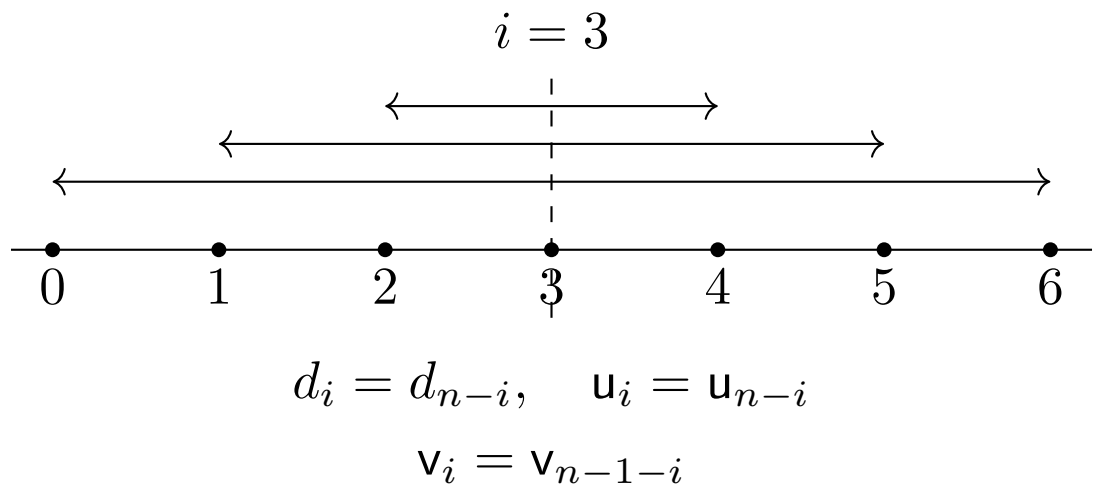

\subsection{From direction separation to a balanced even-path contradiction}\label{sec:paths}
We will prove that the entire even-index group of direction parameters lies strictly on one side of the odd-index group.
Alternation of adjacent comparisons alone is insufficient; we also need the locations of the two subsequences' extrema.

\subsubsection{Reconstructing one interval state from both endpoints}
The two transfers and scalar equation give at the left endpoint
\[
 \frac{r-1}{2}\mathsf v_i
 =\alpha d_i\mathsf u_i+\frac{\mathsf u_{i+1}}{d_i}
 +\frac c{2d_i}+\frac\omega2.
\]
At the right endpoint replace $(d_i,\mathsf u_i,\mathsf u_{i+1})$ by
$(d_{i+1},\mathsf u_{i+1},\mathsf u_i)$.
Put $Y_i=d_i\mathsf u_i$, subtract, and multiply by $d_id_{i+1}$ to obtain
\begin{equation}\label{eq:adjacent-residual}
 (\alpha d_id_{i+1}-1)(Y_{i+1}-Y_i)=\frac c2(d_{i+1}-d_i).
\end{equation}
The coefficient is positive. With $a=d_i,b=d_{i+1},{u_*}=\mathsf u_i,{w_*}=\mathsf u_{i+1}$,
the equation is an affine relation with fixed point $k=c/[2(\alpha ab-1)]$.
We first determine on which side of this fixed point the actual states lie.

\begin{lemma}[A strict edge threshold]\label{lem:threshold}
Put $k_i=c/[2(\alpha d_id_{i+1}-1)]$. Then $\mathsf u_i,\mathsf u_{i+1}>k_i$, and
\begin{equation}\label{eq:reverse-slope}
 Y_{i+1}-Y_i=-\Lambda_i(\mathsf u_{i+1}-\mathsf u_i),\qquad
 \Lambda_i=\frac{k_id_{i+1}}{\mathsf u_i-k_i}>0.
\end{equation}
\end{lemma}
\begin{proof}
Interchange endpoints if necessary to arrange $0<a\le b$. Equation~\eqref{eq:adjacent-residual} is equivalent to
\begin{equation}\label{eq:edgeconvex}
 {w_*}=\frac ab{u_*}+\left(1-\frac ab\right)k.
\end{equation}
If ${u_*}\le k$, then ${u_*}\le {w_*}\le k$.
Let ${a_*}$ be the other neighbor state at the node with parameter $a$. Positivity gives ${a_*}>-\chi$, including an exterior neighbor.
The scalar equation yields $\alpha a^2{u_*}={a_*}+{w_*}+c+\omega a>{u_*}+c-\chi+\omega a$,
so $\alpha a^2>1$. In particular, $\alpha a^2-a/b>0$.
Insert \eqref{eq:edgeconvex} and use ${u_*}\le k$:
\[
 \begin{aligned}
 c-\chi+\omega a&<\alpha a^2{u_*}-{w_*}\\
 &=(\alpha a^2-a/b){u_*}-(1-a/b)k\\
 &\le(\alpha a^2-1)k\le c/2.
 \end{aligned}
\]
The last inequality uses $a\le b$, namely
$(\alpha a^2-1)/(\alpha ab-1)\le1$.
Consequently $a<(\chi-c/2)/\omega=1/(2r)$, which instead gives
$\alpha a^2<(r-1)/(24r)<1$, a contradiction.
Hence ${u_*}>k$, and the convex combination gives ${w_*}>k$.
Finally $b({w_*}-k)=a({u_*}-k)$ implies
$b{w_*}-a{u_*}=(b-a)k=-kb({w_*}-{u_*})/({u_*}-k)$, proving \eqref{eq:reverse-slope}.
When $a=b$, both differences vanish and the identity still holds.
\end{proof}

\subsubsection{Signs and magnitudes of adjacent differences and full parity separation}
\begin{proposition}[Global parity separation]\label{prop:parity}
Let $s=\sgn(\mathsf R-1)$. Then
\begin{equation}\label{eq:parity}
 s(d_{2j}-d_{2k+1})>0
 \qquad(2j,2k+1\in[0,n]).
\end{equation}
\end{proposition}
\begin{proof}
Put $x_i=\mathsf u_{i+1}-\mathsf u_i$, $-1\le i\le n$.
Subtracting the first transfer at adjacent nodes gives
$\mathsf v_{i+1}-\mathsf v_{i-1}=r(Y_{i+1}-Y_i)$;
applying the second transfer to the two adjacent intervals also gives
\[
 \mathsf v_{i+1}-\mathsf v_{i-1}
 =\frac2{r-1}\left(\frac{x_{i+1}+x_i}{d_{i+1}}+
                         \frac{x_i+x_{i-1}}{d_i}\right).
\]
Compare these expressions using \eqref{eq:reverse-slope}:
\[
 \frac{x_{i-1}}{d_i}+
 \left(\frac1{d_i}+\frac1{d_{i+1}}+3\alpha\Lambda_i\right)x_i
 +\frac{x_{i+1}}{d_{i+1}}=0\quad(0\le i<n).
\]
The exterior conditions give $x_{-1}=(1-\mathsf R)(\mathsf u_0+\chi)$ and
$x_n=(\mathsf R-1)(\mathsf u_0+\chi)$.
To make both endpoint values positive and both off-diagonal coefficients negative, define $Z_i=s(-1)^ix_i$.
Since $n$ is even,
\[
 Z_{-1}=Z_n=B:=|\mathsf R-1|(\mathsf u_0+\chi)>0,
\]
the internal equation becomes
\begin{equation}\label{eq:Z-equation}
 -\frac{Z_{i-1}}{d_i}+
 \left(\frac1{d_i}+\frac1{d_{i+1}}+3\alpha\Lambda_i\right)Z_i
 -\frac{Z_{i+1}}{d_{i+1}}=0.
\end{equation}
A negative internal minimum contradicts the extra positive diagonal term $3\alpha\Lambda_i$.
At a zero minimum, the neighboring values also vanish and propagate to a positive endpoint, again a contradiction.
At a positive internal maximum at least $B$, the same extra diagonal term makes the left side strictly positive.
Thus $0<Z_i<B$.

We also need the magnitudes of these differences to be ordered. Write $n=2m$.
State reflection gives $Z_i=Z_{n-1-i}$. Rewrite \eqref{eq:Z-equation} as a difference of quotients,
and put $J_i=(Z_i-Z_{i-1})/d_i$. Then
\begin{equation}\label{eq:J-growth}
 J_{i+1}-J_i=3\alpha\Lambda_iZ_i>0.
\end{equation}
The central equality $Z_m=Z_{m-1}$ gives $J_m=0$. The left difference quotients are therefore strictly negative and the right ones strictly positive:
\begin{equation}\label{eq:Z-nesting}
 B>Z_0>Z_1>\cdots>Z_{m-1}=Z_m<\cdots<Z_{n-1}<B.
\end{equation}
Since $x_i=s(-1)^iZ_i$, adjacent signs alternate, and
\begin{equation}\label{eq:twostep-u}
 \sgn(\mathsf u_{i+2}-\mathsf u_i)=s(-1)^i
 \quad(-1\le i\le m-2).
\end{equation}
The central two-step difference is zero.

We now transfer this ordering to the parameters $d_i$. By \eqref{eq:scalareq}, $d_i$ is the unique positive root of
\begin{equation}\label{eq:rootmonotone}
 \alpha\mathsf u_i t^2-\omega t-(\mathsf u_{i-1}+\mathsf u_{i+1}+c)=0
\end{equation}
The neighbor sum plus $c$ is positive, including at an endpoint: one neighbor is positive, the exterior state exceeds $-\chi$,
and $c-\chi>0$. Write this root as $t=t({\xi},{\eta})$, where $\xi=\mathsf u_i>0$ and $\eta=\mathsf u_{i-1}+\mathsf u_{i+1}$. Its derivatives are
\[
 \frac{\partial t}{\partial {\xi}}=-\frac{\alpha t^2}{2\alpha {\xi}t-\omega}<0,
 \qquad \frac{\partial t}{\partial {\eta}}=\frac1{2\alpha {\xi}t-\omega}>0,
\]
because $2\alpha {\xi}t-\omega=\alpha {\xi}t+({\eta}+c)/t>0$.
Let $A_i=\mathsf u_{i-1}+\mathsf u_{i+1}$. For $0\le i\le m-2$,
\[
 \begin{aligned}
 \mathsf u_{i+2}-\mathsf u_i&=s(-1)^i(Z_i-Z_{i+1}),\\
 A_{i+2}-A_i&=-s(-1)^i[(Z_{i-1}-Z_i)+(Z_{i+1}-Z_{i+2})].
 \end{aligned}
\]
The two parameters change strictly in opposite directions: the first bracket in the second line is always positive, and only the second may vanish in the last comparison.
The two derivatives therefore move the positive root in the same direction, giving
$\sgn(d_{i+2}-d_i)=-s(-1)^i$ for $0\le i\le m-2$.
Also \eqref{eq:edgeconvex} and $\mathsf u_i>k_i$ give
$\sgn(d_{i+1}-d_i)=-s(-1)^i$ for $0\le i<n$.

We now pass from local strict comparisons to separation of the whole groups. First take $s=1$.
On the left half, the even subsequence strictly decreases and the odd subsequence strictly increases; reflection $d_i=d_{n-i}$ determines the right half.
The parity of $m$ decides which group contains the center, so we identify the two group extrema separately.

If $m=2k$ is even, the left-half orders are
\[
 d_{2j+2}<d_{2j}\quad(0\le j<k),\qquad
 d_{2j+3}>d_{2j+1}\quad(0\le j<k-1).
\]
For $k=1$ the second range is empty: the odd group on the left contains only $d_1$.
After reflection, the minimum over all even parameters is $d_m$, and the maximum over all odd parameters is $d_{m-1}=d_{m+1}$.
Use the central adjacent comparison at $i=m-1$. This index is odd, giving $d_m-d_{m-1}>0$.
Hence
\[
 \min_{0\le2j\le n}d_{2j}=d_m
 >d_{m-1}=d_{m+1}
 =\max_{0\le2j+1\le n}d_{2j+1}.
\]

If $m=2k+1$ is odd, the left-half orders instead are
\[
 d_{2j+2}<d_{2j}\quad(0\le j<k),\qquad
 d_{2j+3}>d_{2j+1}\quad(0\le j<k).
\]
The minimum over all even parameters is now $d_{m-1}=d_{m+1}$, and the maximum over all odd parameters is $d_m$.
Again compare at $i=m-1$. This index is even, so $d_m-d_{m-1}<0$.
Thus
\[
 \min_{0\le2j\le n}d_{2j}=d_{m-1}=d_{m+1}
 >d_m=\max_{0\le2j+1\le n}d_{2j+1}.
\]
In both cases every even parameter strictly exceeds every odd parameter, giving the full group comparison in addition to adjacent alternation.

For $s=-1$, all the adjacent and two-step comparisons above reverse.
If $m$ is even, $d_m$ is the even-group maximum and $d_{m-1}=d_{m+1}$ is the odd-group minimum; the central comparison gives $d_m<d_{m-1}$.
If $m$ is odd, the even-group maximum is $d_{m-1}=d_{m+1}$ and the odd-group minimum is $d_m$; the central comparison gives $d_{m-1}<d_m$.
Every even parameter is therefore strictly smaller than every odd parameter in this case.
This proves both sign cases of \eqref{eq:parity}.
\end{proof}

\begin{remark}
For $s=1$, the two central arrangements take the following forms at the two shortest lengths:
\[
 \begin{array}{ll}
 n=4:&d_1=d_3<d_2<d_0=d_4,\\[2pt]
 n=6:&d_1=d_5<d_3<d_2=d_4<d_0=d_6.
 \end{array}
\]
The center belongs to the even group in the first row and to the odd group in the second. Both rows place every odd parameter before every even parameter.
\end{remark}

\subsubsection{A balanced even path from the same actual data}
Direction separation is now proved. Return to the first-round equation and obtain both moments directly from the complete cycle. All three exterior intervals are included, so
\[
 \sum_i\delta_i=0,\qquad
 \sum_it_i\delta_i=\sum_i(t_i-t_{i+1})\kappa_i
 =-\sum_ih_i\kappa_i=-\sum_i(\zeta_{i+1}-\zeta_i)=0.
\]
Multiply \eqref{eq:hermite1} by $1,t_i$ and sum. In the first sum each interval contributes
$(r+1)(\beta_i+\beta_{i+1})/h_i$. In the second its numerator is
\[
 \begin{aligned}
 &t_i(r\beta_i+\beta_{i+1})+t_{i+1}(\beta_i+r\beta_{i+1})\\
 &\quad=(r+1)(t_i\beta_i+t_{i+1}\beta_{i+1})
       +h_i(\beta_i-\beta_{i+1}).
 \end{aligned}
\]
After division by $h_i$, the last terms cancel around the cycle. The exterior data are
\[
 (h_n,h_{n+1},h_{n+2})=(-1,1,-1),\qquad
 (t_n,t_{n+1},t_{n+2},t_0)=(1,0,1,0).
\]
Their contributions to the two sums are $-\beta_0-\beta_n$ and $-\beta_n$, respectively. Hence
\begin{equation}\label{human:weighted-flux}
 \sum_{i=0}^{n-1}\frac{\beta_i+\beta_{i+1}}{h_i}=\beta_0+\beta_n,
 \qquad
 \sum_{i=0}^{n-1}\frac{t_i\beta_i+t_{i+1}\beta_{i+1}}{h_i}=\beta_n.
\end{equation}

Using the original normalized direction data, put $a_i=|\beta_i|$, $z_i=1/a_i$, and
\[
 b_0=1/h_0-1,\quad b_n=1/h_{n-1}-1,\quad
 b_i=1/h_{i-1}+1/h_i\ (0<i<n).
\]
Equation~\eqref{eq:ai} gives $d_i=b_i a_i^2$. Since
$\beta_i=\sgn(\beta_0)(-1)^ia_i$, collecting the two cyclic sum identities by node gives
\begin{equation}\label{eq:moment0}
 M_0:=\sum_{i=0}^n(-1)^id_i z_i=0,
\end{equation}
\begin{equation}\label{eq:moment1}
 M_1:=\sum_{i=0}^n(-1)^id_i z_i t_i=0.
\end{equation}
Thus both zero moments follow from the actual first-round equation without assuming reflection of the directions.

The reciprocal sequence already satisfies $d_iz_i=z_{i-1}+z_{i+1}$ internally, but has $n+1$ positive values.
Multiplying by an affine weight that vanishes at the right endpoint leaves exactly the even number $n$ of positive path values. Define
\begin{equation}\label{eq:actualpath}
 g_i=(1-t_i)z_i\quad(0\le i\le n),\qquad g_{-1}=0.
\end{equation}
Then $g_0,\ldots,g_{n-1}>0$ and $g_n=0$. The internal product difference is
\[
 \begin{aligned}
 g_{i-1}+g_{i+1}
 &=(1-t_i)(z_{i-1}+z_{i+1})+h_{i-1}z_{i-1}-h_iz_{i+1}\\
 &=d_ig_i,
 \end{aligned}
\]
The last two terms cancel because each is $a_i$.
At the first endpoint, $d_0g_0=1/a_1-a_0=(1-h_0)/a_1=g_1$.
The neighbor-sum ratios of this $n$-vertex path are therefore precisely $d_0,\ldots,d_{n-1}$.

The zero moments also give $\sum_{i=0}^{n-1}(-1)^id_ig_i=M_0-M_1=0$.
Let $A_n$ be the path adjacency matrix and $S=\diag(1,-1,\ldots,1,-1)$.
Since $SA_n=-A_nS$,
\[
 \sum_i(-1)^id_ig_i=\mathbf1^{\mathsf T}SA_ng
 =-\mathbf1^{\mathsf T}A_nSg=-\boldsymbol\nu^{\mathsf T}Sg,
 \qquad\boldsymbol\nu=A_n\mathbf1.
\]
this is balance with degree weight one at the endpoints and two internally.
The last geometric node $n$ has become an auxiliary zero endpoint; the positive path vertices are only $0,\ldots,n-1$.

\subsubsection{Overlap on a balanced even path and its nonnegative certificate}
\begin{theorem}[Overlap on a balanced even path]\label{thm:path}
Let $N\ge2$ be even, $g_0,\ldots,g_{N-1}>0$, and $g_{-1}=g_N=0$. Put
\[
 \mu_i=(g_{i-1}+g_{i+1})/g_i,\qquad
 \nu_i=\begin{cases}1,&i=0,N-1,\\2,&1\le i\le N-2.\end{cases}
\]
If
\begin{equation}\label{eq:pathbalance}
 \sum_{i=0}^{N-1}(-1)^i\nu_i g_i=0,
\end{equation}
then the two closed intervals
\begin{equation}\label{eq:pathoverlap}
 [\min_{i\ \text{even}}\mu_i,\max_{i\ \text{even}}\mu_i]
 \quad\text{and}\quad
 [\min_{i\ \text{odd}}\mu_i,\max_{i\ \text{odd}}\mu_i]
\end{equation}
intersect.
\end{theorem}

The required certificate can now be identified. If every even ratio exceeds every odd ratio, choose a threshold $\tau>0$ strictly between them.
Then $b=S(A_N-\tau I)g>0$. For a nonnegative nonzero $q$, we have $q^{\mathsf T}b>0$;
on the other hand, anticommutation gives
$q^{\mathsf T}b=-((A_N+\tau I)q)^{\mathsf T}Sg$.
Balance removes the degree-vector term. Thus a vector $q$ with the next property will contradict separation.

\begin{lemma}[A nonnegative dual for every threshold]\label{lem:pathdual}
Let $A_N$ be the adjacency matrix of an even $N$-vertex path, and put
\[
 S=\diag(1,-1,\ldots,1,-1),\qquad\boldsymbol\nu=A_N\mathbf1.
\]
For every $\tau>0$ there is a nonzero $q\ge0$ such that
\begin{equation}\label{eq:pathdual}
 S((A_N+\tau I)q-\boldsymbol\nu)\ge0
\end{equation}
entrywise.
\end{lemma}
\begin{proof}
For $0<\tau\le1$, take $q=(0,1,\ldots,0,1)$. Even rows equal their degrees and odd rows are at most their degrees.
For $\tau\ge2$, take $q=(1,0,\ldots,1,0)$. Odd rows equal their degrees and even rows are at least their degrees.
Both checks include the degree-one endpoints.

In the remaining range write $\tau=2\cos\theta$, $0<\theta<\pi/3$.
First solve $(A_L+\tau I)q=\boldsymbol\nu_L$ on an even block of length $L$.
Temporarily number its vertices $j=1,\ldots,L$. Introduce auxiliary boundary values equal to one at $j=0,L+1$. The recurrence then has right-hand side two at every node, including the endpoints.
The recurrence has constant particular solution $2/(2+\tau)$ and homogeneous roots $-e^{\pm i\theta}$.
With center $H=(L+1)/2$, the reflection-symmetric homogeneous part has the form
$B(-1)^j\sin((H-j)\theta)$.
These auxiliary boundary values determine $B=\tau/[(2+\tau)\sin(H\theta)]$.
Thus in the original zero-based vector indexing,
\begin{equation}\label{eq:sineblock}
 q^{(L)}_{j-1}=\frac2{2+\tau}
 +\frac\tau{2+\tau}(-1)^j\frac{\sin((H-j)\theta)}{\sin(H\theta)}
 \quad(1\le j\le L).
\end{equation}
The auxiliary boundary values for the certificate are one; the boundary values appended to the positive path $g$ are zero.

For nonnegativity choose even $L$ with $(L-1)\theta\le\pi$.
Then $0<H\theta\le\pi/2+\theta<\pi$, so the denominator is positive.
The sine addition identity verifies internal right side two; removing these two unit boundary terms gives right-hand side one at the endpoints:
\begin{equation}\label{eq:sineblock-system}
 (A_L+\tau I)q^{(L)}=\boldsymbol\nu_L.
\end{equation}
Since $|(H-j)\theta|\le(H-1)\theta\le\pi/2$, the oscillatory term has largest absolute value at an endpoint.
Each component is therefore at least the two equal endpoint values
\begin{equation}\label{eq:sine-end}
 q^{(L)}_0=q^{(L)}_{L-1}
 =\frac{2\sin\theta\cos((H-1)\theta)}{(2+\tau)\sin(H\theta)}\ge0.
\end{equation}
Endpoints may vanish, while the positive right side in \eqref{eq:sineblock-system} ensures that the vector is nonzero.

Let $L_*$ be the largest admissible even length. If $N\le L_*$, take $L=N$.
If $N>L_*$, place an $L_*$ block on the right and prepend $(N-L_*)/2$ pairs $(0,1)$.
Maximality gives $(L_*+1)\theta>\pi$, so this block has $H\theta>\pi/2$.
Its endpoint formula can also be written
\[
 q^{(L_*)}_0=\tan(\theta/2)[\sin\theta+\cos\theta\cot(H\theta)],
\]
Thus $0\le q^{(L_*)}_0\le1-\cos\theta=(2-\tau)/2$.
Only two rows change at the interface. The last prepended odd row has value
$\tau+q^{(L_*)}_0\le(\tau+2)/2<2$, satisfying the odd-row upper bound;
the first even row of the tail has its old value one increased by its new left neighbor one to two, exactly its new internal degree.
Every other row retains its condition. The prepended length is even, so the tail keeps its parity.
This constructs a certificate for every finite even path length.
\end{proof}

\begin{proof}[Proof of Theorem~\ref{thm:path}]
If the intervals are disjoint, first suppose the even ratios are larger. Choose the threshold $\tau$ and certificate $q$, and put
$e=S((A_N+\tau I)q-\boldsymbol\nu)\ge0$. Then
\[
 0<q^{\mathsf T}b
 =-((A_N+\tau I)q)^{\mathsf T}Sg
 =-\boldsymbol\nu^{\mathsf T}Sg-e^{\mathsf T}g
 =-e^{\mathsf T}g\le0,
\]
This is a contradiction. If the odd ratios are larger, reverse the path. Even $N$ exchanges parity while preserving zero balance, giving the same contradiction.
\end{proof}

\begin{remark}\label{rem:pathscope}
Both length and degree balance enter the proof. On the odd path $g=(1,1,1)$, balance is zero,
but the ratios $(1,2,1)$ have separated parity intervals. On the even path $g=(1,2)$, balance is nonzero and the ratios are $(2,1/2)$.
\end{remark}

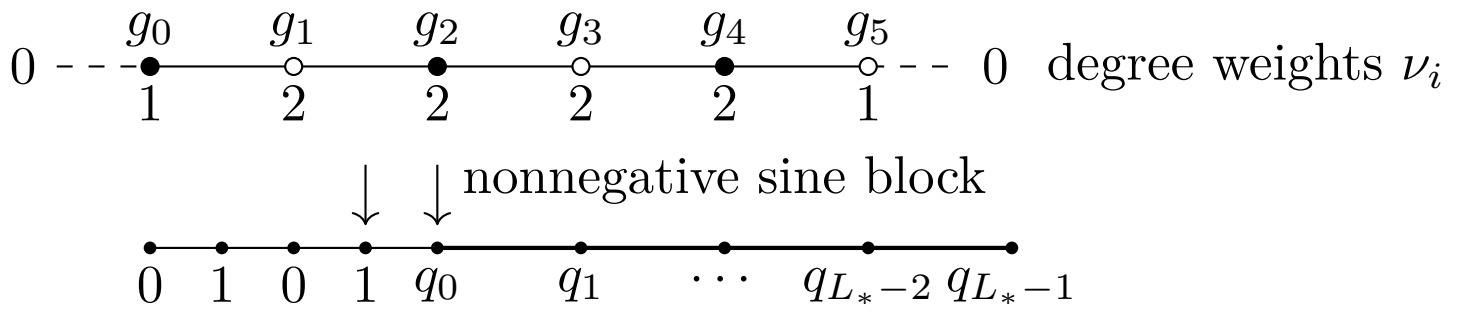
\begin{figure}[tb]
\centering
\begin{tikzpicture}[x=.95cm,y=1cm,every node/.style={font=\small}]
\draw[dashed](-.65,0)--(0,0);\draw(0,0)--(5,0);\draw[dashed](5,0)--(5.65,0);
\foreach \i in {0,...,5}{\draw[fill=white](\i,0)circle(1.6pt);\node[above]at(\i,0){$g_{\i}$};}
\foreach \i in {0,2,4}{\fill(\i,0)circle(1.6pt);}
\node[left]at(-.65,0){$0$};\node[right]at(5.65,0){$0$};
\node[below]at(0,0){$1$};\node[below]at(5,0){$1$};
\foreach \i in {1,...,4}{\node[below]at(\i,0){$2$};}
\node[right]at(6.1,0){degree weights $\nu_i$};
\draw(0,-1.2)--(2,-1.2);\draw[thick](2,-1.2)--(6,-1.2);
\foreach \x/\lab in {0/0,.5/1,1/0,1.5/1,2/{q_0},3/{q_1},4/{\cdots},5/{q_{L_*-2}},6/{q_{L_*-1}}}
 {\fill(\x,-1.2)circle(1.2pt);\node[below]at(\x,-1.2){$\lab$};}
\draw[->](1.5,-.65)--(1.5,-1.05);\draw[->](2,-.65)--(2,-1.05);
\node at(4,-.75){nonnegative sine block};
\end{tikzpicture}
\caption{Top: auxiliary zero endpoints and degree weights of the actual positive path. Bottom: the certificate joins an even checkerboard prefix to a finite sine block;
only the two rows indicated by arrows need an interface check.}\label{fig:path}
\end{figure}

\subsubsection{Completion of the leading-layer theorem}
\begin{theorem}[Vanishing of the normalized leading layer]\label{thm:leading}
For every integer $r\ge3$ and every dual defect class, each $H_v$ may be replaced by an equivalent representative modulo $\A_v$ so that
\begin{equation}\label{eq:leading-vanish}
 H_v-K_\tau\in P_v^{r+2}A_{2r-1}
\end{equation}
at every incident vertex and face. Auxiliary boundary representatives remain zero.
\end{theorem}
\begin{proof}
For odd $r$, the assertion is Lemma~\ref{lem:odd-leading}. For even $r\ge4$, first normalize the two-direction vertices. If support remains, finiteness gives an extreme vertex of the support. Proposition~\ref{prop:oddcorner} and Lemma~\ref{lem:short-corner} leave only even $n\ge4$.
Positivity, branch uniqueness, and reflection of the states from the actual second radial layer give the full strict separation of Proposition~\ref{prop:parity}.
The same actual corner also gives, by \eqref{eq:actualpath}, a balanced positive path with $n$ vertices
and ratios $d_0,\ldots,d_{n-1}$. Removing the last geometric even node $n$ retains strict separation,
contradicting the interval overlap in Theorem~\ref{thm:path}.
Thus the normalized support is empty. Every coefficient was defined from the actual residual divided by the nonzero common scale $\Lambda$,
so its zero leading term is equivalent to actual leading-layer vanishing. Changing the radial unit also only multiplies $P_v$ by a nonzero constant.
This proves \eqref{eq:leading-vanish}. The only representative changes were the initial local annihilator corrections at two-direction vertices.
\end{proof}

\subsubsection{Completion of the dimension theorem}
\begin{proof}[Completion of Theorem~\ref{thm:main}]
The odd orders have been proved in Theorem~\ref{thm:oddcase}. For any even $r\ge4$ and any dual defect class, Theorem~\ref{thm:leading} now eliminates the leading coefficient. Theorem~\ref{thm:from-leading} then makes that class zero. Hence $M_{3r+1}=0$ also for even $r$. Proposition~\ref{prop:dimension}, Lemma~\ref{lem:localdim}, and \eqref{eq:binomialidentity} yield \eqref{eq:main}, completing every $r\ge3$.
\end{proof}

\begin{corollary}[Dimension stability at fixed local direction counts]\label{cor:constancy}
Fix the combinatorics of a disk triangulation. For each $r\ge3$, $\dim S_{3r+1}^r(\Delta)$ is constant on a family of nondegenerate straight-line realizations in which every $t_v$ is fixed.
\end{corollary}
\begin{proof}
Every term on the right side of \eqref{eq:main} is constant in that family.
\end{proof}

\section{Coefficient expansions for the leading layer}\label{app:coefficients}

This appendix gives the full expansions used in Appendices \ref{sec:corner} and \ref{sec:positivity}. We use the center form $U=P_u$, the shared neighbor form $A=L_i$, and the preceding and following forms $W=L_{i-1}$ and $T=L_{i+1}$. The residuals are divided by the common nonzero scalar $\Lambda$; $H_u$ and $K_i$ remain fixed.

\subsection{The two axis restrictions and the Hermite identities}
At the corner node $t_i$, let $x=t-t_i$, and abbreviate
\[
 h=h_{i-1},\quad k=h_i,\quad z=\zeta_i,\quad
 a=\kappa_{i-1},\quad b=\kappa_i,\quad \delta=b-a.
\]
Under the substitution used to define $f_i$, the three neighboring point forms read
\begin{align*}
 A&=-x+zU,\\
 W&=-h-x+(z-ah)U,\\
 T&=k-x+(z+bk)U.
\end{align*}
Thus $T=sW+pU+qA$ with
\begin{equation}\label{eq:append-spq}
 s=-k/h,\qquad p=k\delta,\qquad q=1+k/h.
\end{equation}
Here $p,q$ denote the coefficients of $U,A$, respectively. In the node notation of Section~\ref{body:even}, $p=\xi_i$, $q=q_i$, and $s=s_i$.

Let $C=\sigma_{i-1}/h^r$. Write the face-potential difference as
\[
 \frac{K_i-K_{i-1}}{\Lambda}=U^{r+1}A^{r+1}Q(A,U,W).
\]
Its two axis restrictions are
\begin{align}
 Q(A,0,W)&=C\,\frac{W^r-(W+(q/s)A)^r}{A}=:Q_0(A,W),\label{eq:axis1}\\
 Q(0,U,W)&=-Cq\sum_{j=1}^{r}\binom rj\frac{p^{j-1}}{s^j}
                           U^{j-1}W^{r-j}=:Q_U(U,W).\label{eq:axis2}
\end{align}
The second restriction follows from the neighboring vertex's leading term and \eqref{eq:lambda-coupling}: the neighboring coefficient multiplied by $p$ equals $Cq$, and the binomial numerator has a factor $p$ in every term. Cancelling that factor polynomially gives the displayed $p^{j-1}$ expression, valid also at $p=0$. The two restrictions agree at $A=U=0$. Every full remainder therefore has the form
\begin{equation}\label{eq:UV-remainder}
 Q=Q_0+Q_U-Q_0(0,W)+AU\,T_0,
 \qquad T_0\in A_{r-3},
\end{equation}
with the last term zero if its degree is negative. Indeed, a polynomial with zero restrictions at both axes is divisible by their product $AU$.

First record the orders at which the undetermined remainder can contribute. After removing the central factor $U^{r+1}$ from the face difference, the unknown term from \eqref{eq:UV-remainder} is $A^{r+2}UT_0$. Its coefficient of $U$ contains $A^{r+2}$, which starts at $x^{r+2}$ after substituting $A=-x$. In the second comparison, extraction of $U^2$ followed by division by $A^{r-2}$ gives
\[
 \frac{[U^2](A^{r+2}UT_0)}{A^{r-2}}=A^4[U]T_0.
\]
Thus the remainder affects neither $x^{r-1},x^r,x^{r+1}$ in the first comparison nor $A,A^2,A^3$ in the second. Coefficients of $U$ are extracted with $A,W$ fixed as independent variables.

Subtracting the two leading terms from the face difference and extracting the coefficient of $U$ gives
\begin{align}
 [x^{r-1}](f_i-f_{i-1})&=0,\label{eq:fdiff1}\\
 [x^r](f_i-f_{i-1})&=-r\sigma_{i-1}\delta_i,\label{eq:fdiff2}\\
 [x^{r+1}](f_i-f_{i-1})
 &=\binom r2\sigma_{i-1}\delta_i
           \left(\frac1{h_i}-\frac1{h_{i-1}}\right).\label{eq:fdiff3}
\end{align}
To write out this subtraction, define the following difference. Its three terms are the old leading term, the new leading term, and the known contribution from the axis restrictions:
\[
 \mathcal B=CA^rW^r-Cs^{-r}A^r(sW+pU+qA)^r
       -A^{r+1}\bigl(Q_0+Q_U-Q_0(0,W)\bigr).
\]
At $U=0$, $\mathcal B$ vanishes identically in $A,W$ by \eqref{eq:axis1}. Consequently, the changes of $A,W$ with $U$ contribute zero to its first radial coefficient, and
\[
 \left.\pd_U\mathcal B\right|_{U=0}
 =-Cr\frac psA^r\bigl(W+(q/s)A\bigr)^{r-1}
   +Cq\binom r2\frac p{s^2}A^{r+1}W^{r-2}.
\]
Substitute $A=-x$, $W=-h-x$, and $p/s=-h\delta$ to obtain
\begin{align*}
 f_i-f_{i-1}\equiv{}&-r\sigma_{i-1}\delta x^r(1-x/k)^{r-1}\\
 &-\binom r2\sigma_{i-1}\delta(1/h+1/k)
                 x^{r+1}(1+x/h)^{r-2}\pmod{x^{r+2}}.
\end{align*}
The coefficient of $x^{r+1}$ is
$r(r-1)\sigma_{i-1}\delta/k-\binom r2\sigma_{i-1}\delta(1/h+1/k)$,
which is \eqref{eq:fdiff3}; the other two coefficients give
\eqref{eq:fdiff1}--\eqref{eq:fdiff2}. In particular,
\eqref{eq:fdiff1} identifies the shared endpoint coefficient defining $\beta_i$.

A polynomial of degree at most $2r-1$ with zeros of multiplicity at least $r-1$ at $0,h$ is determined by its endpoint coefficients $a_{\rm L},a_{\rm R}$ of order $r-1$. With $\epsilon=(-1)^{r-1}$, Hermite reconstruction gives
\begin{equation}\label{eq:hermite-polynomial}
 f(x)=x^{r-1}(x-h)^{r-1}
 \left[\frac{a_{\rm R}-\epsilon a_{\rm L}}{h^r}x
       +\frac{\epsilon a_{\rm L}}{h^{r-1}}\right].
\end{equation}
For the two pieces at $t_i$, write the shared endpoint coefficient as $A_i=\sigma_{i-1}\beta_i$. Expansion gives
\[
 \begin{aligned}
 \frac{f_i(t_i+x)}{\sigma_{i-1}}
 ={}&\beta_i x^{r-1}-\frac{r\beta_i+\beta_{i+1}}k x^r\\
 &+(r-1)\frac{r\beta_i/2+\beta_{i+1}}{k^2}x^{r+1}+O(x^{r+2}),\\
 \frac{f_{i-1}(t_i+x)}{\sigma_{i-1}}
 ={}&\beta_i x^{r-1}+\frac{\beta_{i-1}+r\beta_i}h x^r\\
 &+(r-1)\frac{\beta_{i-1}+r\beta_i/2}{h^2}x^{r+1}+O(x^{r+2}).
 \end{aligned}
\]
Subtract and compare with \eqref{eq:fdiff2}--\eqref{eq:fdiff3}. Using $\sigma_i/\sigma_{i-1}=(-1)^r$ yields \eqref{eq:hermite1}--\eqref{eq:hermite2}. The second comparison cancels $r-1$, so these identities require $r\ge2$.

\subsection{Coefficients of the second radial layer}
For even $r$, the two incident derivative kernels give \eqref{eq:jet2}. Its second radial layer has three coefficients per face. The lowest-power matching across a node identifies the endpoint coefficients $\gamma_i$, while the middle coefficient $\eta_i$ remains attached to its face interval.

Use independent variables $A=L_i$, $W=L_{i-1}$, and $U$ with $L_{i+1}=sW+pU+qA$. Divide the $U^2$ coefficient by $A^{r-2}$ and by the old nonzero leading scale, and put $W=1$. The old face contributes
\[
 \gamma_i+\eta_{i-1}A+\gamma_{i-1}A^2.
\]
The contributions from the three orders of the new face in \eqref{eq:jet2} are, respectively,
\begin{align}
 &\binom r2\frac{p^2A^2}{s^2}
           \left(1+\frac qsA\right)^{r-2},\label{eq:second-0}\\
 &-p\left[\frac{(r-1)\beta_{i+1}A^2}{s^2}
             \left(1+\frac qsA\right)^{r-2}
       +\frac{r\beta_iA}{s}\left(1+\frac qsA\right)^{r-1}\right],
                                                        \label{eq:second-1}\\
 &\frac{\gamma_{i+1}A^2}{s^2}\left(1+\frac qsA\right)^{r-2}
   +\frac{\eta_iA}{s}\left(1+\frac qsA\right)^{r-1}
   +\gamma_i\left(1+\frac qsA\right)^r.\label{eq:second-2}
\end{align}
Subtracting the axis restriction \eqref{eq:axis2} adds $q\binom r3p^2A^3/s^3$ to the old-minus-new difference. By the order estimate following \eqref{eq:UV-remainder}, the unknown remainder starts at $A^4$. In independent coordinates, radial terms containing $U^3$ or higher do not contribute to the coefficient of $U^2$.

The constant coefficient cancels. The coefficients of $A,A^2,A^3$ in the reduced difference are
\begin{align}
 &s^{-1}(r\beta_i p-r\gamma_iq+s\eta_{i-1}-\eta_i),\label{eq:rawjet1}\\
 &s^{-2}\bigl[s^2\gamma_{i-1}-\gamma_{i+1}
 +(r-1)\beta_{i+1}p+r(r-1)\beta_ipq\nonumber\\
 &\hspace{34mm}-\tbinom r2\gamma_iq^2-(r-1)\eta_iq-\tbinom r2p^2\bigr],
                                                        \label{eq:rawjet2}\\
 &-\frac{q(r-2)}{6s^3}\bigl[6\gamma_{i+1}+3(r-1)\eta_iq
       +r(r-1)\gamma_iq^2\nonumber\\
 &\hspace{23mm}-3r(r-1)\beta_ipq-6(r-1)\beta_{i+1}p
       +2r(r-1)p^2\bigr].\label{eq:rawjet3}
\end{align}
Setting the first expression to zero gives \eqref{eq:eta}; substituting it in the second gives \eqref{eq:gamma}. The last expression gives \eqref{eq:thirdjet} exactly where $q(r-2)\ne0$. At $r=2$ the corresponding monomial has negative residual $W$ degree and the condition is absent; at an exterior node $q=0$ it is also absent. This explains why the even-order proof uses $r\ge4$ and why the exterior nodes are treated separately throughout.

\section{Further degree and low-order comparisons}\label{app:degree}\label{sec:discussion}

\subsection{How the three degree comparisons fit together}
The proof computes the same objects in two ways. Edge data determine the coefficients of each face potential, while the vertex conditions make a neighboring face difference divisible by both endpoint forms. Degree $3r+1$ provides three distinct pieces of information.

First, the exponents of a monomial of degree $3r+1$ cannot all be at most $r$. Thus the three edge traces read every coefficient of the face potential. Second, the two incident derivative kernels make a vertex residual divisible by $P_v^{r+1}$. There is only one coefficient at this lowest power, and the leading-coefficient theorem rules it out. Third, removing that coefficient raises the factor only to $P_v^{r+2}$; the neighboring face difference can still have a quotient of degree $r-3$. The uniform higher-layer induction raises the endpoint powers until their total degree exceeds the degree of the face potentials.

\begin{center}
\begin{tabularx}{\linewidth}{@{}cc>{\centering\arraybackslash}p{.19\linewidth}>{\centering\arraybackslash}p{.15\linewidth}X@{}}
\toprule
$r$ & $d$ & Power after the leading layer & Jump quotient degree & Further layers\\\midrule
3&10&5&0&$p=1$\\
4&13&6&1&$p=1$\\
5&16&7&2&$p=1,2$\\
6&19&8&3&$p=1,2$\\
\bottomrule
\end{tabularx}
\end{center}
Each step of the uniform induction first excludes straight-line transmission. The common modes are then lifted to the local annihilator, so the improved divisibility still describes the same original functional class. At termination, $h=r+1+\lceil r/2\rceil$ and $2h>3r+1$, which forces every face difference to vanish. The uniform layer indices in this table count radial powers in the original residual.

\subsection{\texorpdfstring{The first change at degree $3r$}{The first change at degree 3r}}
At degree $3r$, all three edge traces miss the monomial $U^rV^rW^r$. More precisely, put
\[
 \mathcal K_{r,d}=\{K\in A_d:\pd_U^{r+1}K=\pd_V^{r+1}K=\pd_W^{r+1}K=0\}.
\]
Reading the allowed monomial exponents gives
\begin{equation}\label{eq:discussion-facekernel}
 \mathcal K_{r,3r+1}=0,\qquad
 \mathcal K_{r,3r}=\Span\{U^rV^rW^r\}.
\end{equation}
Each degree-$3r$ face potential therefore has one free scalar. The vertex residual also starts at $V^rU^rW^r$, so changing that face scalar changes the leading residual coefficient. Our proof first constructs a unique face polynomial $K_\tau$, then compares the polynomials on adjacent faces. One degree lower, the face scalars must first be coordinated before new edge and boundary propagation relations can be established.

Equation~\eqref{eq:discussion-facekernel} identifies the exact change in the local construction. Whether Schumaker's lower bound is attained universally at degree $3r$ depends on the global cokernel and requires a further global argument. For $r=2$, universal attainment in degree six is posed as an open problem in \citet{WD26}.

\subsection{Septic identities from the two axis restrictions}\label{sec:septic-axes}
We derive \eqref{body:septic-exact-jump} and \eqref{body:septic-cubic-difference}. Take $r=2$ and use the center form $U=P_u$ with
\[
 A=L_i,\qquad W=L_{i-1},\qquad T=L_{i+1}=sW+qA+pU,
\]
where $h=h_{i-1}$, $k=h_i$, $s=-k/h$, $p=k\delta_i$, and $q=1+k/h$. Thus $p=\xi_i$, $q=q_i$, and $s=s_i$ in the node notation. Here $F_i$ is the remainder after division by the common scale $\Lambda$:
\[
 \frac{H_u-K_i}{\Lambda}=\frac1{h_i^2}U^3L_i^2L_{i+1}^2+U^4F_i.
\]
Put $C=h^{-2}$ and write $(K_i-K_{i-1})/\Lambda=U^3A^3Q(A,U,W)$.
At $r=2$, \eqref{eq:axis1} and \eqref{eq:axis2} give
\[
 \begin{aligned}
 Q(A,0,W)&=-\frac{2Cq}sW-\frac{Cq^2}{s^2}A,\\
 Q(0,U,W)&=-\frac{2Cq}sW-\frac{Cpq}{s^2}U.
 \end{aligned}
\]
Since $Q$ is homogeneous of degree one, these axis restrictions determine every coefficient:
\[
 Q=-\frac{Cq}{s^2}(2sW+qA+pU)
   =-\frac q{k^2}(sW+T).
\]
This proves \eqref{body:septic-exact-jump}. Insert it into the difference of the two residuals and cancel $U^3$:
\[
 U(F_i-F_{i-1})=\frac{A^2}{k^2}
 \bigl[qA(sW+T)-(T^2-s^2W^2)\bigr].
\]
Since $T-sW=qA+pU$, the bracket is $-pU(sW+T)$. Cancel $U$ and use $p=k\delta_i$, $s=-k/h$ to obtain
\[
 F_i-F_{i-1}=-\frac p{k^2}A^2(sW+T)
 =-\delta_iA^2\left(\frac Tk-\frac Wh\right),
\]
which is \eqref{body:septic-cubic-difference}. There is no undetermined remainder in this calculation.

\begingroup\small\raggedright

\endgroup
\end{document}